\documentclass[12pt,a4paper]{article}
\usepackage{bbm}

\usepackage[colorlinks,linkcolor=blue,anchorcolor=green,citecolor=red]{hyperref}
\usepackage{lipsum}
\usepackage{amsfonts,amsmath,amssymb,amsthm}
\usepackage{graphicx}
\usepackage{epstopdf}
\usepackage{algorithmic}
\usepackage{enumerate}
\usepackage{mathrsfs}
\usepackage{enumitem}

\usepackage{bm}
\usepackage{stmaryrd}
\SetSymbolFont{stmry}{bold}{U}{stmry}{m}{n} 
\usepackage{todonotes}

\usepackage{tikz}
\usetikzlibrary{intersections, arrows.meta}
\usetikzlibrary{matrix}

\newcommand{\mr}[1]{\mathrm{#1}}

\newcommand{\bn}{\bm{n}}
\newcommand{\bx}{\bm{x}}
\newcommand{\bxb}{\bar{\bm{x}}}
\newcommand{\bxi}{\bm{\xi}}

\newcommand{\e}{\mathrm{e}}
\newcommand{\s}{\mathrm{s}}
\newcommand{\ga}{\gamma_J}

\newcommand{\EE}{{E}_{\mathrm{T}}}
\newcommand{\DD}{{D}_{\mathrm{T}}}
\newcommand{\JJ}{{J}_{\mathrm{T}}}
\newcommand{\II}{{I}_{\mathrm{T}}}

\newcommand{\CC}{\mathbb{C}}
\newcommand{\NN}{\mathbb{N}}
\newcommand{\RR}{\mathbb{R}}

\newcommand{\cB}{{B}_{\Gamma}}
\newcommand{\Br}{{B_r}}
\newcommand{\Bc}{{B_c}}

\newcommand{\cL}{\mathcal{L}}
\newcommand{\cF}{\mathcal{F}}

\newcommand{\NF}{\hat{N}_{F}}

\newcommand{\Pp}{{P}_p}

\newcommand{\Qomega}{\hat{Q}_{\bm{\omega}}}
\newcommand{\QF}{{Q_{F}}}
\newcommand{\cR}{{R}} \newcommand{\cRh}{\hat{\cR}}

\newcommand{\rmq}{{q}}

\newcommand{\cU}{\mathcal{U}}

\newcommand{\rH}{\mathrm{H}}

\newcommand{\symp}{\hat{P}}
\newcommand{\symq}{\hat{Q}}
\newcommand{\syma}{\hat{A}}

\newcommand{\scrF}{\mathscr{F}}
\newcommand{\scrL}{\mathscr{L}}
\newcommand{\scrS}{\mathscr{S}}

\newcommand{\balpha}{\bm{\alpha}}

\newcommand{\bbeta}{\bm{\beta}}

\newcommand{\bc}{\bm{c}}
\newcommand{\bd}{\bm{d}}
\def\bf{\bm{f}}
\newcommand{\bg}{\bm{g}}
\newcommand{\bfg}{\bm{g}}\newcommand{\bfgh}{\hat{\bfg}}

\newcommand{\bk}{\bm{k}}
\newcommand{\bu}{\bm{u}}
\newcommand{\bv}{\bm{v}}\newcommand{\bvh}{\hat{\bv}}
\newcommand{\bw}{\bm{w}}\newcommand{\bwh}{\hat{\bw}}
\newcommand{\by}{\bm{y}}\newcommand{\byh}{\hat{\by}}
\newcommand{\bh}{\bm{h}}\newcommand{\bhh}{\hat{\bh}}

\newcommand{\psih}{\hat{\Psi}}

\newcommand{\bomega}{\bm{\omega}}
\newcommand{\bzeta}{\bm{\zeta}}
\newcommand{\brho}{{\bm{\rho}}}
\newcommand{\brhoh}{\hat{\bm{\rho}}}
\newcommand{\brhopm}{{\bm{\rho}_{\pm}}}
\newcommand{\brhohpm}{{\hat{\bm{\rho}}_{\pm}}}

\usepackage{amsopn}
\DeclareMathOperator{\diag}{diag}
\DeclareMathOperator{\dist}{dist}
\DeclareMathOperator{\supp}{supp}

\newcommand{\Mhat}{\hat{M}}
\newcommand{\Vhat}{\hat{V}}
\newcommand{\Mom}{\hat{M}_{\bomega}}
\newcommand{\Vom}{\hat{V}_{\bomega}}
\newcommand{\NJ}{\hat{N}_{J}}
\newcommand{\Phat}{\hat{P}}
\newcommand{\Dhat}{\hat{D}}
\newcommand{\Ohat}{\hat{O}}
\newcommand{\That}{\hat{T}}
\newcommand{\Tom}{\hat{T}_{\bm{\omega}}}

\renewcommand{\imath}{{\mathrm{i}\mkern1mu}}
\DeclareMathOperator{\OpsS}{Ops-S}

\newcommand{\Sop}{\hat{S}_{\Gamma}}
\renewcommand{\Re}[1]{\mathrm{Re}{#1}}
\renewcommand{\Im}[1]{\mathrm{Im}{#1}}

\usepackage{mathtools}
\usepackage{etoolbox}
\DeclareFontFamily{U}{matha}{\hyphenchar\font45}
\DeclareFontShape{U}{matha}{m}{n}{
<-6> matha5 <6-7> matha6 <7-8> matha7
<8-9> matha8 <9-10> matha9
<10-12> matha10 <12-> matha12
}{}
\DeclareSymbolFont{matha}{U}{matha}{m}{n}

\DeclareFontFamily{U}{mathx}{\hyphenchar\font45}
\DeclareFontShape{U}{mathx}{m}{n}{
<-6> mathx5 <6-7> mathx6 <7-8> mathx7
<8-9> mathx8 <9-10> mathx9
<10-12> mathx10 <12-> mathx12
}{}
\DeclareSymbolFont{mathx}{U}{mathx}{m}{n}

\DeclareMathDelimiter{\vvvert} {0}{matha}{"7E}{mathx}{"17}%
\DeclarePairedDelimiterX{\mnorm}[1]
{|}{|}
{\ifblank{#1}{\:\cdot\:}{#1}}

\usepackage{booktabs} 

\graphicspath{{figures/}} 

\ifpdf
	\DeclareGraphicsExtensions{.eps,.pdf,.png,.jpg}
\else
	\DeclareGraphicsExtensions{.eps}
\fi

\def\sqr#1#2{{\vcenter{\vbox{\hrule height.#2pt
              \hbox{\vrule width.#2pt height#1pt \kern#1pt \vrule width.#2pt}
              \hrule height.#2pt}}}}

\newtheorem{theorem}{\hskip 1.3em Theorem}[section]
\newtheorem{definition}[theorem]{\hskip 1.3em Definition}

\newtheorem{corollary}[theorem]{\hskip 1.3em Corollary}
\newtheorem{lemma}[theorem]{\hskip 1.3em Lemma}
\newtheorem{remark}[theorem]{\hskip 1.3em Remark}

\newtheorem{assumption}[theorem]{\hskip 1.3em Assumption}

\makeatletter
   
   \@addtoreset{equation}{section}
\makeatother

\begin{document}

\title{Kreiss stability analysis of Hagstrom-Warburton nonreflecting boundary conditions 
for the first-order time-dependent Maxwell equations
}

\author{Kaifang Liu\thanks{
School of Mathematics, Yangzhou University, Yangzhou 225009, Jiangsu, China. 
{\small\textit{e-mail}:} {\small \texttt{kliu@yzu.edu.cn}}.},
~~~Matthias Schlottbom\thanks{
Department of Applied Mathematics, University of Twente, P.O. Box 217, 7500 AE Enschede, The Netherlands.  
{\small\textit{e-mail}:} {\small \texttt{m.schlottbom@utwente.nl}}. },~~~
J.J.W. van der Vegt\thanks{
Department of Applied Mathematics, University of Twente, P.O. Box 217, 7500 AE Enschede, The Netherlands. 
{\small\textit{e-mail}:} {\small\texttt{j.j.w.vandervegt@utwente.nl}}. },~~~
  and~~~
Yan Xu\thanks{Corresponding author. 
School of Mathematical Sciences, University of Science and Technology of China, Hefei 230026, Anhui, China. 
Laoshan Laboratory, Qingdao 266237, P.R. China
{\small\textit{e-mail}:} {\small  \texttt{yxu@ustc.edu.cn}}.}
}

\date{}
\maketitle

\begin{abstract}
    We analyze the stability of the
    Hagstrom-Warburton nonreflecting boundary conditions (HW-NRBCs) for the first-order time-dependent Maxwell equations.
    The HW-NRBCs enjoy very small reflection coefficients and do not use high-order derivatives or nonlocal boundary operators, which makes them well-suited for high-order accurate numerical discretizations.
    The main result of this paper is an $L^2$ a-priori bound on the solution in terms of initial and boundary data and volume sources in a half-space.
    To obtain this result, we first derive a mapping that takes outgoing components to ingoing components of the solution of the Maxwell equations with HW-NRBCs, which shows that the Kreiss condition does not hold uniformly. Next, to prove stability, several symmetrizers are constructed, which establishes well-posedness in a generalized sense.
\end{abstract}



\maketitle

\section{Introduction}
Scattering problems arise in many application fields, such as seismic wave simulations, acoustic and electromagnetic wave transmission, and are modeled mathematically as partial differential equations defined on an unbounded domain.
In practical computations, the unbounded domain is frequently truncated to a finite computational domain with an artificial boundary enclosing all sources and scatterers.
In order to minimize unphysical reflections at the boundary two popular approaches emerged, nonreflecting boundary conditions (NRBCs) and perfectly matched
layers (PMLs). 

The perfectly matched layer \cite{Berenger1994}, originally introduced for the Maxwell equations, surrounds the computational domain with an absorbing layer that attenuates outgoing waves. Its effectiveness depends on the choice of absorption profile and layer thickness. After discretization, however, residual reflections from the outer boundary typically occur due to numerical dispersion and the finite PML thickness \cite{Pled2022}, but see \cite{Chern2019,Hojas2024} for essentially reflectionless discrete PMLs on structured grids for the linear scalar wave equation. Although easily implemented within standard schemes such as the finite element method and finite difference time-domain method, the PML introduces a substantial number of additional degrees of freedom—particularly in three dimensions—and requires careful parameter tuning to minimize reflections \cite{Cohen2002,Appelo2006}.

With NRBCs, the governing PDEs are solved on the truncated domain without modification, and the quality of the approximation depends on how effectively the NRBCs suppress spurious reflections at the artificial boundary.
A classical nonlocal NRBC formulation is the Dirichlet-to-Neumann (DtN) map \cite{Keller1989}; see \cite{Givoli2004,Grote1995,Grote1998,Hagstrom2007} for other types of nonlocal conditions. Although exact, such methods require solving a global problem and are therefore computationally expensive.
To reduce cost, local approximate NRBCs were introduced by Engquist and Majda \cite{Engquist1977}, with related formulations such as the Sommerfeld-type \cite{Givoli1992} and Bayliss–Turkel \cite{Bayliss1980} conditions. These lower-order NRBCs are exact in one dimension but less accurate in higher dimensions, as they perfectly absorb only waves incident normal to the boundary. To keep reflections small, the artificial boundary must therefore be placed sufficiently far from scatterers and sources, which limits their practicality in numerical simulations.
Several locally defined NRBCs with high-order accuracy have since been developed.
Key contributions include those of Engquist and Majda \cite{Engquist1977}, Mur \cite{Mur1981}, Collino \cite{Collino1993}, Higdon \cite{Higdon1986,Higdon1994}, Givoli and Neta \cite{Givoli2003}, and Hagstrom and Warburton \cite{Hagstrom2004}.
While the Engquist–Majda conditions become cumbersome at high order due to the presence of many derivatives, Collino \cite{Collino1993} introduced local high-order boundary conditions using auxiliary variables, greatly simplifying their implementation. The Givoli–Neta \cite{Givoli2003} and Hagstrom–Warburton (HW) \cite{Hagstrom2004} formulations further advanced this approach.
Among these, the HW conditions are particularly attractive: they exhibit high-order decay of reflection coefficients compared with the Higdon \cite{Higdon1994} and Givoli–Neta \cite{Givoli2003} NRBCs, and they facilitate the treatment of corners and edges when multiple artificial boundaries are present. These advantages extend to general hyperbolic systems, including the elastic wave \cite{Higdon1990}, shallow water \cite{Givoli2003High}, and Maxwell equations \cite{Hagstrom2004}. For recent developments on the HW-NRBCs, see \cite{Hagstrom2019,Hagstrom2023,Hagstrom2020,Hagstrom2009}.

A central issue is the well-posedness of initial boundary value problems with NRBCs.
Energy estimates are commonly used to establish well-posedness for dissipative boundary conditions \cite{Kreiss2004} but may fail in more general settings \cite[Section 7.4]{Kreiss2004}. 
Energy estimate techniques have been used for the Higdon-NRBCs and the HW-NRBCs for the scalar wave equation \cite{Baffet2011}, 
but a stability analysis for the time-dependent Maxwell equations with HW-NRBCs is not available.

An alternative approach to prove well-posedness relies on the uniform Kreiss condition -- also known as the uniform Lopatinskii condition or the Kreiss-Sakamoto condition -- \cite{Kreiss1970,Majda1975,Higdon1986}, which provides necessary and sufficient algebraic criteria for the strong well-posedness of hyperbolic systems with inhomogeneous boundary conditions. 
In \cite{Kreiss1970}, Kreiss established the theory for strictly hyperbolic systems with noncharacteristic boundaries, while Majda and Osher \cite{Majda1975} extended it to characteristic boundaries, encompassing systems such as the time-dependent Maxwell and shallow water equations. Higdon \cite{Higdon1986} later offered a physical interpretation of this theory. For comprehensive treatments and recent developments, see \cite{BenzoniGavage2006,Kreiss2004,Motamed2019}.
Failure of the uniform Kreiss condition leads to weaker notions of well-posedness, typically accompanied by a loss of derivatives. 
Examples of such behavior arise, for instance, in \cite{Eller2009} for Maxwell’s equations with (conservative) perfectly conducting boundary conditions, in \cite{Eller2012} for general hyperbolic systems with conservative boundary conditions, and in \cite{Eller2018} for constant-coefficient hyperbolic problems posed on half-spaces. 
An approach based on symmetrizers for noncharacteristic boundaries is used in \cite[Section~2]{Coulombel2004}.

In this paper, we analyze the Hagstrom–Warburton nonreflecting boundary conditions for the first-order time-dependent Maxwell equations with general (non-smooth) material coefficients in the interior of the half-space $\Omega=\{\bx\in \RR^3:\,x_1<0\}$. 
Our approach uses the energy method for the interior domain \cite{BenzoniGavage2006} where the coefficients can be non-smooth. The energy method requires the boundary conditions to be maximally dissipative, which is not the case for the HW-NRBCs. In the neighborhood of the boundary, where the coefficients are constant, we use therefore the Laplace–Fourier domain techniques developed in \cite{Kreiss1970,Majda1975}.
A key element of the proof is the construction of eigensolutions of the interior Maxwell system via the bilateral Laplace–Fourier transform and converting the HW-NRBCs —originally time-dependent PDEs on the artificial boundary— into an algebraic system in the frequency domain.
Our analysis reveals that the uniform Kreiss condition, associated with this algebraic system, is not satisfied; consequently, one of the key assumptions in \cite{Majda1975} fails, and their results on strong well-posedness with inhomogeneous boundary data do not apply. The failure of the uniform Kreiss condition arises from the presence of generalized eigenvalues. 
For strictly hyperbolic systems with noncharacteristic boundaries, well-posedness can nevertheless be established for homogeneous boundary conditions, provided that the associated generalized eigenvalues satisfy certain additional structural properties \cite{Motamed2019}. 
We show that the generalized eigenvalues corresponding to the HW--NRBCs for the Maxwell equations also satisfy the properties assumed in \cite{Motamed2019}.
In the case of inhomogeneous boundary data, well-posedness may still hold, though typically with a loss of derivatives \cite{Eller2018}. 
While \cite{Eller2018} obtains a well-posedness result through the construction of explicit solutions, which, however, does not generalize to variable coefficients and cannot handle lower-order terms, we construct appropriate symmetrizers to obtain stability bounds also suitable for  more general cases.
Particular attention is devoted to points on the imaginary axis where the Laplace–Fourier-transformed Maxwell matrix is non-diagonalizable or singular, significantly complicating the symmetrizer construction.
After the construction of a symmetrizer it is straightforward --using G\r{a}rding's inequality-- to obtain a bound on the solution in terms of the data, which proves the stability and well-posedness of the time-dependent Maxwell equations with HW-NRBCs.
As a byproduct of our analysis, we derive explicit reflection coefficients for the HW-NRBCs, extending analogous results for the scalar wave equation and the equations for acoustics \cite{Givoli2006, Hagstrom2008, Hagstrom2004,Hagstrom2009}. 

The organization of the paper is as follows. Section~\ref{sec:HW hyperbolic} introduces the Maxwell system together with the HW--NRBCs and states our main a priori estimate. 
Sections~\ref{sec:analysis-HW}--\ref{sec:half-space} are devoted to the analysis of half-space problems. Section~\ref{sec:analysis-HW} constructs eigensolutions of the constant-coefficient Maxwell system and derives the associated reflection matrix, allowing us to express the HW--NRBCs as a mapping from outgoing to incoming modes in the Laplace--Fourier domain (Section~\ref{sec:ukc}). This analysis also reveals the failure of the uniform Kreiss condition and shows that three distinct cases must be treated separately. In Section~\ref{sec:symmetrizer}, we construct symmetrizers adapted to these cases, which are then used in Section~\ref{sec:half-space} to establish an a priori bound in half-spaces. In the appendices, we provide a detailed derivation of the HW-NRBCs, several technical lemmas, and G{\aa}rding inequalities used in the analysis.

\section{Maxwell equations with HW-NRBCs}
\label{sec:HW hyperbolic}
 
The Maxwell equations for linear materials on the half-space $\Omega=\{\bx\in \RR^3:\,x_1<0\}$ are defined as
\begin{subequations}
    \label{eq:Maxwell time-dependent}
\begin{alignat}{3}
    \begin{pmatrix}
        \varepsilon_r &  0\\ 0 &\mu_r
    \end{pmatrix}
    \frac{\partial }{\partial t} \begin{pmatrix}
        \bm{E} \\ \bm{H}
    \end{pmatrix} &= \begin{pmatrix}
        0 & \nabla\times \\ -\nabla\times & 0
    \end{pmatrix}  \begin{pmatrix}
        \bm{E} \\ \bm{H}
    \end{pmatrix}+\begin{pmatrix}
        \bm{j}\\\mathbf{0}
    \end{pmatrix}, \quad && \text{in } \Omega \times (0,\infty),\label{eq:Maxwell time-dependent 1} \\
    \nabla\cdot (\varepsilon_r\bm{ E}) &= \varrho,&& \text{in } \Omega \times (0,\infty),
    \label{eq:Maxwell time-dependent 3} \\
    \nabla\cdot (\mu_r\bm{H}) &=0, && \text{in } \Omega \times (0,\infty),\label{eq:Maxwell time-dependent 4} 
\end{alignat}
\end{subequations}
with initial conditions
\begin{equation}
    \label{eq:Maxwell-time-dependent-4}
    \bm{E}(\cdot,0)=\bm{E}_0,\quad
    \bm{H}(\cdot,0)=\bm{H}_0\quad \text{in } \Omega,
\end{equation}
and boundary conditions that will be specified in the next section.
Here, \(\bm{E}\), \(\bm{H}:\Omega \times
(0,\infty)\to \mathbb{R}^3\) denote, respectively, the electric field
and magnetic field with initial fields $\bm{E}_0, \bm{H}_0$;
\(\nabla \) is the nabla operator, $\nabla\times$ the curl operator, $\nabla\cdot$ the divergence operator and \(t\) is time.
The parameter functions \(\varepsilon_r,\,  \mu_r\) represent, respectively, the relative electric
permittivity and relative magnetic permeability.
The current density, denoted by \(\bm{j}:\Omega \times
(0,\infty)\to \mathbb{R}^3\), is related to the charge density
\(\varrho:\Omega \times (0,\infty)\to \mathbb{R}\) by the conservation
law
\begin{equation}
    \label{eq:Maxwell time-dependent 5}
    \frac{\partial \varrho}{\partial t}=\nabla\cdot \bm{j}. 
\end{equation}
We introduce the scaled variables
\begin{align*}
    \bm{u}=\left[ \sqrt{\bar{\varepsilon}_r \bar{\mu}_r}\bm{E}^\mathrm{T},\bar{\mu}_r \bm{H}^\mathrm{T} \right]^{\mathrm{T}},\quad
    \bm{f}=\left[ \sqrt{\bar{\mu}_r/\bar{\varepsilon}_r }\bm{j}^\mathrm{T},\bm{0}\right]^{\mathrm{T}},
\end{align*}
with \(\bar{\varepsilon}_r, \,\bar{\mu}_r\) reference values and use the transformation  \(\tilde{\bm{x}}=\sqrt{\bar{\varepsilon}_r \bar{\mu}_r} \bm{x}\), which implies that \(\tilde{\nabla} = (\bar{\varepsilon}_r \bar{\mu}_r)^{-1/2} \nabla \). 
The Maxwell equations \eqref{eq:Maxwell time-dependent} can then be formulated as the symmetric hyperbolic system
\begin{subequations}
    \label{eq:hyperbolic-sys}
    \begin{alignat}{3}
        &L \bm{u} =
        B\frac{\partial \bm{u}}{\partial t}  + \sum_{i=1}^3 A_i \frac{\partial \bm{u}}{\partial \tilde{x}_i}  = \bm{f}, &\quad& \text{in } \Omega\times (0,\infty),\label{eq:hyperbolic-sys-1} \\
        &\bm{u}(\cdot,0) = \bm{u}_0, && \text{in } \Omega, \label{eq:hyperbolic-sys-2}
    \end{alignat}
\end{subequations}
with the symmetric matrices 
\begin{equation}
    A_i=\begin{pmatrix} 0& -C_i\\ C_i & 0
    \end{pmatrix}\in \mathbb{R}^{6\times 6} \quad \text{and} \quad
    B(\tilde{\bm{x}})=\begin{pmatrix} \frac{\varepsilon_r(\tilde{\bm{x}})}{\bar{\varepsilon}_r} I_3 & \\  & \frac{\mu_r(\tilde{\bm{x}})}{\bar{\mu}_r} I_3
    \end{pmatrix}\in \mathbb{R}^{6\times 6},
    \label{eq:Ai-B}
\end{equation}
where
\[
C_1 = \begin{pmatrix}
    0 & 0 & 0 \\ 0 & 0 & -1\\ 0 & 1 & 0
\end{pmatrix},\  C_2 = \begin{pmatrix}
    0 & 0 & 1 \\ 0 & 0 & 0\\ -1 & 0 & 0
\end{pmatrix},\  C_3 = \begin{pmatrix}
    0 & -1 & 0 \\ 1 & 0 & 0\\ 0 & 0 & 0
\end{pmatrix},
\]
and \(I_n \in \mathbb{R}^{n \times n}\), \(n \in \mathbb{N}\), the identity matrix.
If $\varepsilon_r$ and $\mu_r$ are constant in $\Omega$ then $B(\tilde{\bm{x}})=I_6$.
In order to simplify notation, we will use \(\bm{x}\) instead of \(\tilde{\bm{x}}\) for the Cartesian coordinates in the remainder.

The HW-NRBCs were introduced by Hagstrom and Warburton in \cite{Hagstrom2004} for the scalar wave equation and hyperbolic systems with non-characteristic boundary.
In Appendix \ref{app:deriv-H-W-new}, the HW-NRBCs for the Maxwell equations are derived for a semi-infinite domain $\Omega_n$ with boundary $\Gamma_n$ and external normal vector $\bm{n}$. For the half-domain $\Omega$, we only need to take $\bm{n}=(1,0,0)^\top$, and obtain the corresponding HW-NRBCs along the boundary $\Gamma=\{\bm{x}:\ x_1=0\}$. 
In fact, define \(\bm{\phi}_0:\Gamma \times (0,\infty) \to \mathbb{R}^6,\) with \(\bm{\phi}_0:=\bm{u}|_{\Gamma}\), and introduce the auxiliary variables \(\bm{\phi}_j:\Gamma\times (0,\infty) \to \mathbb{R}^6, j=1,2,\dots,J+1\). The sequence of HW-NRBCs for the auxiliary variables \(\bm{\phi}_j\) can be stated at the boundary \(\Gamma\) as
\begin{subequations}
    \label{eq:hw-nrbcs}
    \begin{alignat}{3}
        A_1\frac{\partial \bm{\phi}_{1}}{\partial t}-\mathcal{L}_0\bm{\phi}_0
        &=\cL_0 \bm{g},&& \text{at } \Gamma \times (0,\infty), \label{eq:rho2-max-3-1} \\
        \mathcal{L}_j\bm{\phi}_j - 
        \mathcal{K}_j\bm{\phi}_{j+1} &=\bm{0},\quad j=1,2, \dots,J,
        \quad && \text{at } \Gamma \times (0,\infty),\label{eq:rho2-max-3-2} \\
        \cB^-\bm{\phi}_{J+1} &= \bm{0},&& \text{at } \Gamma \times (0,\infty),\label{eq:rho2-max-3-3}\\
        \cB^0\mathcal{L}_1 \bm{\phi}_1 &= \bm{0}, && \text{at } \Gamma \times (0,\infty),\label{eq:rho2-max-3-4}\\
        \bm{\phi}_j(\bm{x},0)&=\bm{0}, \quad j=1,2,\dots,J+1,\quad && \text{at } \Gamma .
    \end{alignat}
\end{subequations}
Here, \(\mathcal{L}_j, \mathcal{K}_j\) are linear operators defined as
\begin{subequations}
\label{eq:KjLj}
    \begin{align}
    \mathcal{L}_j &= (a_jA_1 -I_6)\frac{\partial }{\partial t} -A_2 \frac{\partial }{\partial x_2}
        -A_3 \frac{\partial }{\partial x_3}, \quad j=0,1,2,\dots,J,\label{eq:operatorLj}\\
        \mathcal{K}_j &= (a_jA_1 +I_6)\frac{\partial }{\partial t} +A_2 \frac{\partial }{\partial x_2}
        +A_3 \frac{\partial }{\partial x_3}, \quad j=1,2,\dots,J,\label{eq:operatorKj}
\end{align}
\end{subequations}
and the parameters \(0<a_0\leq 1,\, 0<a_j< 1\) are to be specified to minimize wave reflections at the boundary \(\Gamma\). 
The matrix $A_1$ has the following eigen-decomposition
\begin{equation}
    \label{eq:matA1}
    A_1 = S \Lambda S^\mathrm{T},\quad \text{with}\ 
    \Lambda=\begin{bmatrix} \Lambda_+ & &\\ & \Lambda_0 &\\ & & \Lambda_- \\
    \end{bmatrix}=\text{diag}( 1,1,0,0,-1,-1 ),
\end{equation}
and
\begin{equation}
    \label{eq:Smat}
    S=\begin{bmatrix}
        &  & 1 &  & &  \\
        & \frac{\sqrt{2}}{2} & & & & -\frac{\sqrt{2}}{2} \\
        -\frac{\sqrt{2}}{2}&  &  & & \frac{\sqrt{2}}{2} & \\
        &  & & 1 &  & \\
        \frac{\sqrt{2}}{2} & & &  & \frac{\sqrt{2}}{2} & \\
        & \frac{\sqrt{2}}{2} & &  & & \frac{\sqrt{2}}{2}
    \end{bmatrix},
\end{equation}
where $S$  is a unitary matrix, hence \(S^{-1} = S^\top\). We assume that the boundary condition contains only incoming information, i.e., $\cB^0\bg=\bm{0}$ and $\cB^+\bg=\bm{0}$. The matrices $\cB^-,\,\cB^0,\,\cB^+$ consist, respectively, of rows $1$-$2$, $3$-$4$ and $5$-$6$ of  $S^\top$.
See Appendix \ref{sec:app:mat-decomp} for the decomposition of $A_1$ for a plane with a normal vector $\bm{n}$. 

Using the matrix $S$, we transform the Maxwell equations into characteristic variables.  Define $\brho=S^{\top}\bm{\bu}$ and write \eqref{eq:hyperbolic-sys-1}  as
\begin{equation}
\label{eq:hyperbolic-sys-1-trans}
    \frac{\partial \brho}{\partial t}  + \Lambda \frac{\partial \brho}{\partial x_1}+\sum_{i=1}^2 \widetilde{A}_i \frac{\partial \brho}{\partial \tilde{x}_i}  = S^{\top}\bm{f}.
\end{equation}
Here, the matrices $\widetilde{A}_2,\,\widetilde{A}_3 \in \mathbb{R}^{6 \times 6}$ have the form
\begin{equation}
    \label{eq:A2A3-tilde}
    \widetilde{A}_2=S^{\top }A_2 S=\begin{pmatrix}
    & -\EE & \\ -\EE& & -\DD \\ & \DD &
    \end{pmatrix} ,\quad  \widetilde{A}_3=S^{\top }A_3 S=\begin{pmatrix}
        & \JJ & \\ \JJ& & \II \\ & \II &
    \end{pmatrix}
\end{equation}
with 
\[
    \EE=\frac{\sqrt{2}}{2}\begin{pmatrix}0 & 1\\ 1 & 0\end{pmatrix},\ 
    \DD=\frac{\sqrt{2}}{2}\begin{pmatrix}0 & 1\\ -1 & 0\end{pmatrix},\  
    \II=\frac{\sqrt{2}}{2}\begin{pmatrix}1 & 0\\ 0 & 1\end{pmatrix},\  
    \JJ=\frac{\sqrt{2}}{2}\begin{pmatrix}1 & 0\\ 0 & -1\end{pmatrix}.
\]

\begin{remark}
    If the initial conditions \(\bm{E}_0\) and
\(\bm{H}_0\) in \eqref{eq:Maxwell-time-dependent-4} satisfy \(\nabla\cdot (\varepsilon_r\bm{E}_0)=\varrho\) and
\(\nabla\cdot \mu_r\bm{H}_0=0\) in \(\Omega\), then using \eqref{eq:Maxwell time-dependent 5} these conditions
are automatically satisfied by the solution to the Maxwell equations for all later
times. 
Therefore, conditions \eqref{eq:Maxwell time-dependent 3}, \eqref{eq:Maxwell time-dependent 4} do not need to be considered explicitly in
the stability analysis.
Since the properties of the source term $\bf,$ initial data $\bu_0$ and boundary data $\bg_-$ do not show up in the analysis we will consider in the remainder general $\bf,\bu_0$ and $\bg_-$ that satisfy Assumption \ref{ass:1} below.
\end{remark}

\subsection{Main results of stability analysis}
\label{ssec:main-results}
To state our main result, we first introduce some function spaces and norms.
Let $\Sigma\subseteq\RR^n$ be an open subset. We denote by $L^p(\Sigma)$, $1 \leq p \leq \infty$, the standard Lebesgue spaces. 
The space of infinitely differentiable functions with compact support in $\Sigma$ is denoted by $C_0^{\infty}(\Sigma)$.
For a nonnegative integer $k$, coefficient $\eta\in \RR$, $\eta>0$, and space-time domain \(V\) with smooth boundary satisfying the segment property, see \cite{Adams2003, Wloka1987}, the space $H_\eta^k(V)$ is the completion of $C_0^\infty(\overline{V})$ in the norm
\begin{equation}
    \label{eq:space-time-norm}
    \|v\|_{k,V, \eta}^2
    = \sum_{|\balpha|\leq k}\int_V\e^{-2 \eta t}|D^{\balpha}v|^2,
\end{equation}
see \cite{Majda1975}.
The space-time domain $V$ is of the form
\begin{equation}
    \label{eq:space-V}
    V=\Sigma\times (0,\infty)\subset \RR^{n+1}\quad\text{or}\quad 
    V=\partial\Sigma\times (0,\infty)\subset \RR^n,
\end{equation}
with $\bx\in \Sigma$ or $\bx \in \partial \Sigma$ and $t\in (0,\infty)$.
The related inner product on $H_\eta^k(V)$ is denoted by $(\cdot,\cdot)_{V,\eta}$, and for $\bv=(v_1,...,v_\ell)^\top\in \RR^\ell$, $\ell\in\NN,$ we have $\|\bv\|_{k,V,\eta}^2=\sum_{i=1}^\ell \|v_i\|_{k,V,\eta}^2$. 
If $\eta=0$ we drop the $\eta$ subscript on the norm and inner-product.

In the remainder we will consider solutions  $\bu\in (H_\eta^1(\Omega\times(0,\infty)))^6$, with trace in $(H_\eta^1(\Gamma\times(0,\infty)))^6$, to the Maxwell equations \eqref{eq:hyperbolic-sys} and HW-NRBCs \eqref{eq:hw-nrbcs} that hold pointwise almost everywhere.

Our approach is based on a separate treatment of the interior of $\Omega$ and a neighborhood of its boundary $\Gamma$.
To state our results, let $\epsilon>0$ and denote with
$$
\Omega_{\epsilon}=\{x\in{\Omega}: \dist(x,\Gamma)<\epsilon\},
$$
an $\epsilon$-neighborhood of the boundary.
We will always assume that the following assumption holds.

\begin{assumption}\label{ass:1}
(i) The parameters $\varepsilon_r,\,\mu_r\in L^\infty(\Omega)$ are positive and real in $\Omega$.\\
(ii) The source term $\bf\in (L^2(\Omega\times(0,\infty)))^6$, the initial data $\bu_0\in (L^2(\Omega))^6$ and the boundary data $\bg_-\in (H_\eta^1(\Gamma\times(0,\infty)))^6$.\\
(iii) There is $\epsilon>0$ such that for all $t>0$, $\bf(\bx,t)=\bm{0}$ and $\bu_0(\bx)=\bm{0}$ for all $\bx\in\Omega_{\epsilon}$.
Moreover, $\varepsilon_r=\bar{\varepsilon}_r,\, \mu_r=\bar{\mu}_r$ in $\Omega_\epsilon$ with $\bar{\varepsilon}_r,\,\bar{\mu}_r$ constant reference values.
\end{assumption}

\begin{theorem}[$L^2$-Stability]
    \label{thm:main-thm}
    Let $\eta\in \RR,\,\eta>0$ be sufficiently large and suppose Assumption 2.2 holds.
    Assume that the solution $\bu$ of the Maxwell equations \eqref{eq:hyperbolic-sys} 
    with HW-NRBCs \eqref{eq:hw-nrbcs} satisfies $\bu\in (H_\eta^1(\Omega\times(0,\infty)))^6$ with trace in $(H_\eta^1(\Gamma\times(0,\infty)))^6$.
    Then there exists a positive constant $C$, independent of $\bu$ and $\eta$, such that 
    the solution $\bm{u}$ satisfies the $L^2$-stability estimate
    \begin{align*}
     \|\bm{u}\|_{0,\Gamma \times [0,\infty),\eta}^2
        &+\eta \|\bm{u}\|_{0,\Omega \times [0,\infty),\eta}^2 \\
        &\leq C  \Big(\frac{1}{\eta}\|\bm{f}\|_{0,\Omega \times [0,\infty),\eta}^2 + \|\bu_0\|_{0,\Omega}^2 
         + \|\bm{g}_-\|_{1,\Gamma\times[0,\infty),\eta}^2
        \Big).
    \end{align*}
\end{theorem}

The proof of Theorem \ref{thm:main-thm} follows directly from Theorem \ref{thm:L2-stability} and proceeds by first analyzing the half-space model by constructing eigenmodes and the associated reflection matrix, expressing the boundary conditions as a map from outgoing to incoming modes in the Laplace-Fourier domain \cite{Kreiss1970,Majda1975}. Since the analysis shows that the uniform Kreiss condition fails \cite{Kreiss1970,Majda1975}, we split the analysis into three distinct cases depending on the spectral properties of the Laplace-Fourier transformed Maxwell operator, which implies the matrix symbol $\Mhat$ \eqref{eq:matrix-M} below. For each case, we build a suitable symmetrizer, whose construction only depends on $\Mhat$ but must satisfy a condition related to the NRBCs. Then, we derive stability estimates in $\Omega_\epsilon$ using G{\aa}rding's inequality. Finally, after combining this result with energy estimates in $\Omega\setminus\Omega_\epsilon$ an a-priori bound is obtained in the half-space. 

In the analysis we use the theory of pseudo-differential operators, in particular, several versions of G{\aa}rding's inequality that can deal with the singularity in the Laplace domain. This general approach is suitable to extend the analysis to general domains with smooth boundaries, see e.g. \cite{Kreiss2004}, Theorem 4 or \cite{BenzoniGavage2006}. This extension will be part of future work.

We emphasize that the symmetrizers constructed in this work are not tied to the specific structure of the HW–NRBCs and therefore apply to a broader class of boundary conditions for the Maxwell system. 
As can be seen from Section~\ref{sec:symmetrizer}, only the verification of the boundary-related estimates—specifically the statements of Lemma~\ref{lem:case1_symmetrizer_fulfills_property_ii}, Lemma~\ref{lem:case2_property_ii_holds} and Lemma~\ref{lem:case3_property_ii_holds}—needs to be adapted when considering alternative boundary conditions. Consequently, the symmetrizer framework developed here provides a flexible tool for analyzing the stability of a wide range of boundary conditions for the Maxwell equations beyond HW–NRBCs.

\section{Solution of constant coefficient Maxwell equations with HW-NRBCs on a semi-infinite domain}
\label{sec:analysis-HW}

\subsection{Definitions, function spaces}
\label{ssec:defs}
Let $\Sigma\subseteq \RR^n$ be an open set. The space of infinitely differentiable functions on $\Sigma$ is denoted by $C^\infty(\Sigma)$. The space of $C^\infty$-functions with bounded derivatives on $\overline{\Sigma}$ is denoted by $C_b^{\infty}(\overline{\Sigma})$. The space $C_{(0)}^\infty(\overline{\Sigma})$ is defined as 
$$C_{(0)}^{\infty}(\overline{\Sigma})=\{v\in C^\infty(\overline{\Sigma}):\ \mathrm{supp}(v)\, \mathrm{compact}\, \subset \overline{\Sigma}\}.
$$

We recall the definitions of the Fourier transform, the Schwartz and $L^2$-Bessel potential spaces, and the pseudo-differential operators and symbols used in this paper.

\begin{definition}[Fourier Transform]
    For $v \in L^1(\mathbb{R}^n)$, the  Fourier transform is defined as
    \[
        \hat{v}(\bomega)=\scrF[v](\bomega) = \int_{\mathbb{R}^n}v (\bm{x} )
        \e^{-\imath \bomega \cdot \bm{x} } \,\mathrm{d}  \bm{x},
    \]
    with \(\bomega \in \mathbb{R}^n\) the dual variable related to \(\bm{x} \in \mathbb{R}^n\), and \(\imath=\sqrt{-1}\).
    For $\hat{v}\in L^1(\mathbb{R}^n)$ the inverse Fourier transform is given by
    \[
        v(\bm{x}) = \scrF^{-1}[\hat{v}](\bm{x}) =\frac{1}{(2\pi)^n} 
        \int_{\mathbb{R}^n}\hat{v}(\bomega) \e^{\imath \bomega \cdot \bm{x}} \,\mathrm{d}  \bomega.
    \]
\end{definition}
\begin{definition}
    The space $\scrS(\mathbb{R}^n)$ of rapidly decreasing (Schwartz) functions is the set of all smooth functions $v:\mathbb{R}^n \to \mathbb{C}$ such that for all multi-indices $\balpha  \in \mathbb{N}_0^n$ and $N\in\NN_0$ with $\NN_0=\NN\cup \{0\}$, there is a constant $C_{\balpha,N}$ such that
    \[
        |\partial_{\bx}^{\balpha} v(\bm{x})| \leq C_{\balpha,N}(1+|\bm{x}|)^{-N}\quad\text{uniformly for all }\bm{x}\in\RR^n.
    \]
    The space of tempered distributions, $\mathscr{S}^\prime (\mathbb{R}^n)=(\mathscr{S}(\mathbb{R}^n))^\prime $ is the space of linear functionals $\ell:\mathscr{S}(\mathbb{R}^n)\to \mathbb{C}$ such that there is $m\in\NN_0$ and $C>0$ such that $|\ell(v)|\leq C |v|_{m,\scrS}$ for all $v\in\scrS(\RR^n)$, where $|v|_{m,\scrS}=\sup_{|\balpha|+|\bbeta|\leq m}\sup_{\bx\in\RR^n}|\bx^{\balpha}\partial_{\bx}^{\bbeta }v(\bx)|$.
\end{definition}

Next, we extend the Fourier transform to the space of tempered distributions.

\begin{definition}
        Let $v \in \scrS^\prime (\mathbb{R}^n)$. The Fourier transform $\scrF[v]$ and its inverse $\scrF^{-1}[v]$ of $v$ are defined as the distributions
        \begin{alignat*}{3}
            \langle \scrF[v],\phi \rangle&= \langle v, \scrF[\phi]\rangle&\quad &\forall \phi \in \scrS(\mathbb{R}^n),\\
            \langle \scrF^{-1}[v],\phi \rangle&= \langle v, \scrF^{-1}[\phi]\rangle&& \forall \phi \in \scrS(\mathbb{R}^n),
        \end{alignat*}
    with $\langle v,\phi\rangle=v(\phi)$, $\phi \in \scrS(\mathbb{R}^n)$ the duality product.
\end{definition}
\begin{definition}[Bilateral Laplace Transform]
    For $g:[0,\infty)\to\RR$ let 
    \[
        \bar{g}(t)=\begin{cases}
            g(t),&t>0,\\
            \frac{1}{2}g(t),&t=0,\\
            0,&t<0.
        \end{cases}
    \]
    If $t\mapsto \e^{-\eta t}\bar{g}(t) \in\scrS^\prime(\mathbb{R})$ for some $\eta\in \RR,\,\eta>0$, then the bilateral Laplace transform \cite{Kreiss2004} is given by
    \[
        \hat{g}(\xi)=\scrL[g](s) 
        = \scrF[\e^{-\eta t}\bar{g}](\xi)
        = \int_{-\infty }^{+\infty }\e^{-st}\bar{g}(t)\,\mathrm{d}t,\quad \xi\in\RR,\ s=\eta+\imath \xi.
    \]
\end{definition}
\begin{definition}[Inverse bilateral Laplace-Fourier transform]
    Given $\hat{v}(\bomega,s)\in \scrS(\RR^{n+1}) $ with $\bomega\in \RR^n$, $s=\eta+\imath \xi \in \CC$, $\eta>0, $ $\eta,\,\xi\in \RR$. The inverse bilateral Laplace-Fourier transform is defined as
    $$
    (\cL\cF)^{-1}_{n+1}[\hat{v}](\bx,t;\eta)=\frac{1}{(2\pi)^{n+1}}\int_{\RR^{n+1}}e^{\imath \xi t+\bomega\cdot\bx}\hat{v}(\bomega,\imath \xi+\eta)\mathrm{d}\bomega \mathrm{d}\xi
    $$
    with $\hat{v}(\bomega,\imath \xi+\eta)=\int_{\RR^{n+1}}e^{-\imath \xi-\bomega\cdot\bx-\eta t}v(\bx,t)\mathrm{d}\bx \mathrm{d}t$, where $v(\bx,t)=0$ for $t<0$.
\end{definition}
\begin{definition}
    The $L^2$-Bessel potential spaces $H^s(\RR^n)$, $s \in \RR$ and norm $\|\cdot \|_s$ are, respectively, defined by
    \[
        H^s(\RR^n)= \left\{ v \in \scrS^\prime (\RR^n):\ \langle D_{x} \rangle^s v \in L^2(\mathbb{R}^n)  \right\},
    \]
    \[
        \|v\|_s^2= \int_{\mathbb{R}^n} \left| \langle D_{\bx} \rangle^s v(\bm{x}) \right|^2 \,\mathrm{d}\bm{x} 
        = \frac{1}{(2 \pi)^n}\int_{\mathbb{R}^n} \langle \bomega \rangle^{2s}\left|  \hat{v}(\bomega) \right|^2 \,\mathrm{d}\bomega, 
    \]
    where $\langle \bomega\rangle = (1+|\bomega|^2)^{1/2}$ and $\langle D_{\bx} \rangle^s v = \scrF^{-1}[\langle\bomega\rangle^s\hat v]$ for all $v\in\scrS^\prime(\RR^n)$.
\end{definition}

For the definition of $H^s(\Sigma),\,s\in \RR$, $\Sigma$ an open subset in $\RR^n$, and associated norm $\|\cdot\|_{s,\Sigma}$, see \cite[Page 77]{McLean2000}.

We also define the related weighted Bessel potential spaces $\widetilde{H}_\eta^s(V)$ with $s\in \RR,$ $\eta\in \RR$, $\eta>0$ with norm
\begin{equation}
    \label{eq:space-time-norm-s}
    \|v\|_{s,\eta}^2 = \frac{1}{(2\pi)^{n+1}}\int_{\RR^{n+1}}\langle \bomega\rangle_\eta^{2s}|\hat{v}(\bomega)|^2 \,\mathrm{d}\bomega,
\end{equation}
where $\langle \bomega\rangle_\eta=(1+\eta^2+|\bomega|^2)^{1/2}$. 
The restriction of \eqref{eq:space-time-norm-s} to the space-time domain $V$ \eqref{eq:space-V} is denoted $\|\cdot\|_{s,V,\eta}$.
Note, for $s=k$, integer, using $e^{-\eta t}\frac{\partial \bv}{\partial t}=\left(\eta+\frac{\partial}{\partial t}\right)(e^{-\eta t}\bv)$ we obtain that the norms \eqref{eq:space-time-norm} and \eqref{eq:space-time-norm-s}, applied to the space-time domain \eqref{eq:space-V} are identical. 
Using the fact that functions in $C_{(0)}^{\infty}(\overline{V})$ are dense in $H^s(V)$ with $V=\RR^{n+1}$, $\RR^{n+1}_+$, or a smooth open set $V\subset \RR^{n+1}$ with the segment property, see e.g. \cite{Grubb1996}, Theorem 4.10, we obtain that for $k=s$, $H_\eta^k(V)\cong\widetilde{H}_\eta^s(V)$. In the remainder, we will therefore drop the $\ \widetilde{}\ $ on $\widetilde{H}_\eta^s(V)$ and use $H_\eta^s(V)$ with $s\in \RR$.

\begin{definition}
    A smooth function $\symp:C^\infty(\Sigma \times \RR^n)$, with $\Sigma\subseteq\RR^N$ an open set, is a symbol of class $S_{1,0}^m(\Sigma\times\RR^n)$ and order $m \in \mathbb{R}$, if for any multi-indices $\balpha \in \mathbb{N}_0^n$, $\bbeta \in \mathbb{N}_0^N$, and any compact set $K\subset \Sigma$, there exists a constant $C_{\balpha,\bbeta,K}$, such that
    \[
        \left| \partial_{\bomega}^{\balpha}\partial_{\bm{x}}^{\bbeta}\symp(\bm{x},\bomega) \right|
        \leq C_{\balpha,\bbeta,K} \langle \bomega \rangle^{m-|\balpha|}.
    \]
\end{definition}

\begin{definition}
    A linear operator $P: C_0^\infty(\Sigma)\to C_0^\infty(\Sigma)$, with $\Sigma\subseteq \RR^N$ an open set, is a pseudo-differential operator with symbol $\symp(\bm{x},\bomega)\in S_{1,0}^m(\Sigma \times \mathbb{R}^n)$, $m\in \mathbb{R}$, if $Pv$ can be represented by
    \[
        Pv(\bm{x}) = \frac{1}{(2 \pi)^n} \int_{\mathbb{R}^n} \e^{\imath \bm{x}\cdot \bomega}
        \symp(\bm{x},\bomega) \hat{v}(\bomega) \,\mathrm{d}\bomega .
    \]
\end{definition}
In this case, we write $P=P(\bm{x},\imath D_{\bm{x}})\in \OpsS_{1,0}^m(\Sigma\times \RR^n)$ where $D_{\bm{x}}=(D_{x_1},\dots,D_{x_N})$ with $\imath D_{x_j}=\partial_{x_j}$, $\imath=\sqrt{-1}$.

\begin{definition}
    A linear operator $P=P(\bx,t,\imath D_{\bx},\imath D_t+\eta)$ $\in \OpsS_{1,0}^{m_1,m_2}(\Sigma\times \RR^n)$ if $s^{-m_1}\symp((\bx,t),(\bomega,s))\in S_{1,0}^{m_2}(\Sigma\times \RR^n)$ for $\eta>0$ fixed, with $m_1,m_2\in \RR$, $\imath D_t=\frac{\partial}{\partial t}$.
\end{definition}

Given a function $\syma \in C_b^{\infty}(\Sigma\times \RR^n\setminus\{0\})$ that is homogeneous of degree $m$ in $\bxi$, which implies, $\hat{A}(\bx,t\bxi)=t^m \hat{A}(\bx,\bxi)$ for $\bxi\in \RR^n$, $t\in \RR, \,t>0$.  Given a function $\chi \in C_b^\infty(\RR^n)$ that vanishes in a neighborhood of $\bm{0}$ and $\chi(\bxi)=1$ if $|\bxi|\geq 1$. Then $\widetilde{A}(\bx,\bxi)=\chi(\bxi)\syma(\bx,\bxi)\in S_{1,0}^m(\Sigma\times\RR^n)$ modulo functions in $S_{1,0}^{-\infty}(\Sigma\times\RR^n)=\bigcap_{m\in\RR}S_{1,0}^m(\Sigma\times\RR^n)$, which relate to negligible pseudo-differential operators \cite{BenzoniGavage2006,Grubb1996}.

In the remainder, we will not make a distinction between $\hat{A}$ and $\tilde{A}$. In case we need to explicitly consider a singularity, we will use a pseudo-differential operator with a symbol in $S_{1,0}^{m_1,m_2}(\Sigma\times\RR^n)$.

\subsection{Eigensolutions of  Maxwell equations with constant coefficients}
\label{ssec:gene-sol}
We derive the eigensolutions of the bilateral Laplace-Fourier transformed Maxwell
equations \eqref{eq:hyperbolic-sys-1} with source term $\bm{f}=\bm{0}$. 
We introduce
\begin{equation}
    \label{eq:U}
    \cU
    =\left\{ (\bomega,s)\in\RR^2\times\CC:\, \Re(s)>0  \right\}.
\end{equation}
Here and in the following, we write $\Re(s)=\eta$ and $\Im(s)=\xi$ for the real and imaginary parts of a complex number $s\in\mathbb{C}$, and $\bomega=(\omega_2,\omega_3)$.
After multiplying with $e^{-\eta t},\,\eta>0$ and applying the bilateral Laplace-Fourier transform,
the Maxwell equations in characteristic variables $\brho$
\eqref{eq:hyperbolic-sys-1-trans} transform into the following
system of differential-algebraic equations
\begin{equation}
    \label{eq:Maxwell-form-3}
    \Lambda \partial_{x_1} \hat{\brho} + 
    \left[ s I_6 + \imath(\omega_2\widetilde{ A_2 }
    + \omega_3\widetilde{ A_3}) \right] \hat{\brho}=\bm{0}  \quad \text{for } x_1<0,
\end{equation}
with $\Lambda$ defined in \eqref{eq:matA1}.

We write
$\hat{\brho}=[\hat{\brho}_+^\top ;\hat{\brho}_{\s}^\top ;\hat{\brho}_-^\top ]^\top $,
with the standing component $\hat{\brho}_{\s}$ corresponding to
the kernel of $\Lambda$. Using \(\Re(s)>0\) in \eqref{eq:Maxwell-form-3} yields
\begin{equation}
    \label{eq:rho-o}
    \hat{\brho}_{\s}=
    \frac{1}{s}(\imath \omega_2 \EE-\imath \omega_3 \JJ)\hat{\brho}_+
    +\frac{1}{s}(\imath \omega_2 \DD-\imath \omega_3 \II)\hat{\brho}_-.
\end{equation}
Eliminating $\hat{\brho}_{\s}$ from \eqref{eq:Maxwell-form-3} we obtain
\begin{equation}
    \label{eq:Maxwell-form-4}
    \partial_{x_1}\brhohpm + \Mhat(\bomega,s) \brhohpm =\bm{0} \quad \text{for }x_1<0,
\end{equation}
with the one-homogeneous matrix symbol $\Mhat\in (S_{1,0}^1(\RR^3))^{4\times 4}$  given by
\begin{equation}
    \label{eq:matrix-M}
    \Mhat(\bomega,s)=\frac{1}{2s}\begin{pmatrix}
        (\kappa^2+ s^{2})I_2 & -\Mom \\ \Mom &
        -(\kappa^2+ s^{2})I_2
    \end{pmatrix}
\end{equation}
with
\begin{align}\label{eq:Momega}
    \Mom=\begin{pmatrix}
        \omega_2^2-\omega_3^2 & 2 \omega_2 \omega_3 \\
        2 \omega_2 \omega_3 & -(\omega_2^2-\omega_3^2) 
    \end{pmatrix}.
\end{align}
Using \(\Mom^2= |\bomega|^4 I_2\) and \(|\bomega|^2=\omega_2^2+\omega_3^2\),
it is straightforward to check that the matrix symbol \(\Mhat\) \eqref{eq:matrix-M} satisfies the eigendecomposition
\begin{equation}
    \label{eq:pdp}
    \Mhat = \Phat \Dhat \Phat^{-1},
\end{equation}
where
\begin{align}
    \label{eq:matrix-P}
    \Phat&=\frac{1}{2\sqrt{\kappa s}(\kappa+s)}
    \begin{pmatrix}
        \Mom & (\kappa+s)^2 I_2\\ 
        (\kappa+s)^2 I_2 & \Mom
    \end{pmatrix} ,\\
        \Phat^{-1}&= \frac{1}{2\sqrt{\kappa s}(\kappa+s)}
        \begin{pmatrix}
            -\Mom & (\kappa+s)^2 I_2\\ 
            (\kappa+s)^2 I_2 & -\Mom
        \end{pmatrix}, \label{eq:matrix-Pinv}
\end{align}
and
\[
    \Dhat = \begin{pmatrix}
        -\kappa I_2 & \\ & \kappa I_2
    \end{pmatrix},
\]
with $\kappa=\sqrt{|\bomega|^2+s^2}$ and the matrix symbols $\Phat,\,\Phat^{-1} \in (S_{1,0}^0(\RR^3))^{4\times 4}$ and $\Dhat \in (S_{1,0}^1(\RR^3))^{4\times 4}$ since they are, respectively, zero-homogeneous and one-homogeneous. 
Here, we take the branch of the square root that is well-defined for $\arg(z)\in(-\pi,\pi)$ for $z\in\mathbb{C}$ and such that $\sqrt{z}>0$ for   $z\in\mathbb{R}$ with $z>0$.
Note that $\Re(\kappa)\geq\Re(s)>0$ for $(\bomega,s)\in\cU$ by Lemma~\ref{lem:Rek}.
The eigendecomposition \eqref{eq:pdp} will be very useful in the next sections.
The eigensolutions of \eqref{eq:Maxwell-form-3} are stated in the following lemma.
\begin{lemma}[Eigensolutions]
    \label{lem:decaying_sol}
    For $(\bomega,s)\in\cU$ 
    all solutions to \eqref{eq:Maxwell-form-3} 
    are given by 
    \begin{equation}
        \label{eq:planewave solution-2}
          \hat{\brho}(x_1) = \Psi_L \bm{c}_0 \e^{\kappa x_1}+\Psi_R \bm{d}_0 \e^{-\kappa x_1},
    \end{equation}
    with  \(\Psi_L=\left[\Psi_1, \Psi_2\right]\in \CC^{6\times 2}, \Psi_R=\left[\Psi_3, \Psi_4\right]\in \CC^{6\times 2}\), 
    $\bm{c}_0, \bm{d}_0 \in \mathbb{C}^2$,
    and \(\Psi=[\Psi_L, \Psi_R]\) given by 
    \begin{equation}
        \label{eq:Psi}
        \Psi=\frac{1}{c_{\Psi}}
        \begin{pmatrix}
            \Mom & (\kappa+s)^2 I_2 \\
            \imath\sqrt{2}(\kappa+s) \hat{O}_L &  \imath\sqrt{2}(\kappa+s) \hat{O}_R\\
            (\kappa+s)^2 I_2 & \Mom 
        \end{pmatrix},
    \end{equation}
    with \(c_{\Psi}=\sqrt{|\kappa+s|^4+|\bomega|^4}\) and 
    \(\Ohat_L = \begin{pmatrix}
        -\omega_3 & \omega_2 \\
        -\omega_2 & -\omega_3
    \end{pmatrix}\), \(\Ohat_R = \begin{pmatrix}
        -\omega_3 & \omega_2 \\
        \omega_2 & \omega_3
    \end{pmatrix}\).
\end{lemma}

\begin{proof}
Introducing the transformed variable \(\hat{\bm{v}}=\Phat^{-1}\brhohpm\), \eqref{eq:Maxwell-form-4} becomes
\begin{equation}
    \label{eq:diag-sys}
    \partial_{x_1}\hat{\bm{v}} +  \Dhat\hat{\bm{v}}=\bm{0} \quad \text{for }x_1<0,
\end{equation}
which has general solution
\[
    \hat{\bm{v}}(x_1) = \begin{pmatrix}
        \bm{c}_0 \e^{\kappa x_1}\\
        \bm{d}_0 \e^{-\kappa x_1}
    \end{pmatrix},\qquad \bc_0,\ \bd_0\in\CC^2.
\]
The eigensolution of the form \eqref{eq:planewave solution-2} with $\Psi$ stated in \eqref{eq:Psi} can then be obtained using \(\brhohpm =\Phat\bvh\), \eqref{eq:rho-o} for $\brhoh_{\s}$, and normalization by $c_\Psi$.
\end{proof}

\begin{remark}
    The eigensolutions \eqref{eq:planewave solution-2} cover
    the zero-frequency case $\bomega=\bm{0}$; compare to \cite{Feng1999} where the case $\bomega=\bm{0}$ was treated separately.
\end{remark}

\subsection{Solution of HW-NRBCs for constant coefficients}
For constant coefficients and $\bg=\bm{0}$, the HW-NRBCs \eqref{eq:hw-nrbcs} reduce to the
following system of algebraic equations 
in the bilateral Laplace-Fourier transformed auxiliary variables $\hat{\brho}_{j}=\scrL\scrF[\brho_j]$ and $\hat{\brho}_0=\scrL\scrF[\brho_0]$, where $\brho_0(x_2,x_3,t)=S^\top \bm{u}(0,x_2,x_3,t)$ and $\brho_j(x_2,x_3,t)=S^\top \bm{\phi}_j(x_2,x_3,t)$ with $S$ given by \eqref{eq:Smat},
\begin{subequations}
    \label{eq:hw-nrbcs-lf}
    \begin{alignat}{3} 
        s\Lambda \hat{\brho}_{1}-\widehat{\mathcal{L}}_0\hat{\brho}_0
        &=\bm{0},\quad \label{eq:hw-nrbcs-lf-1}\\
        \widehat{\mathcal{L}}_j\hat{\brho}_j - \widehat{\mathcal{K}}_j\hat{\brho}_{j+1}&=\bm{0},&&\label{eq:hw-nrbcs-lf-2}\\ 
        \hat{\brho}_{J+1,-} &= \bm{0},&&  \label{eq:hw-nrbcs-lf-3} \\
        \hat{\brho}_{1,\s} &=
    \frac{1}{s}(\imath \omega_2 \EE-\imath \omega_3 \JJ)\hat{\brho}_{1,+}
    +\frac{1}{s}(\imath \omega_2 \DD-\imath \omega_3 \II)\hat{\brho}_{1,-}, &&\label{eq:hw-nrbcs-lf-4}
    \end{alignat}
\end{subequations}
for $(\bomega,s)\in\cU$ and $1\leq j\leq J$.
The operators $\widehat{\mathcal{L}}_j$ and $\widehat{\mathcal{K}}_j$ are obtained from the bilateral Laplace-Fourier transform of \eqref{eq:KjLj} and are equal to
\begin{subequations}
    \label{eq:hw-nrbcs-lf_operators}
\begin{align}
    \widehat{\mathcal{L}}_j &=\widehat{S^{\top} \mathcal{L}_j S} = (a_j \Lambda -I_6)s
    -\imath \omega_2\widetilde{A}_2 -\imath \omega_3\widetilde{A}_{3} ,\\
    \widehat{\mathcal{K}}_j &=\widehat{S^{\top} \mathcal{K}_j S} = (a_j \Lambda +I_6)s
    +\imath \omega_2\widetilde{A}_2 +\imath \omega_3\widetilde{A}_{3} 
\end{align}
\end{subequations}
for $j=0,1,\ldots,J$.
We note that \eqref{eq:hw-nrbcs-lf-4} results from \eqref{eq:rho2-max-3-4}, which is necessary because \eqref{eq:hw-nrbcs-lf-1} determines $\hat{\brho}_1$ only up to the standing component, see Remark~\ref{rem:nonuniqueness} below.
In the following lemma, we state that the auxiliary variables in the bilateral Laplace-Fourier transformed HW-NRBCs \eqref{eq:hw-nrbcs-lf} are linear combinations of $\Psi_L$ and $\Psi_R$ given in Lemma~\ref{lem:decaying_sol},
which allows us to compute the reflection coefficient of the HW-NRBCs in Lemma \ref{lem:reflection} below.
\begin{lemma}
    \label{lem:coeff_aux_variables}
    Suppose $(\bomega,s)\in\cU$,
    $0 < a_0 \leq 1$, $0<a_j<1$ for $1 \leq j \leq J$.
    Let $\hat{\brho}_j$, $0\leq j\leq J+1$, satisfy \eqref{eq:hw-nrbcs-lf}.
    Then there exist \(\bm{c}_j(\bomega,s),\, \bm{d}_j(\bomega,s) \in \mathbb{C}^2\) such that
    \begin{equation}
        \label{eq:rhoLR}
        \hat{\brho}_j=\Psi_L \bm{c}_j+\Psi_R \bm{d}_j\quad \text{for all } 0\leq j\leq J+1,
    \end{equation}
    with
    \begin{equation}
        \label{eq:cjdj}
            \bm{c}_j = \frac{\gamma^+_0}{s}\prod_{l=1}^{j-1} \frac{\gamma_l^{+}}{\gamma_l^{-}} \bm{c}_0 , \qquad
        \bm{d}_j = \frac{\gamma^-_0}{s}\prod_{l=1}^{j-1} \frac{\gamma_l^{-}}{\gamma_l^{+}}\bm{d}_0, \quad 1\leq j\leq J+1,
    \end{equation}
    where $\gamma_j^{ \pm}=a_j s \pm \kappa$ and \(\prod_{i=1}^{J-1}\cdot =1\) if \(J=1\). 
\end{lemma}
\begin{proof}
Denote the projection onto the kernel of \(\Lambda\) with
$\Pi=\operatorname{diag}(0,0,1,1,0,0)$.
From Lemma \ref{lem:decaying_sol}, see also \eqref{eq:rho-o}, the columns of $\Psi_L$ and $\Psi_R$ form a basis for the space of all vectors $\hat{\bm{w}}\in \mathbb{C}^6$ satisfying $\Pi \widehat{\mathcal{M}} \hat{\bm{w}}=\bm{0}$, with 
$\widehat{\mathcal{M}}=s I_6+\imath (\omega_2 \widetilde{A}_2+\omega_3\widetilde{A}_3)$.
We observe the following relations
\begin{align}\label{eq:relations_M_Lj_Kj}
    \widehat{\mathcal{M}} \Psi_L=-\kappa\Lambda\Psi_L,\quad 
    \widehat{\mathcal{M}} \Psi_R=\kappa\Lambda\Psi_R,\quad
    \widehat{\mathcal{L}}_j =a_js\Lambda-\widehat{\mathcal{M}},\quad
    \widehat{\mathcal{K}}_j =a_js\Lambda+\widehat{\mathcal{M}},
\end{align}
which follow from the definition of the eigensolutions $\Psi_L$ and $\Psi_R$, \eqref{eq:Maxwell-form-3}, and the definition of $\widehat{\mathcal{M}}$.
In view of \eqref{eq:hw-nrbcs-lf-1} and since
 $\Pi\widehat{\mathcal{M}}=-\Pi \widehat{\mathcal{L}}_j= \Pi \widehat{\mathcal{K}}_j$,
 we observe that $\Pi\widehat{\mathcal{M}}\hat{\brho}_0=\bm{0}$.
Noticing that \eqref{eq:hw-nrbcs-lf-4} just states that $\Pi\widehat{\mathcal{M}}\hat{\brho}_1=\bm{0}$, we deduce recursively from \eqref{eq:hw-nrbcs-lf-2} that $\Pi\widehat{\mathcal{M}}\hat{\brho}_j=\bm{0}$ for all $0\leq j\leq J+1$.
Thus, there exist $\bm{c}_j(\bomega,s), \bm{d}_j(\bomega,s) \in \mathbb{C}^2$ such that
    $\hat{\brho}_j=\Psi_L \bm{c}_j+\Psi_R \bm{d}_j$, for all $0 \leq j \leq J+1$, which proves \eqref{eq:rhoLR}.\\
    From \eqref{eq:hw-nrbcs-lf-1}, \eqref{eq:hw-nrbcs-lf-2} and \eqref{eq:relations_M_Lj_Kj} we deduce that
    \begin{subequations}\label{eq:coeff_relation}
    \begin{align}
        s\Lambda(\Psi_L \bm{c}_1+\Psi_R \bm{d}_1) &= \Lambda\left( \gamma_0^+ \Psi_L {\bm{c}}_0 + \gamma_0^-\Psi_R \bm{d}_0\right),\\
        \Lambda \left(\gamma_j^+ \Psi_L\bm{c}_j+\gamma_j^-\Psi_R\bm{d}_{j}\right)& =\Lambda \left(\gamma_j^- \Psi_L \bm{c}_{j+1} +\gamma_j^+\Psi_R\bm{d}_{j+1}\right), 
    \end{align}
    \end{subequations}
    for $ 1 \leq j \leq J $. 
    Using that $\Mom^2=|\bomega|^4 I_2$ and $|\bomega|^2=\kappa^2-s^2$, we see that
    \begin{align*}
        \det\begin{pmatrix}
            \Mom & (\kappa+s)^2 I_2 \\
            -(\kappa+s)^2 I_2 & -\Mom
        \end{pmatrix}= 16 (s\kappa)^2(\kappa+s)^4\neq 0,
    \end{align*}
 because $\Re(\kappa)\geq\Re(s)>0$ by Lemma~\ref{lem:Rek}.
 The columns of $\Lambda\Psi_L$ and $\Lambda\Psi_R$ are thus linearly independent, and \eqref{eq:cjdj} is proven by comparing the coefficients in \eqref{eq:coeff_relation}.
\end{proof}

\begin{remark}
    Note that for $(\bomega,s)\in\cU$ we have $\gamma_j^{ \pm} \neq 0$ for $1 \leq j \leq J$ if $0<a_j<1$. This result follows from Lemma \ref{lem:Rek}, which states that \(\Re (\kappa)\geq \Re (s)\). Hence,
    $$
    \Re{(\gamma_j^+)}=a_j \Re{(s)}+\Re(\kappa) \geq a_j \Re(s)+\Re(s)=\left(1+a_j\right) \Re(s)>0 ,
    $$
    which gives \(\gamma_j^+\neq 0\). 
    For $0<a_j < 1$, we have
    $$
    \Re{(\gamma_j^-)}=a_j \Re(s)-\Re(\kappa) \leq a_j \Re(\kappa)-\Re(\kappa)=\left(a_j-1\right) \Re(\kappa)<0,
    $$
    which gives \(\gamma_j^-\neq 0\).
    In addition, if $a_j=1$, then $\gamma^-_j=0$ if and only if $\bomega=\bm{0}$.
    In that case, \eqref{eq:coeff_relation} shows that $\bm{c}_j=\bm{d}_{j+1}=\bm{0}$ for $1\leq j\leq J$.
\end{remark}
\begin{remark}
    Lemma~\ref{lem:coeff_aux_variables}   holds for an arbitrary set of vectors $\byh_j\in\mathbb{C}^6$, $0\leq j\leq J+1$ that satisfies the set of equations \eqref{eq:hw-nrbcs-lf}. In particular, $\byh_0$ does not have to be the trace of $\brhoh_0$.
\end{remark}
\begin{remark}\label{rem:nonuniqueness}
Suppose $\bvh\in\mathbb{C}^6$ is such that $\bvh_-=\bvh_+=\bm{0}$ and $\bvh_\s$ is arbitrary. If $\brhoh_j$, $0\leq j\leq J+1$ are a solution to \eqref{eq:hw-nrbcs-lf-1}--\eqref{eq:hw-nrbcs-lf-3}, then $\hat{\bm{\psi}}_j$, $0\leq j\leq J+1$, defined by
$\hat{\bm{\psi}}_0=\brhoh_0$ and $\hat{\bm{\psi}}_j = \brhoh_j+(-1)^j\hat{\bm{v}}$, $j\geq 1$, also solve \eqref{eq:hw-nrbcs-lf-1}--\eqref{eq:hw-nrbcs-lf-3}.
This follows from $\widehat{\mathcal{L}}_j\bvh=-\widehat{\mathcal{K}}_j\bvh$ and $(\brhoh_{j}+(-1)^j\bvh)_{\pm}=\brhoh_{j,\pm}$.
Therefore, without condition \eqref{eq:hw-nrbcs-lf-4},  the boundary conditions \eqref{eq:hw-nrbcs-lf-1}--\eqref{eq:hw-nrbcs-lf-3} have multiple solutions.
\end{remark}
\subsection{Reflection coefficients of HW-NRBCs}
The interaction of the solution with the boundary is often described by the reflection coefficient, see \cite{Higdon1994,Givoli2003} for example.
The reflection coefficient derived in the next result is instrumental in obtaining a Fourier-Laplace domain description of the HW-NRBC. Recall the definition of $\Mom$ given in \eqref{eq:Momega}.
\begin{lemma}[Reflection coefficients]
    \label{lem:reflection}
    Given the conditions stated Lemma~\ref{lem:coeff_aux_variables} we have   $\bm{c}_0 =\mr{R} \bm{d}_0$ with $\mr{R}\in\mathbb{C}^{2\times 2}$ given by
    \begin{equation}
        \label{eq:reflection-coef}
         \mr{R}=-\frac{\ga}{(\kappa+s)^2}\Mom\quad\text{with}\quad\ga=\frac{\gamma_0^{-}}{\gamma_0^{+}} \prod_{j=1}^J\left(\frac{\gamma_j^{-}}{\gamma_j^{+}}\right)^2.
    \end{equation}
\end{lemma}
\begin{proof}
    Lemma \ref{lem:coeff_aux_variables} shows that
    \begin{align*}
     \hat{\brho}_{J+1}&=\Psi_L \bm{c}_{J+1}+\Psi_R \bm{d}_{J+1} \quad\text{with} \\
     \bm{c}_{J+1} &= \frac{\gamma_0^{+}}{s} \prod_{l=1}^J \frac{\gamma_l^{+}}{\gamma_l^{-}} \bm{c}_0, \quad\text{and}\quad
\bm{d}_{J+1}=\frac{\gamma_0^{-}}{s} \prod_{l=1}^J \frac{\gamma_l^{-}}{\gamma_l^{+}}\bm{d}_0.
    \end{align*}
    The truncation condition \eqref{eq:hw-nrbcs-lf-3}, $\hat{\brho}_{J+1,-}=\bm{0}$, can be written as
    $\Pi_{-} \hat{\brho}_{J+1}=\bm{0}$, with  $\Pi_{-}=\operatorname{diag}(0,0,0,0,1,1)$.
    Since
    $\Pi_{-} \Psi_L=\frac{1}{c_\Psi}(\kappa+s)^2 I_2$ and  $\Pi_{-} \Psi_R=\frac{1}{c_\Psi} \Mom$,
    the truncation condition becomes after multiplication by \(s c_\Psi\)
    \[
    (\kappa+s)^{2} \gamma_{0}^{+} \prod_{l=1}^{J} \frac{\gamma_{l}^{+}}{\gamma_{l}^{-}} \bm{c}_0+\gamma_{0}^{-} \prod_{l=1}^{J} \frac{\gamma_{l}^{-}}{\gamma_{l}^{+}} \Mom \bm{d}_0=\bm{0}.
    \]
    Hence \(\bm{c}_{0}=\mathrm{R} \bm{d}_{0}\) with $\mathrm{R}$ as claimed.
\end{proof}

\begin{remark}\label{rmk:plainwave} 
In view of Lemma \ref{lem:decaying_sol}, we have for $s=-\imath \omega_f,\, \omega_f\in \RR$,
\begin{align*}
    \brho(\bx,t)= &A_L(\bx) \e^{\imath(\bk_L\cdot \bx-\omega_f t)}\Psi_L\bc_0+A_R(\bx) \e^{\imath(\bk_R\cdot \bx-\omega_f t)}\Psi_R\bd_0 ,
\end{align*}
with amplitude functions $A_L(\bx)=e^{\Re(\kappa)x_1}$, $A_R(\bx)=e^{-\Re(\kappa)x_1}$, and wave vectors $\bk_L=(\Im(\kappa),\omega_2,\omega_3)^\top$, $\bk_R=(-\Im(\kappa),\omega_2,\omega_3)^\top$.
By Lemma \ref{lem:Rek}, we have $\Im(s)\Im(\kappa)\geq 0$, which shows that $\brho(\bx,t)$ is a left-moving wave if $\bd_0=\bm{0}$, while $\brho(\bx,t)$ is a right-moving wave if $\bc_0=\bm{0}$. Since $\kappa=\sqrt{|\bomega|^2-\omega_f^2}$, we observe that the contribution of $\Psi_L$ is an evanescent plane wave if  $|\bomega|^2>\omega_f^2$. If  $|\bomega|^2<\omega_f^2$, then $\kappa$ is purely imaginary and $\brho(\bx,t)$ is a propagating plane wave, and 
$\alpha=\kappa/s=\sqrt{1-|\bomega|^2/\omega_f^2}$ is a real number. Lemma \ref{lem:reflection} then implies that the reflection coefficient of the HW-NRBCs for the Maxwell equations is
    \begin{equation*}
        \mathtt{R}= \frac{1}{\omega_f^2 (1+\alpha)^2}\frac{a_0-\alpha}{a_0+\alpha}  
        \left( \prod_{j=1}^{J}\frac{a_j-\alpha}{a_j+\alpha} \right)^2 \Mom
        \quad \text{for }J\geq 0.
    \end{equation*}
Since the cosine of the angle between the wave vector of the plane wave and the $x_1$-axis is $\alpha$, we see that the plane wave is annihilated at $\Gamma$ if there is $j$ such that $\alpha=a_j$. 
Note, however, that for $\alpha=0$ (waves parallel to the boundary $\Gamma$), $\mathtt{R}=\frac{1}{|\bomega|^2} \Mom$, which gives $|R|=1$. Waves nearly parallel to $\Gamma$ require therefore larger values of $J$ to obtain a small reflection coefficient $\mathtt{R}$.
\end{remark}

\section{On Kreiss' condition for stability}\label{sec:ukc}
As pointed out by Agmon, see \cite{Kreiss1970,Motamed2019}, a necessary condition for the well-posedness of the Maxwell equations with the HW-NRBCs on the half-space is that
solutions of the form 
$\brho(\bx,t) =  \e^{\kappa x_1} 
        \e^{\imath \bomega\cdot \bar{\bm{x}} +st} \Psi_L \bm{c}_0$
must vanish, i.e., $\bm{c}_0=\bm{0}$.
As Lemma~\ref{lem:reflection} shows, this necessary condition for well-posedness is satisfied.
A sufficient condition for (strong) well-posedness is the uniform Kreiss condition, see \cite{Kreiss1970,Majda1975}.
To formulate the  Kreiss condition in our context and to investigate its validity, we next derive the mapping that takes outgoing components to ingoing components of the solution of the Maxwell equations with HW-NRBCs.
This mapping is also instrumental in the following sections.
\begin{lemma}\label{lem:operator-S}
    Let $(\bomega,s)\in\cU$, \(0<a_0 \leq 1\) and \(0<a_j<1\), \(j=1,\dots,J\).
    Let $\hat{\brho}_j$, $0\leq j\leq J+1$ satisfy \eqref{eq:hw-nrbcs-lf}.
     Then 
    \begin{equation}
    \hat{\brho}_{0,-}=\Sop \hat{\brho}_{0,+}
    \quad \text{with}\quad 
    \Sop=\frac{1}{\beta_J} \Mom \text{ and } \beta_J=\frac{(\kappa+s)^2-(\kappa-s)^2\ga}{1-\ga} ,
    \label{eq:rho-m-rho-p} 
     \end{equation}
    where $\ga$ is defined in \eqref{eq:reflection-coef}, and $\Mom$ in \eqref{eq:Momega}.
\end{lemma}
\begin{proof}
    An application of Lemma \ref{lem:reflection} shows that
    \[
        \hat{\brho}_{0}=\Psi_{L} \mr{R}\bm{d}_{0}+\Psi_{R} \bm{d}_{0},
    \]
    where $\mr{R}$ is given in \eqref{eq:reflection-coef}.
    The above relation implies that
    \[
        \hat{\brho}_{0,+}=\left( \Mom \mr{R} +(\kappa+s)^2 I_2 \right) \bm{d}_0
        \quad \text{and}\quad
        \hat{\brho}_{0,-}=\left(  (\kappa+s)^2 \mr{R} + \Mom\right) \bm{d}_0,
    \]
    which leads to
    \begin{align*}
        \hat{\brho}_{0,-} 
        &= \left(  (\kappa+s)^2  \mr{R} + \Mom\right) \left( \Mom \mr{R} +(\kappa+s)^2 I_2 \right)^{-1}\hat{\brho}_{0,+} \\
        &= \frac{1-\ga}{(\kappa+s)^2 \left( 1-\frac{(\kappa-s)^2}{(\kappa+s)^2}\ga \right)}\Mom\hat{\brho}_{0,+}.
    \end{align*}
    Here, we have used that $\Mom^2 = |\bomega|^4 I_2 = (\kappa+s)^2 (\kappa-s)^2 I_2$.
    Rearranging terms, we arrive at the claim.
\end{proof}

The uniform Kreiss condition requires that 
\begin{align*}
    \inf_{(\bomega,s)\in\cU} |\det \NJ(\bomega,s)|>0,
\end{align*}
where the matrix symbol $\NJ(\bomega,s)$ is defined by
$$
\NJ(\bomega,s)=\Sop\Psi_{L,+}-\Psi_{L,-} \in \mathbb{C}^{2 \times 2}.
$$
The discussion preceding Lemma~\ref{lem:operator-S} shows that $|\det \NJ(\bomega,s)|>0$ for all $(\bomega,s)\in\cU$.
In order to investigate the behaviour of $\NJ(\bomega,s)$ in the limiting case, it is convenient to introduce scaled variables,
\begin{align}\label{eq:scaled_variables}
    s^\prime =\frac{s}{|\bzeta|},\quad
    \bomega'=\frac{\bomega}{|\bzeta|},
\end{align}
where \(\bzeta=(\omega_2,\omega_3,\eta,\xi)\) and
\(|\bzeta|=\sqrt{\eta^2+\xi^2+|\bomega|^2}\) for  $(\bomega,s)\in\cU$. 
The matrix symbol \(\NJ\) is homogeneous of degree zero,
that is $\NJ(\bomega,s)=\NJ(\bomega',s')$.
The limiting cases, can thus be investigated by considering $\Re(s')\to 0$. We provide the representation of $\NJ$ in the next result.
\begin{lemma}
    \label{lem:NJ}
    For $(\bomega,s)\in\cU$ we have that
    \begin{align}\label{eq:NJ}
    \NJ(\bomega,s)=-\frac{4 \kappa s}{ c_\Psi \left[ 1-\frac{(\kappa-s)^2}{(\kappa+s)^2}\ga \right]} I_{2}.
    \end{align}
\end{lemma}
\begin{proof}
    Noticing
    $
    \Psi_{L,-}=\frac{(\kappa+s)^{2}}{c_\Psi} I_{2}$,
    $\Psi_{L,+}=\frac{1}{c_\Psi} \Mom$
    and \(\Mom^{2}=|\bomega|^{4} I_{2}\), we obtain that
    \[
    \NJ(\bomega,s)=\frac{1}{c_\Psi}
    \left(\frac{|\bomega|^{4}(1-\ga) }{ (\kappa+s)^2 \left( 1-\frac{(\kappa-s)^2}{(\kappa+s)^2}\ga \right) }-(\kappa+s)^{2}\right) I_{2}.
    \]
    Using $|\bomega|^{2}=(\kappa-s)(\kappa+s)$ proves the claim by a direct calculation.
\end{proof}
We note that no division by zero occurs in \eqref{eq:NJ} for $(\bomega,s)\in \cU$, because Lemma~\ref{lem:bound_as-k} ensures that
$\left| \frac{(\kappa-s)^2}{(\kappa+s)^2}\ga\right|$
is strictly smaller than one for $a_j\Re(s)>0$, $0\leq j\leq J$. 
If
$\Re(s')\to 0$, then a division by zero may occur in \eqref{eq:NJ}. 
The following asymptotic expansions of $\ga$ allow a more precise investigation.

\begin{lemma}\label{lem:kappa_asymptotics}
    Assume that $a_0=1, \, 0<a_j<1$ for $j=1,..., J$. The following asymptotics hold for  $|\kappa/s|\ll 1$, with $|\kappa/s| < \min_{0\leq j\leq J} a_j$,
\begin{alignat*}{3}
\ga &= 1-c_0(\kappa/s)\frac{\kappa}{s} + O(\left|\frac{\kappa}{s}\right|^2),                            \quad&&\text{with } c_0(z)=\sum_{j=0}^J \frac{2(2-\delta_{0j})}{a_j +z}, \quad &&z\in\CC,\\
\frac{\kappa-s}{\kappa+s} \ga &= -1+c_1(\kappa/s)\frac{\kappa}{s}+ O(\left|\frac{\kappa}{s}\right|^2),\quad&&\text{with } c_1(z)=\sum_{j=0}^J \frac{4}{a_j +z},\quad &&z\in\CC,\\
\frac{(\kappa-s)^2}{(\kappa+s)^2}\ga &= 1-c_2(\kappa/s)\frac{\kappa}{s} +O(\left|\frac{\kappa}{s}\right|^2),\quad &&\text{with } c_2(z)=\sum_{j=0}^J \frac{2(2+\delta_{0j})}{a_j +z},\quad &&z\in\CC.
\end{alignat*}
\end{lemma}
\begin{proof}Follows directly from 
    $\frac{as-\kappa}{as+\kappa}=1-\frac{2}{a+\kappa/s}\kappa/s$
    and the definition of $\ga$.
\end{proof}
\begin{lemma}\label{lem:s_asymptotics}
    Assume that $a_0=1, \, 0<a_j<1$ for $j=1,..., J$. The following asymptotics hold for $|s/\kappa|\ll 1$
\begin{align*}
\ga &= -1+d_0(s/\kappa)\frac{s}{\kappa} + O(\left|\frac{s}{\kappa}\right|^2),                            \quad\text{with } d_0(z)=\sum_{j=0}^J \frac{2(2-\delta_{0j})a_j}{a_j z+1},\ &&z\in\CC,\\
\frac{\kappa-s}{\kappa+s} \ga &= -1+d_1(s/\kappa) \frac{s}{\kappa}+ O(\left|\frac{s}{\kappa}\right|^2),\quad\text{with } d_1(z)=\sum_{j=0}^J \frac{4a_j}{a_j z+1}, &&z\in\CC,\\
\frac{(\kappa-s)^2}{(\kappa+s)^2} \ga &= -1+d_2(s/\kappa)\frac{s}{\kappa} + O(\left|\frac{s}{\kappa}\right|^2),\quad\text{with } d_2(z)=\sum_{j=0}^J \frac{2(2+\delta_{0j})a_j}{a_j z+1}, &&z\in\CC.
\end{align*}
\end{lemma}
\begin{remark} 
    The assumption $a_0=1$ may be replaced by $0<a_0\leq 1$ in Lemmas~\ref{lem:kappa_asymptotics} and~\ref{lem:s_asymptotics} with slight modifications of $c_0(z)$, $c_1(z)$, and $c_2(z)$, etc.
\end{remark}
In view of the previous lemmas, we have that
\begin{alignat*}{3}
    \ga&=1-c_0(\kappa/s)\kappa/s + O(|\kappa/s|^2)\qquad&&\text{for}\quad\kappa\to 0,\\
    \ga&=-1+d_0(s/\kappa)s/\kappa + O(|s/\kappa|^2)&&\text{for}\quad s\to 0,
\end{alignat*}
Let $\{(\bomega_k,s_k)\}_k$ be a sequence in $\cU$.
Then for the scaled variables defined in \eqref{eq:scaled_variables}, we have that
\begin{align*}
    \lim_{k\to\infty} s_k'\to 0\quad&\implies\quad \lim_{k\to\infty}\kappa_k' =1,\\
    \lim_{k\to\infty} s_k'\to \pm \imath\frac{1}{\sqrt{2}} \quad&\iff\quad \lim_{k\to\infty}\kappa_k' =0.
\end{align*}
Using Lemma \ref{lem:kappa_asymptotics} we obtain from \eqref{eq:NJ}
\begin{align*}
    \lim_{\kappa'\to 0} |\det \NJ(\bomega',s')|\neq 0,
\end{align*}
but Lemma \ref{lem:s_asymptotics} gives
\begin{align*}
    \lim_{s'\to 0} |\det \NJ(\bomega',s')|= 0.
\end{align*}
Therefore, the uniform Kreiss condition does not hold. Using the terminology of \cite{Kreiss1970,Motamed2019}, 
the points $(\bomega_0',0)$, where $\bomega_0'\in\mathbb{R}^2$ with $|\bomega_0'|=1$ is arbitrary, are called generalized eigenvalues. Since the corresponding $\kappa_0'=1$, these generalized eigenvalues are also called surface generalized eigenvalues, see \cite[Definition 3.1]{Motamed2019}.
In \cite[Theorem~3.1]{Motamed2019}, it is shown that strictly hyperbolic systems are well-posed in the generalized sense in this case.
The Maxwell equations with HW-NRBCs are, however, not strictly hyperbolic and contain characteristic boundaries, so the results of \cite{Motamed2019} do not apply.
In the next section, we construct suitable symmetrizers for our problem, which allow us to prove well-posedness in the generalized sense of the Maxwell equations with HW-NRBCs
which are not strictly hyperbolic and contain characteristic boundaries.

\section{Symmetrizer for the Maxwell equations with Hagstrom-Warburton nonreflecting boundary conditions}
\label{sec:symmetrizer}
The symmetrizer $\cR\in (\OpsS_{1,0}^{m_1,m_2}(V\times\RR^3))^{\ell\times \ell}$, $m_1,m_2 \in \RR$, $\ell\in \mathbb{N}$ for the Maxwell equations with HW-NRBCs on the space-time domain $V=\Omega\times(0,\infty)\subset \RR^{4}$ or $V=\Gamma\times(0,\infty)\subset \RR^3$ is a zero homogeneous pseudo-differential operator and is defined as
\begin{equation*}
    \begin{split}
     &   \cR(\bx,t,\imath D_{\bxb}, \imath D_t+\eta)
    (e^{-\eta t}\bv(\bx,t))  \\
    &=\frac{1}{(2\pi)^3} \int_{\RR^3} e^{\imath\xi t+\imath \bomega\cdot \bxb} \hat{\cR}(\bomega,\imath\xi+\eta) 
    \widehat{\bv} (x_1,\bomega,\imath\xi+\eta) \,\mathrm{d}\bomega \,\mathrm{d}\xi,\\
    &=(\cL\cF)_3^{-1}[\cRh(\bomega,s)\bvh(x_1,\bomega,s)](\bx,t;\eta),\quad \forall \bv(\bx,t) \in (C_{(0)}^\infty(\overline{V}))^\ell,
    \end{split}
\end{equation*}
with matrix symbol $\hat{\cR}\in (S_{1,0}^{m_1,m_2}(\RR^3))^{\ell\times \ell}$ and $\hat{\bv} (x_1,\bomega,s)
=\int_{\RR^3} e^{-\imath \xi t-\imath\bomega\cdot\bxb-\eta t}\bv(\bx,t)\,\mathrm{d}\bxb \,\mathrm{d}t $, $\bxb=(x_2,...,x_n)^\top$, $\bomega\in \RR^2$, where
 $\bv(\bx,t)=\bm{0}$ for $(\bx,t)\in \RR^4\setminus \overline{V}$.
The main challenge in proving stability is to obtain a suitable symmetrizer symbol $\cRh$, which is the main topic of this section.

\begin{figure}[t]
    \centering
    \begin{tikzpicture}
        \filldraw[fill=black!5, draw=black] (0,-1.5) arc[start angle=-90, end angle=90, radius=1.5];
        \filldraw[fill=red!10, draw=black] (0,-0.25) arc[start angle=-90, end angle=90, radius=0.25];
        \filldraw[fill=blue!10, draw=black](0,-1.25) arc[start angle=-90, end angle=90, radius=0.25];
        \filldraw[fill=blue!10, draw=black] (0,0.75) arc[start angle=-90, end angle=90, radius=0.25];
        \filldraw[red] (0,0) circle (0.03)node[left, xshift=-5mm]{$s'=0$}; 
        \filldraw[fill=blue, draw=black] (0,1.) circle (0.03)node[left]{$\xi'=|\bomega'|$}; 
        \filldraw[fill=blue, draw=black] (0,-1.) circle (0.03) node[left]{$\xi'=-|\bomega'|$}; 
         \draw[-,dashed] (0,0) -- (1.5,0);
        \draw[->] (-0.5,0) -- (3,0) node[below] {$\eta'$};

        \draw[-] (0,-2.2) -- (0,-1.5);
        \draw[->] (0,-1.5) -- (0,2.2) node[left] {$\xi'$};
        
    \end{tikzpicture}
    \caption{Sketch of the sets defined in \eqref{eq:cases} for the construction of the local symmetrizers:  
    $\cU_2$: $\kappa'\leq\nu$ centred around \(\imath\xi'=\pm \imath|\bomega'|=\pm \frac{\imath}{\sqrt{2}}\) in blue;
    $\cU_3$: \(|s'|\leq \nu\) centred around $s'=0$ in red. 
    }
    \label{fig:points}
\end{figure}
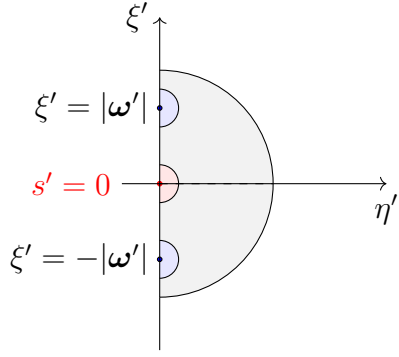

The construction of a symmetrizer \(\cRh \) for the matrix symbol \(\Mhat(\bomega,s)\) given in \eqref{eq:matrix-M} strongly depends on the transformation of \(\Mhat(\bomega,s)\)
into simple block forms, see \cite[Lemmas 2.3 and 
2.4]{Kreiss1970}, \cite[Assumption 1.9]{Majda1975}.
Since the matrix symbol \(\Mhat(\bomega, s)\) \eqref{eq:matrix-M} has a singularity at \(s=0\) and is also not diagonalizable for values of \( (\bomega,s)\in\cU\) for which $\kappa=0$, the construction of the symmetrizer will be done using the following sets
\begin{subequations}
    \label{eq:cases}
    \begin{align}
        \cU_1&=\left\{ (\bomega,s)\in \cU: |s'|\geq \nu/2,\ |\kappa'|\geq \nu/2\right\}, \label{eq:U1}\\ 
        \cU_2 &=\left\{ (\bomega,s)\in \cU:\  |\kappa'|<\nu\right\}, \label{eq:U2}\\
        \cU_3 &=\left\{ (\bomega,s)\in \cU:\  |s'|<\nu \right\} \label{eq:U3},
    \end{align}
\end{subequations}
where $0<\nu<1/(2\sqrt{2})$, as illustrated in Figure~\ref{fig:points}. Note that we have used scaled variables in the conditions that define the different sets \eqref{eq:cases}.

In the following, we construct three local symmetrizers corresponding to the three sets \eqref{eq:cases} and then combine them using a partition of unity subordinate to $\{\cU_j:1\leq j\leq 3\}$.
For $(\bomega,s)\in\cU_1$ the canonical transformation of $\Mhat(\bomega,s)$ into block form is the spectral decomposition \eqref{eq:pdp}. Since the matrix symbol \(\Phat^{-1}\) is uniformly bounded on $\cU_1$, this allows us to follow the approach taken by Kreiss \cite{Kreiss1970}.
For $(\bomega,s)\in\cU_2$, the matrix symbol \(\Phat^{-1}\) becomes singular since \(\kappa\) can approach $0$, which results in a double pair of zero eigenvalues for the matrix symbol $\Mhat(\bomega,s)$ and the matrix symbol \(\Phat^{-1}\) has rank two. In order to obtain the proper block form, we apply a Taylor series expansion of \(\Mhat(\bomega,s)\), which yields an approximate Jordan decomposition. 
For $(\bomega,s)\in\cU_3$, the matrix symbol \(\Mhat(\bomega,s)\) is singular for \(s=0\), and the construction of the symmetrizer requires special care, as described in \cite{Majda1975}.

We now summarize the key conditions defining the symbols of the symmetrizer $\cRh$ on $\cU_1\cup\cU_2$; the situation for $\cU_3$ is different, and the corresponding properties will be stated in the lemmas in Section~\ref{subs:case3-new}.
Recall the operator \(\Sop\) introduced in \eqref{eq:rho-m-rho-p}, which specifies the relation between the in- and out-going components of the solution to the Maxwell equations at the boundary.
Let us denote by $\bfgh$ the Laplace-Fourier transform of the boundary data $S^\top\bg$ in \eqref{eq:hw-nrbcs}.
Since we assumed that $\bfgh_0=\bfgh_+=\bm{0}$, the considerations that led to Lemma~\ref{lem:operator-S} yield verbatim
\begin{align}
\label{eq:hw-nrbc-g}
    \brhoh_{0,-}=\Sop\brhoh_{0,+}-\bfgh_-.
\end{align}
Recall the definition of $\bzeta=(\omega_2,\omega_3, \eta,\xi)$ in \eqref{eq:scaled_variables}, where $\eta=\Re(s)$, $\xi=\Im(s)$ and $\bomega=(\omega_2,\omega_3)$.
For notational convenience, we partition a vector $\byh\in\CC^4$ as $\byh=[\byh_+^\top,\byh_-^\top] = [\byh_1^\top,\byh_2^\top]^\top$. We use $\pm$ indexing if we relate to out- and in-going components and Arabic indexing otherwise. We denote the entries of $\byh_{j}$ by $\hat y_{jk}$, $j,k=1,2$.
\begin{definition}[Symmetrizer]
    \label{def:symmetrizer}
    Let $\mathcal{V}\subset\cU$.
    A matrix symbol \(\cRh:\cU\to (S_{1,0}^0(\RR^3))^{4\times 4}\) is called the symbol of  a symmetrizer $\cR$ for the matrix symbol $\Mhat(\bomega,s)\in (S_{1,0}^{1}(\RR^3))^{4\times 4}$ on $\mathcal{V}$ if there exists $\eta_0>0$ such that the following properties hold:
    \begin{enumerate}[label=(\roman*),ref=(\roman*),wide] 
        \item $\cRh$ is homogeneous of degree zero and a uniformly bounded smooth function of \(\bzeta\) for $\eta\geq\eta_0$ and $(\bomega,s)\in\mathcal{V}$. Moreover, $\cRh(\bomega,s)$ is Hermitian. \label{enum:symm-1}
        \item There exist $\delta_1>0$ and $C>0$ such that \(\byh^{\mathrm{H}} \cRh\byh\geq \delta_1 |\byh|^2 - C|\bfgh_-|^2\) for all $(\bomega,s)\in\mathcal{V}$ and
         \(\byh\in \mathbb{C}^4\), $\bfgh_-\in\CC^2$ satisfying \(\byh_-=\Sop\byh_+ -\bfgh_-\);\label{enum:symm-2}
        \item There exists $\delta_2>0$ such that \(\Re ((\cRh \Mhat)(\bomega,s)) \geq \delta_2 \eta I_4\) for all $(\bomega,s) \in\mathcal{V}$. \label{enum:symm-3}
    \end{enumerate}
\end{definition}

\subsection{Symmetrizer symbol for \texorpdfstring{$\Mhat$}{} on \texorpdfstring{$\cU_1$}{U1}}
\label{subs:caseI}
We employ the matrix symbol $\Phat^{-1}$ arising from the spectral decomposition \eqref{eq:pdp} and proceed according to the method described in \cite[Section 4]{Kreiss1970}. For a positive parameter $d\in \RR,\, d\geq 1$, whose value will be specified later, we define the symmetrizer symbol $\cRh_{1}$ as follows:
\begin{equation}
    \label{eq:symmetrizer-P}
    \cRh_{1} = \Phat^{-{\mathrm{H}}} Q_d \Phat^{-1} \quad \text{with}\quad Q_d=\begin{pmatrix}
        - I_2 & \\ & d I_2
    \end{pmatrix}.
\end{equation}
Since $\cRh_{1}$ is zero-homogeneous we have $\cRh_{1}\in  (S_{1,0}^0(\RR^3))^{4\times 4}$.
It is important to determine the spectral properties of $\Phat^{-1}$, which we state next.
Recall the definition of $\cU_1$ with parameter $\nu>0$ stated in \eqref{eq:U1}.
\begin{lemma}\label{lem:eigsPHP}
    Given $\Phat\in  (S_{1,0}^0(\RR^3))^{4\times 4}$, stated in \eqref{eq:matrix-Pinv}.
  For $(\bomega,s)\in\cU$, the eigenvalues of $\Phat^{-\rH}\Phat^{-1}$ and the eigenvalues of $\Phat^{\rH} \Phat$ are the same, and they are given by $\hat{\lambda}_1 =|s|/|\kappa|$ and $\hat{\lambda}_2=|\kappa|/|s|$. Moreover,
  \begin{align*}
      \frac{\sqrt{\nu}}{\sqrt{2}}|\byh|\leq  |\Phat^{-1}\byh|\leq \frac{\sqrt{2}}{\sqrt{\nu}}|\byh|,\quad\text{for all } \byh \in \mathbb{C}^4\text{ and } (\bomega,s)\in\cU_1.
  \end{align*}
  where $|\cdot|$ denotes the spectral norm induced by the Euclidian norm.
\end{lemma}
\begin{proof}
 The matrix symbol $\Mom$ defined in \eqref{eq:Momega}, $\bomega\neq \bm{0}$, can be diagonalized as
    \begin{align}\label{eq:diag_M_omega}
        \Mom \Vom
        =\Vom
        \begin{pmatrix}|\bomega|^2 &0\\ 0& -|\bomega|^2\end{pmatrix}, \quad 
        \Vom=\begin{pmatrix}
        \omega_2 &-\omega_3\\ \omega_3 & \omega_2
        \end{pmatrix}.
    \end{align}
If $\bomega=\bm{0}$, we set $V_0=I_2$ and the same decomposition holds.
    The following decompositions hold 
    \begin{align*}
    \begin{pmatrix}
        \Mom & (\kappa+s)^2 I_2\\(\kappa+s)^2 I_2 & \Mom
        \end{pmatrix} \begin{pmatrix}
        \Vom & \Vom\\-\Vom & \Vom
        \end{pmatrix}
        &=
        \begin{pmatrix}
        \Vom & \Vom\\-\Vom & \Vom
        \end{pmatrix}
        \diag(-\hat{\lambda}_{\bomega}^-, -\hat{\lambda}_{\bomega}^+, \hat{\lambda}_{\bomega}^+, \hat{\lambda}_{\bomega}^-),
        \\
        \begin{pmatrix}
        \Mom & (\bar\kappa+\bar s)^2 I_2\\(\bar\kappa+\bar s)^2 I_2 & \Mom
        \end{pmatrix} \begin{pmatrix}
        \Vom & \Vom\\-\Vom & \Vom
        \end{pmatrix}
        &=
        \begin{pmatrix}
        \Vom & \Vom\\-\Vom & \Vom
        \end{pmatrix}
        \diag(
        -\overline{\hat{\lambda}_{\bomega}^-}, 
        -\overline{\hat{\lambda}_{\bomega}^+}, \overline{\hat{\lambda}_{\bomega}^+}, 
        \overline{\hat{\lambda}_{\bomega}^-}),
    \end{align*}
    with $\hat{\lambda}_{\bomega}^+ = (\kappa+s)^2 + |\bomega|^2 $ and $\hat{\lambda}_{\bomega}^-=(\kappa+s)^2-|\bomega|^2$.
    Therefore, $\Phat$ and $\Phat^{\rH}$ share the same eigenvectors and the eigenvalues of $P^{\rH} P$ are given by
    \begin{align*}
        \frac{|(\kappa+s)^2+|\bomega|^2 |^2}{4|\kappa| |s| |\kappa+s|^2},\quad \frac{| (\kappa+s)^2-|\bomega|^2|^2}{4|\kappa| |s| |\kappa+s|^2},
    \end{align*}
    each with algebraic and geometric multiplicities of two.
    Since 
    \begin{align*}
    (\kappa+s)^2 = \kappa^2+2\kappa s+s^2 = 2s^2 + 2\kappa s + |\bomega|^2,
    \end{align*}
    we further obtain
    \begin{align*}
    \Re( (\kappa+s)^2 -|\bomega|^2 ) &= 2\, \Re( s(\kappa+s)),\\
    \Im( (\kappa+s)^2 \pm |\bomega|^2) &= 2\, \Im( s(\kappa+s)) = 2\, \Im(\kappa(\kappa+s)),\\
    \Re( (\kappa+s)^2+|\bomega|^2) &= \Re(2 s(\kappa+s)+2|\bomega|^2) = 2\, \Re(\kappa(\kappa+s)),
    \end{align*}
    where we used that $\Im(\kappa^2)=\Im(s^2)$ and $|\bomega|^2=(\kappa-s)(\kappa+s)$.
    Hence, we obtain
    \begin{align*}
        \hat{\lambda}_1 &= \frac{|(\kappa+s)^2-|\bomega|^2 |^2}{4|\kappa| |s| |\kappa+s|^2} = \frac{\Re( s(\kappa+s))^2 +  \Im( s(\kappa+s))^2}{|\kappa| |s| |\kappa+s|^2} = \frac{|s|}{|\kappa|},\\
        \hat{\lambda}_2 &= \frac{|(\kappa+s)^2+|\bomega|^2 |^2}{4|\kappa| |s| |\kappa+s|^2} = \frac{\Re( \kappa(\kappa+s))^2 +  \Im( \kappa(\kappa+s))^2}{|\kappa| |s| |\kappa+s|^2} = \frac{|\kappa|}{|s|}.
    \end{align*}
    Since the eigenvalues of \(\Phat^{\rH}\Phat\) coincide with those of \(\Phat\Phat^{\mathrm{H}}\), and the eigenvalues of \(\Phat^{-\mathrm{H}}\Phat^{-1}\) are the reciprocals of the eigenvalues of \(\Phat\Phat^{\mathrm{H}}\), the assertion concerning the eigenvalues follows immediately. 
    Moreover, as $|\Phat^{-1}\widehat{\bm{y}}|^2 = (\Phat^{-\mathrm{H}}\Phat^{-1}\widehat{\bm{y}},\,\widehat{\bm{y}})$, estimating \(|\kappa/s|\) and \(|s/\kappa|\) with the bounds \(\nu/2 \le |s'|,\,|\kappa'| \le 1\) valid on \(\cU_1\) by \eqref{eq:U1} yields the desired result.
\end{proof}

The next sequence of results verifies that $\cRh_{1}$ is a symmetrizer symbol  on $\cU_1$ for $\Mhat(\bomega,s)$ given by \eqref{eq:matrix-M}.
\begin{lemma}
    \label{lem:case1_symmetrizer_fulfills_property_i}
    For any $\eta\geq \eta_0>0$ the symmetrizer symbol \(\cRh_{1}\in  (S_{1,0}^0(\RR^3))^{4\times 4}\) defined in
    \eqref{eq:symmetrizer-P} satisfies property \ref{enum:symm-1} on $\cU_1$ and $|\cRh_{1}|\leq 2d/\nu$ for all $(\bomega,s)\in\cU_1$.
\end{lemma}
\begin{proof}
Clearly, $\cRh_{1}$ is Hermitian and smooth and $\Phat^{-1}$ is zero-homogeneous. Then, relation \eqref{eq:symmetrizer-P} implies that $\cRh_{1}$ is also zero-homogeneous. The boundedness of \(\cRh_{1}\) follows from Lemma~\ref{lem:eigsPHP} and $d\geq 1$.
\end{proof}
\begin{lemma}
    \label{lem:case1_symmetrizer_fulfills_property_ii}
    For any $\eta\geq \eta_0>0$ the symmetrizer symbol \(\cRh_{1}\) defined in
    \eqref{eq:symmetrizer-P} satisfies property \ref{enum:symm-2} on $\cU_1$ with $\delta_1=\nu/2$ provided $d\geq 5$.
\end{lemma}
\begin{proof}
Introduce \(\byh=\Phat\hat{\bm{v}}\).
The boundary condition $\byh_-=\Sop\byh_+ -\bfgh_-$ then becomes in these variables
    \begin{align*}
         (\Phat_{21}- \Sop \Phat_{11}) \bvh_1 = (\Sop \Phat_{12}-\Phat_{22}) \bvh_2 -\bfgh_-,
    \end{align*}
    where $\Phat_{ij}$, $i,j=1,2$ denote the corresponding blocks of $\Phat$.
Using $\Sop=\beta_J^{-1}\Mom$ we compute
    $$\Phat_{21}- \frac{1}{\beta_J}\Mom  \Phat_{11} =\frac{\beta_J(\kappa+s)^2-|\bomega|^4}{2\beta_J(\kappa+s)\sqrt{s\kappa}}I_2.$$
    Recalling the definition of $\beta_J$ given in \eqref{eq:rho-m-rho-p} 
    and that $|\bomega|^2=(\kappa+s)(\kappa-s)$, we compute
    \begin{align*}
        (\beta_J(\kappa+s)^2-|\bomega|^4)^{-1} = \frac{1-\ga}{4(\kappa+s)^2s\kappa},
    \end{align*}
    and obtain
    \begin{align*}
       \mnorm{{(\Phat_{21}- \frac{1}{\beta_J}\Mom  \Phat_{11})^{-1}}}
        = \left|\frac{(1-\gamma_J)\beta_J}{2(\kappa+s)\sqrt{s\kappa}}\right|
        = \left|\frac{(\kappa+s)(1-(\frac{\kappa-s}{\kappa+s})^2\gamma_J)}{2\sqrt{s\kappa}}\right|,
    \end{align*}
    The above term
     is bounded by some constant $C>0$ on $\cU_1$ using Lemma \ref{lem:bound_as-k}.
    We next compute
    \begin{align*}
        (\Phat_{21}- \frac{1}{\beta_J}\Mom  \Phat_{11})^{-1} (\frac{1}{\beta_J}\Mom  \Phat_{12}-\Phat_{22}) = \frac{\beta_J-(\kappa+s)^2}{|\bomega|^4-(\kappa+s)^2\beta_J } \Mom.
    \end{align*}
Similar as before, we compute
    \begin{align*}
        \frac{\beta_J-(\kappa+s)^2}{|\bomega|^4-(\kappa+s)^2\beta_J }
        &=   \frac{(1-\ga)\beta_J-(1-\ga)(\kappa+s)^2}{(1-\ga)(\kappa+s)^2(\kappa-s)^2-(\kappa+s)^2\big( (\kappa+s)^2-(\kappa-s)^2 \ga \big)}\\
        &=    \frac{(\kappa+s)^2-(\kappa-s)^2
\ga-(1-\ga)(\kappa+s)^2}{(\kappa+s)^2(\kappa-s)^2-(\kappa+s)^2(\kappa+s)^2}\\
        &=  -\frac{1}{(\kappa+s)^2}  \frac{\ga \Big((\kappa+s)^2-(\kappa-s)^2\Big)}{(\kappa+s)^2-(\kappa-s)^2}\\
        &=  -\frac{\ga}{(\kappa+s)^2}.
    \end{align*}
    Because $\mnorm{\Mom}= |\bomega|^2=\kappa^2-s^2$,
    we have that
    \begin{align*}
         \left|  \frac{\beta_J-(\kappa+s)^2}{|\bomega|^4-(\kappa+s)^2\beta_J }\right| \mnorm{\Mom}=
          \left|  \frac{\kappa-s}{\kappa+s}  \ga\right|\leq 1,
    \end{align*}
    by Lemma~\ref{lem:bound_as-k},
    therefore, with $C>0$ introduced above we have that
    \begin{align*}
        |\bvh_1| \leq |\bvh_2| + C |\bfgh_-|\quad\text{for all }(\bomega,s)\in\cU_1.
    \end{align*}
This allows us to verify property \ref{enum:symm-2} as follows,
\begin{align*}
    \byh^{\mathrm{H}}\cRh_{1}\byh = \hat{\bm{v}}^{\mathrm{H}} \begin{pmatrix}
        - I_2 & \\ & dI_2
    \end{pmatrix}\hat{\bm{v}}
    &= - |\hat{\bm{v}}_1|^2+ d |\hat{\bm{v}}_2|^2 \\
    &= |\hat{\bm{v}}|^2 +(d-1) |\hat{\bm{v}}_2|^2- 2 |\hat{\bm{v}}_1|^2 \\
    &\geq  |\hat{\bm{v}}|^2 +(d-5) |\hat{\bm{v}}_2|^2  -4C^2|\bfgh_-|^2\\
    &\geq  |\hat{\bm{v}}|^2  -4C^2|\bfgh_-|^2,
\end{align*}
provided $d\geq 5$.
The claim follows from Lemma \ref{lem:eigsPHP}.
\end{proof}

\begin{lemma}\label{lem:case1_symmetrizer_fulfills_property_iii}
For any $\eta\geq \eta_0>0$ the symmetrizer symbol \(\cRh_{1}\) defined in
    \eqref{eq:symmetrizer-P} satisfies property \ref{enum:symm-3} on $\cU_1$ with $\delta_2=\nu/2$ provided $d\geq 1$.
\end{lemma}
\begin{proof}
A direct calculation using \eqref{eq:pdp} and writing $\Mhat=\Mhat(\bomega,s)$ results in
\begin{align*}
    2\Re (\cRh_{1}\Mhat)&=\cRh_{1}\Mhat + \Mhat^{\mathrm{H}}\cRh_{1} 
    = (\Phat^{-1})^{\mathrm{H}} \left( Q_d \Dhat+\Dhat^{\mathrm{H}}Q_d  \right) \Phat^{-1} \\
    &\geq (\Phat^{-1})^{\mathrm{H}} \begin{pmatrix}
        2\Re (\kappa) I_2 & \\ & 2 d \Re (\kappa) I_2
    \end{pmatrix}  \Phat^{-1}\\
    &\geq 2\Re (\kappa) \min \left\{ 1,d \right\}(\Phat^{-1})^{\mathrm{H}} \Phat^{-1}\\
    &=  2 \Re (\kappa) \min\left\{\frac{|s|}{|\kappa|},\frac{|\kappa|}{|s|}\right\} I_4\\
    &\geq  \eta \nu I_4,
\end{align*}
where we have used Lemma~\ref{lem:eigsPHP} and \(d\geq 1\) in the third step, and Lemma~\ref{lem:Rek} in the last step.
The result follows.
\end{proof}

\subsection{Symmetrizer symbol for \texorpdfstring{$\Mhat$}{} on \texorpdfstring{$\cU_2$}{}}\label{subs:case2}
In view of Lemma~\ref{lem:eigsPHP},
the matrix symbol $\Phat^{-1}$ becomes singular if the scaled variable $\kappa'$ is close to $0$,
which occurs in a neighborhood of the points $(\bomega_0',s_0')$, where $s_0' = \pm \imath |\bomega_0'|$ and $|\bomega_0'| = 1/\sqrt{2}$, which follows from $\kappa'=0$ and $|\bomega'|^2+|s'|^2=1 $ and Lemma~\ref{lem:lower_bound_omega}.
Therefore, another approach must be taken to transform the matrix symbol $\Mhat(\bomega,s)$ into block form if $(\bomega,s)\in\cU_2$.
We note that the points \( (\bomega_0,s_0)\) are called ``special points'' in \cite[Section 2.1]{Kato1966}.
Different from the procedure described in \cite[Lemma~2.3]{Kreiss1970} (see also \cite[Chapter~2]{Kato1966}), which relies on an expansion of the spectral projection operator to derive a block decomposition of the matrix symbol $\Mhat(\bomega,s)$, we provide an explicit construction of the symmetrizer on $\cU_2$.
Using the matrix decompositions as in \cite[Section 2.1]{Feng1999}, compare also to \cite[Section 2, pp.626ff]{Majda1975}, we find that the matrix symbol \(\Mhat \in  (S_{1,0}^1(\RR^3))^{4\times 4}\) defined in \eqref{eq:matrix-M} can be decomposed as
\begin{align}\label{eq:matrix-M-1}
    \Mhat(\bomega,s) = \Qomega \NF \Qomega^{\mathrm{H}}
\end{align}
with
\[
    \NF = \begin{pmatrix}
        & & & s\\
        & & -\frac{\kappa^2}{s} &\\
        & -s& & \\
        \frac{\kappa^2}{s} & & & \\
    \end{pmatrix},\  
    \Qomega=\frac{1}{\sqrt{2}|\bomega|}\begin{pmatrix}
        - \omega_{3} & - \omega_{2} & \omega_{2} & - \omega_{3}\\ 
        \omega_{2} & - \omega_{3} & \omega_{3} & \omega_{2}\\
        - \omega_{3} & - \omega_{2} & - \omega_{2} & \omega_{3}\\
        \omega_{2} & - \omega_{3} & - \omega_{3} & -\omega_{2}
    \end{pmatrix},
\]
where \(\Qomega\) is unitary. Introducing the matrices
\begin{align*}
    \Bc = 
    \begin{pmatrix}
        1 & 0 & 1 & 0 \\
        0 & 1 & 0 & -1\\
        0 & 0 & 1 & 0\\
        0 & 0 & 0 & 1
    \end{pmatrix},\quad
    \Br = 
    \begin{pmatrix}
        1 & 0 & 0 & 0\\
        0 & 1 & 0 & 0\\
        0 & 0 & 1 & 0\\
        0 & 1 & 0 & 1
    \end{pmatrix},
\end{align*}
and defining 
\begin{equation}
    \Tom=\Qomega\Bc \Br, \label{eq:T-mat}
\end{equation}
we obtain the following block decomposition,
\begin{align}\label{eq:case2_block_decomp_M}
    \Tom^{-1} \Mhat(\bomega,s) \Tom = I_2 \otimes \Mhat_1,\qquad \Mhat_1=\begin{pmatrix} 0& s\\ \kappa^2/s & 0\end{pmatrix},
\end{align}
where $\otimes$ denotes the usual Kronecker product of matrices.
We note that the matrix symbol $\Mhat_1$ is homogeneous of degree one, hence $\Mhat_1 \in  (S_{1,0}^1(\RR^3))^{2\times 2}$, while $\Tom$ is homogeneous of degree zero.
The eigenvalues of $\Mhat_1$ are clearly $\pm\kappa$ and in the limiting case $\kappa\to 0$, $\Mhat_1$ becomes proportional to a Jordan block
$$
J=\begin{pmatrix}0&1\\0&0\end{pmatrix}.
$$
Hence, an eigendecomposition of $\Mhat_1$ will not be useful for the construction of a symmetrizer. Instead, we apply a Taylor series expansion in $\eta$ of the matrix symbol \(\Mhat_1(\bomega,\eta+\imath\xi)\) around \(\eta=0\) to arrive at an approximate decomposition.

\begin{lemma}\label{lem:case2_taylor_jordan}
The matrix symbol $\Mhat_1(\bomega,s)\in  (S_{1,0}^1(\RR^3))^{2\times 2}$ has the following expansion
\begin{align*}
\Mhat_1(\bomega,\eta+\imath \xi) = \imath \xi H-\imath\frac{|\bomega|^2}{\xi} J^\top  + \left( H +\frac{|\bomega|^2}{\xi^2} J^\top \right) \eta - \sum_{k=2}^\infty \imath^{k+1} \frac{|\bomega|^2}{\xi^{k+1}} \eta^k J^\top
\end{align*}
for $(\bomega,s)\in\cU_2$ and $H=J+J^\top$.
\end{lemma}
\begin{proof}
Because $\kappa^2=s^2+|\bomega|^2$, we observe that
\begin{align*}
\Mhat_1(\bomega,s) 
= (\eta+\imath\xi) \begin{pmatrix}
    0 & 1\\1 &0
\end{pmatrix}
+ \frac{|\bomega|^2}{\eta+\imath\xi}\begin{pmatrix}0 & 0\\1 &0\end{pmatrix}.
\end{align*}
The result then follows from the Taylor expansion around $\eta=0$ and 
\begin{align*}
\frac{\partial^k}{\partial\eta^k} \frac{1}{\eta+\imath\xi} 
= (-1)^k \frac{k!}{(\eta+\imath \xi)^{k+1}}
\end{align*}
and noting that $\eta^k/\xi^{k}=(\eta'/\xi')^k$ is smaller than one for $\nu$ sufficiently small by Lemma~\ref{lem:lower_bound_omega}. In fact, using Lemma~\ref{lem:lower_bound_omega}, there exists a universal constant $C>0$ such that
\begin{align*}
    \left|\sum_{k=2}^\infty \imath^{k+1} \frac{|\bomega|^2}{\xi^{k+1}} \eta^k \right|\leq \frac{|\bomega|^2}{\xi^2 |\xi-\eta|}\eta^2
    =\frac{|\bomega'|^2}{(\xi')^2 |\xi'-\eta'|}\eta'\, \eta
    \leq C\nu\eta
\end{align*}
for all $0<\nu<1/(2\sqrt{2})$, where we used a geometric series in the first step.
\end{proof}

The eigenvalues of $\Tom^{-\rH} \Tom^{-1}$ can be computed explicitly.
For later reference, we state the next result.
\begin{lemma}\label{lem:case2:lowerbound_Vinv}
For $\Tom $ defined in \eqref{eq:T-mat}, 
the eigenvalues of $\Tom^{-\rH} \Tom^{-1}$ are $(3\pm\sqrt{5})/2$, each with a multiplicity of two.
In particular,
\begin{align*}
\frac{3-\sqrt{5}}{2} |\byh|^2\leq |\Tom^{-1} \byh|^2 \leq \frac{3+\sqrt{5}}{2} |\byh|^2 \quad\text{for all }\byh\in\mathbb{C}^4.
\end{align*}
\end{lemma}
The expansion of $\Mhat_1(\bomega,s)$ derived in Lemma~\ref{lem:case2_taylor_jordan} is the basis for the construction of the symmetrizer on $\cU_2$, cf. \cite{Kreiss1970}.
We consider the following ansatz for the symmetrizer symbol $\hat{\cR}_2 \in  (S_{1,0}^0(\RR^3))^{4\times 4}$:
\begin{align}\label{eq:case2_def_symmetrizer_R}
\cRh_{2}=\Tom^{-\rH}(I_2\otimes \hat{R}_s) \Tom^{-1}, 
\quad \hat{R}_s=D-\imath \frac{\eta}{\xi} F,
\end{align}
with 
$
D=\begin{pmatrix}
    0&1\\1&d_{22}
\end{pmatrix}$,
 $F=\begin{pmatrix}0 &-f_{12}\\f_{12}&0\end{pmatrix}$ and $d_{22}, f_{12}\in \RR$.
The next sequence of results verifies that $\cRh_{2}$ is a symmetrizer symbol for $\Mhat(\bomega,s)$ on $\cU_2$.
We will make repeated use of the following:
By Lemma~\ref{lem:lower_bound_omega}, there exists a constant $C>0$ such that 
\begin{align}\label{eq:case2_eta_over_xi_is_1}
    \left|1-\left(\frac{|\bomega|}{\xi}\right)^2\right| + \left|\frac{\eta}{\xi}\right|\leq C\nu \quad\text{for all }(\bomega,s)\in\cU_2.
\end{align}
\begin{lemma}\label{lem:case2_symmetrizer_prop1}
    The symmetrizer symbol $\cRh_{2}$ defined in \eqref{eq:case2_def_symmetrizer_R} satisfies property {\ref{enum:symm-1}}. In particular, there exists a positive constant $C>0$ that depends on $d_{22}$ and $f_{12}$ such that
    $| \cRh_{2}| \leq C$ for all $(\bomega,s)\in\cU_2$.
\end{lemma}
\begin{proof}
Since $D$ is symmetric and $F$ is skew-symmetric, the matrix symbol $\cRh_s$ is Hermitian.
Moreover, $\cRh_{2}$ is zero-homogeneous and smooth for $(\bomega,s)\in\cU_2$.
Boundedness of $\cRh_{2}$ follows from Lemma~\ref{lem:case2:lowerbound_Vinv} and \eqref{eq:case2_eta_over_xi_is_1}.
\end{proof}

Next, we turn to the verification of property \ref{enum:symm-2}, which requires the following auxiliary result.

\begin{lemma}\label{lem:case2_bounding_on_boundary}
    For $\nu$ sufficiently small, there is a constant $C_{2}>0$ such that for any $\byh\in\mathbb{C}^4$ and $\bfgh_-\in\CC^2$ that satisfy $\byh_-=\Sop\byh_+ -\bfgh_-$ and $\bvh=\Tom^{-1}\byh=(v_{11},v_{12},v_{21},v_{22})^\top$ we have the following
    \begin{align*}
         \sum_{j=1}^2  |v_{j,1}|^2  \leq C_{2} \left(|\bfgh_-|^2+\sum_{j=1}^2 |v_{j,2}|^2\right) \quad\text{for all }(\bomega,s)\in\cU_2.
    \end{align*}
\end{lemma}

\begin{proof}
    Denote by $\Pp$ the permutation matrix
    \begin{align*}
        \begin{pmatrix}1\\&0&1\\&1&0\\&&&1\end{pmatrix}.
    \end{align*}
    Introduce $\bwh=\Pp\bvh$. Since $\Tom\Pp\bwh=\byh$, the boundary condition $\byh_-=\Sop\byh_+-\bfgh_-$ then becomes in these variables
    \begin{align*}
        (\That^p_{21}- \Sop \That^p_{11}) \bwh_1 = (\Sop \That_{12}^p-\That_{22}^p) \bwh_2-\bfgh_-,
    \end{align*}
    where $\That^p_{ij}$, $i,j=1,2$ denote the corresponding blocks of $\Tom\Pp$.
    Solving this linear system, we obtain
    \begin{align*}
        \bwh_1 = -\frac{|\bomega|^2-\beta_J}{|\bomega|^2+\beta_J}\bwh_2-(\That^p_{21}- \Sop \That^p_{11})^{-1}\bfgh_-.
    \end{align*}
    We compute
    \begin{align*}
        (\That^p_{21}- \Sop \That^p_{11})^{-1} = \frac{\beta_J\sqrt{2}}{|\bomega|(|\bomega|^2+\beta_J)}\begin{pmatrix}\omega_2-\omega_3 & \omega_2+\omega_3\\-\omega_2 & -\omega_3\end{pmatrix},
    \end{align*}
    from which we deduce that
    \begin{align*}
        \mnorm{(\That^p_{21}- \Sop \That^p_{11})^{-1}} =\left|\frac{\beta_J}{|\bomega|^2+\beta_J} \right| \frac{\sqrt{5}+3}{\sqrt{2}}.
    \end{align*}
    Recalling the definition of $\beta_J$ given in \eqref{eq:rho-m-rho-p}, we observe that
    \begin{align}\label{eq:case2_aux2}
        (1-\gamma_J)\beta_J = (\kappa+s)^2 \left( 1- \frac{(\kappa-s)^2}{(\kappa+s)^2}\gamma_J\right). 
    \end{align}
    Similarly, some algebra yields that
    \begin{equation}\label{eq:case2_aux0}
         \begin{split}
        (1-\ga)(|\bomega|^2-\beta_J)&=(1-\ga)(\kappa-s)(\kappa+s)-(\kappa+s)^2+(\kappa-s)^2 \ga\\
        &=-2s(\kappa+s)\left(1+\frac{\kappa-s}{\kappa+s}\ga\right),
        \end{split}
    \end{equation}
    and
    \begin{equation}\label{eq:case2_aux1}
        \begin{split}
        (1-\ga)(|\bomega|^2+\beta_J)&=(1-\ga)(\kappa-s)(\kappa+s)+(\kappa+s)^2 -(\kappa-s)^2 \ga\\
        &=2\kappa(\kappa+s)\left(1-\frac{\kappa-s}{\kappa+s}\ga\right).
        \end{split} 
    \end{equation}
    Using the asymptotic expansions of $\ga$ for $|\kappa'|\ll 1$ given in Lemma~\ref{lem:kappa_asymptotics}, we obtain
            \begin{align*}
                \frac{|\bomega|^2-\beta_J}{|\bomega|^2+\beta_J} 
                = -\frac{c_1(\kappa/s)}{2-c_1(\kappa/s)\kappa/s} + O(|\kappa/s|^2),
        \end{align*}
    and
    \begin{align*}
        \frac{\beta_J}{|\bomega|^2+\beta_J} = \frac12(1+\kappa/s) \frac{c_2(\kappa/s)}{2-c_1(\kappa/s)\kappa/s} + O(|\kappa/s|^2).
    \end{align*}
    Since $|\bwh_2|^2=|v_{12}|^2+|v_{22}|^2$ and $|\bwh_1|^2=|v_{11}|^2+|v_{21}|^2$ the assertion is proven.
\end{proof}

\begin{lemma}\label{lem:case2_property_ii_holds}
Let $f_{12}= 1/(C\nu)$, where $C>0$ is the constant from \eqref{eq:case2_eta_over_xi_is_1},
and $d_{22}=3+3C_2$, where $C_2$ is the constant from Lemma~\ref{lem:case2_bounding_on_boundary}.
Then the symmetrizer symbol $\cRh_{2}$ defined in \eqref{eq:case2_def_symmetrizer_R} satisfies property \ref{enum:symm-2} on $\cU_2$ with $\delta_1=(3-\sqrt{5})/2$.
\end{lemma}
\begin{proof}
Writing $\bvh=\Tom^{-1} \byh$, using the definition of $\cRh_{2}$ in \eqref{eq:case2_def_symmetrizer_R}, the bounds $|\Im(z_1z_2)|\leq |z_1||z_2| $ and $|\Re(z_1z_2)|\leq |z_1||z_2|$ for any $z_1,z_2\in\mathbb{C}$ as well as the Cauchy-Schwarz inequality, we have that
\begin{align*}
    \byh^{\rH}\cRh_{2} \byh &=  \sum_{j=1}^2 \bvh_j^{\rH} \hat{R}_s \bvh_j\\
&= \sum_{j=1}^2 2\,\Re(v_{j,1}\overline{v_{j,2}})+d_{22}|v_{j,2}|^2  +2\frac{\eta}{\xi}f_{12}\Im(v_{j,1}\overline{v_{j,2}})\\
        &\geq \sum_{j=1}^2 -(1+\frac{\eta}{\xi}f_{12})|v_{j,1}|^2 + (d_{22}-1-\frac{\eta}{\xi} f_{12})|v_{j,2}|^2 \\
        &\geq \sum_{j=1}^2 -2|v_{j,1}|^2 + (d_{22}-2)|v_{j,2}|^2=(*),
\end{align*}
where we used the assumption that $f_{12}\leq 1/(C\nu)$ in the last step.
We continue
\begin{align*}
        (*)&=\sum_{j=1}^2  |\hat{\bm{v}}_j|^2 -3|v_{j,1}|^2 + (d_{22}-3)|v_{j,2}|^2 \\
        &\geq \sum_{j=1}^2 |\hat{\bm{v}}_j|^2 -3C_{2}|v_{j,2}|^2 + (d_{22}-3)|v_{j,2}|^2 -3C_2|\bfgh_-|^2,
\end{align*}
 where we used Lemma~\ref{lem:case2_bounding_on_boundary} in the last step.
    By choosing $d_{22}$ such that
     $d_{22}\geq 3+3C_{2}$, we  have thus shown that
    \begin{align*}
     \byh^{\rH}\cRh_{2} \byh\geq |\bvh|^2 -C|\bfgh_-|^2\quad\text{for all } (\bomega,s)\in\cU_2,
    \end{align*}
    for some $C>0$.
    The claim then follows from Lemma~\ref{lem:case2:lowerbound_Vinv}.
\end{proof}
\begin{lemma}\label{lem:case2_symmetrizer_prop2}
For sufficiently small $\nu>0$, there exists $\delta_2>0$ such that the symmetrizer symbol $\cRh_{2}$ defined in
    \eqref{eq:case2_def_symmetrizer_R} satisfies property {\ref{enum:symm-3}} in $\cU_2$.
\end{lemma}
\begin{proof}
Due to the block structure of $\Tom^{-1}\Mhat(\bomega,s)\Tom$, the definition of $\cRh_{2}$  \eqref{eq:case2_def_symmetrizer_R}, and Lemma~\ref{lem:case2_taylor_jordan}, it suffices to investigate $\hat{R}_s \Mhat_1$.
We have that
\begin{align*}
    \hat{R}_s\Mhat_1&= \imath \left(\xi DH - \frac{|\bomega|^2}{\xi}DJ^\top\right)+ \left( (D+F)H +  \frac{|\bomega|^2}{\xi^2} (D-F)J^\top\right)\eta + O(\nu\eta).
\end{align*}
Noting that $DH$, $DJ^\top$ and $\xi$, $|\bomega|$ are real, we have
\begin{align*}
\Re(\hat{R}_s\Mhat_1) = \left( (D+F)H +  \frac{|\bomega|^2}{\xi^2} (D-F)J^\top\right)\eta + O(\nu\eta),
\end{align*}
and therefore, with $\gamma=(|\bomega|/\xi)^2$,
\begin{align*}
\frac{1}{2}\Re(\hat{R}_s\Mhat_1+\Mhat_1^{\rH}\hat{R}_s)
=\begin{pmatrix} \gamma(f_{12}+1) +1-f_{12} & \frac{1}{2}(\gamma+1)d_{22}\\
\frac{1}{2}(\gamma+1)d_{22} & 1+f_{12}
 \end{pmatrix} 
+ O(\nu\eta).
\end{align*}
The eigenvalues of the last matrix are given by
$$
1+ \frac{\gamma}{2}(1+f_{12}) \pm \frac{1}{2} \sqrt{ (\gamma-2)^2 f_{12}^2+2\gamma(\gamma-2)f_{12} +(\gamma+1)^2d_{22}^2+\gamma^2}.
$$
Note that the term within the square-root is positive for sufficiently small $\nu$, because then $f_{12}=1/(C\nu)$ as chosen in Lemma~\ref{lem:case2_property_ii_holds} is large.
Choosing $\nu$ such that
$$
\frac{\gamma}{2}(1+f_{12}) \geq  \frac{1}{2} \sqrt{ (\gamma-2)^2 f_{12}^2+2\gamma(\gamma-2)f_{12} +(\gamma+1)^2d_{22}^2+\gamma^2}$$
shows that the eigenvalues are bounded from below by $1$.
In fact, the latter bound is equivalent to
\begin{align}\label{eq:case2_tmp}
2(2\gamma-1) f_{12}^2 + 4\gamma f_{12}\geq (\gamma+1)^2 d_{22}^2.
\end{align}
In view of \eqref{eq:case2_eta_over_xi_is_1} and upon decreasing $\nu$ if necessary, we have that $|1-\gamma|\leq 1/2$ for all $(\bomega,s)\in\cU_2$.
Therefore, we find $f_{12}=1/(C\nu)$ such that \eqref{eq:case2_tmp} holds for all $(\bomega,s)\in\cU_2$.
With this choice, we have that 
\begin{align*}
    \frac{1}{2}\Re(\hat{R}_s\Mhat_1+\Mhat_1^{\rH} \hat{R}_s)\geq  (1-C\nu)\eta I_2,
\end{align*}
where $C>0$ is some constant.
Upon decreasing $\nu$ even further if necessary, we conclude that property \ref{enum:symm-3} holds.
\end{proof}

\subsection{Symmetrizer symbol for \texorpdfstring{$\Mhat$}{} on \texorpdfstring{$\cU_3$}{}} 
\label{subs:case3-new}
We recall the decomposition
\begin{align*}
    \Mhat(\bomega,s) = \Qomega \NF \Qomega^{\mathrm{H}},
\end{align*}
where \(\Qomega\) is unitary, see \eqref{eq:matrix-M-1}. 
On the subdomain $\cU_3$, the matrix symbol $\Mhat(\bomega,s)$ has a singularity for $s'\to 0$. The block decomposition \eqref{eq:case2_block_decomp_M} is not suitable in this case because of an unfortunate coupling in the equation for $\brho_s$, compare to \eqref{eq:case3_expression_rho_o} below. 
Otherwise, the layout of this section is similar to the previous ones.
Introducing the unitary matrices
\begin{align}
\QF= \frac{1}{\sqrt{2}}\begin{pmatrix}
    I_2 & I_2\\\imath I_2 & -\imath I_2
\end{pmatrix}
\label{eq:case3-QF}
\end{align}
and
\begin{equation}
    \label{eq:T-mat-2}
    \Tom=\Qomega \QF \in  (S_{1,0}^0(\RR^3))^{4\times 4}
\end{equation}
we will use the following decomposition in this section,
\begin{align}\label{eq:case3_block_decomp_M}
    \Mhat(\bomega,s) = \Qomega \NF \Qomega^{\mathrm{H}} = \Tom \begin{pmatrix} \Mhat_1&\\&\Mhat_2\end{pmatrix}\Tom^{\mathrm{H}},
\end{align}
with
\begin{align*}
    \Mhat_1 = \imath \begin{pmatrix}0&s\\-\kappa^2/s& 0\end{pmatrix},\quad \Mhat_2 = - \Mhat_1.
\end{align*}
The matrices $\Mhat_1$ and $\Mhat_2$ can be diagonalized as follows,
\begin{align}
    \Mhat_1 \Vhat_1 = \Vhat_1 \diag(-\kappa,\kappa),\quad 
    \Mhat_2 \Vhat_2 = \Vhat_2 \diag(-\kappa,\kappa),\label{eq:M1V1}
\end{align}
with
\begin{align*}
\Vhat_1 &= \begin{pmatrix} -\imath s/\kappa & \imath s/\kappa\\ 1 & 1\end{pmatrix},\quad
\Vhat_2 = \begin{pmatrix} \imath s/\kappa & -\imath s/\kappa\\ 1 & 1\end{pmatrix},\\
\Vhat_1^{-1} &= \frac{1}{2}\begin{pmatrix} \imath\kappa/s & 1\\ -\imath\kappa/s & 1\end{pmatrix},\quad
\Vhat_2^{-1} = \frac{1}{2}\begin{pmatrix} -\imath\kappa/s & 1\\ \imath \kappa/s & 1\end{pmatrix}.
\end{align*}
We note that precisely the first column of $\Vhat_1^{-1}$ and $\Vhat_2^{-1}$ is singular. This will become important in the following.
For a real number $d\geq 1$ to be specified below, we consider the following symmetrizer symbol $\hat{\cR}_3$ for $\Mhat$ on $\cU_3$
\begin{align}\label{eq:case3_def_symmetrizer_R}
    \cRh_{3} = \Tom^{-\rH} \diag(\cRh^1,\cRh^2) \Tom^{-1}, 
\end{align}
with singular blocks $\cRh^j = \Vhat_j^{-\rH}\diag(-1,d) \Vhat_j^{-1}$, $j=1,2$, which are equal to 
\begin{subequations}
\label{eq:V1H}
    \begin{align}
    \Vhat_1^{-\rH}\diag(-1,d)\Vhat_1^{-1} = 
    \begin{pmatrix} \frac{|\kappa|^2(d-1)}{4|s|^2} & \imath\frac{\overline{\kappa}(d+1)}{4\overline{s}} \\ -\imath\frac{\kappa(d+1)}{4{s}} & \frac{d-1}{4}\end{pmatrix}, \label{eq:V1H-1} \\
    \Vhat_2^{-\rH}\diag(-1,d)\Vhat_2^{-1} = 
    \begin{pmatrix} \frac{|\kappa|^2(d-1)}{4|s|^2} & -\imath\frac{\overline{\kappa}(d+1)}{4\overline{s}} \\ \imath\frac{\kappa(d+1)}{4{s}} & \frac{d-1}{4}\end{pmatrix}. \label{eq:V1H-2}
\end{align}
\end{subequations}
Using \eqref{eq:V1H} we can decompose $\cRh_3$ as 
\begin{equation}
    \byh^{\rH} \cRh_3 \byh = \bwh_1^{\rH}\cRh_{3,1} \bwh_1+ 
    \bwh_2^{\rH}\cRh_{3,2} \bwh_2 + 2\Re(\bwh_2^{\rH}\cRh_{3,3} \bwh_1)
    \label{eq:yR3y}
\end{equation}
with
\begin{subequations}
    \label{eq:R3-submat}
    \begin{alignat}{3}
    \cRh_{3,1}&\in (S_{1,0}^{-2,2}(\RR^3))^{2\times 2},\quad
    &&\cRh_{3,1}=\frac{1}{4}(d-1)\left|\frac{\kappa}{s}\right|^2 I_2, \label{eq:R3-1}\\
    \cRh_{3,2}&\in (S_{1,0}^{0}(\RR^3))^{2\times 2},\quad
    &&\cRh_{3,2}=\frac{1}{4}(d-1) I_2, \label{eq:R3-2}\\
    \cRh_{3,3}&\in (S_{1,0}^{-1,1}(\RR^3))^{2\times 2},\quad
    &&\cRh_{3,3}=\frac{1}{4}\imath (d+1)\frac{\kappa}{s}\begin{pmatrix}
        -1 & 0 \\ 0 & 1
    \end{pmatrix} , \label{eq:R3-3}
\end{alignat}
\end{subequations}
and $\bwh_1=(\hat{y}_1,\hat{y}_3)^\top$, $\bwh_2=(\hat{y}_2,\hat{y}_4)^\top$.

The following sequence of results verify counterparts of properties \ref{enum:symm-1}--\ref{enum:symm-3} for the singular case $|s'|<\nu$ with $\nu>0$ sufficiently small; see also \cite{Majda1975}. We begin with a generalization of property \ref{enum:symm-1}.
\begin{lemma}\label{lem:case3_property_i_holds}
The symmetrizer symbol $\cRh_{3}$ is zero homogeneous and Hermitian. Moreover,
\begin{align*}
      \Re(\cRh_{3}\byh,\bwh) \leq  \frac{d}{2}|\diag(|\kappa/s|,1,|\kappa/s|,1) \Tom^{-1}\byh| |\bwh|,\quad \text{for all } \bwh,\, \byh\in\CC^4.
\end{align*}
\end{lemma}
\begin{proof}
Clearly, $\cRh_{3}$ is homogeneous of degree zero and Hermitian pointwise. An application of the Cauchy-Schwarz inequality and using $\Tom$ is unitary shows that
\begin{align*}
\Re(\cRh_{3}\byh,\bwh)
    =\Re(\diag(\cRh^1,\cRh^2)\Tom^{-1} \byh ,\Tom^{-1} \bwh)\leq \frac{d}{2}|\diag(\cRh^1,\cRh^2) \Tom^{-1}\byh| |\bwh|.
\end{align*}
Since
    \begin{align*}
        \cRh^j = \Vhat_j^{-\rH} \diag(-1,d) \Vhat_j^{-1} \leq d \Vhat_j^{-\rH}\Vhat_j^{-1} = \frac{d}{2}\diag(|\kappa|^2/|s|^2,1),
    \end{align*}
    the result follows.
\end{proof}
For the verification of a generalization of the property \ref{enum:symm-2} to the singular case $s'\approx 0$, we prove the following auxiliary result.
\begin{lemma}\label{lem:case3_bounding_on_boundary}
    There exists a constant $C_3>0$ such that for all $\byh\in\mathbb{C}^4$ and $\bfgh_-\in\CC^2$ satisfying the boundary relation $\byh_-=\Sop\byh_+-\bfgh_-$, the transformed variable $\bvh=\Tom^{-1}\byh=(v_{11},v_{12},v_{21},v_{22})^\top$ satisfies the estimate
    \begin{align}
         \sum_{j=1}^2  |\frac{\kappa}{s}v_{j,1}|^2  \leq C_3\left(|\frac{\kappa}{s}\bfgh_-|^2+\sum_{j=1}^2 |\frac{s}{\kappa}v_{j,2}|^2\right)  \quad\text{for all }(\bomega,s)\in\cU_3. \label{eq:vj1-vj2}
    \end{align}
\end{lemma}
\begin{proof}
The proof is similar to the proof of Lemma~\ref{lem:case2_bounding_on_boundary}; recall the permutation matrix $\Pp$.
Introduce $\bwh=\Pp\bvh$. Since $\Tom\Pp\bwh=\byh$, the boundary condition $\byh_-=\Sop\byh_+-\bfgh_-$ then becomes in these variables
    \begin{align*}
        (\That^p_{21}- \Sop \That^p_{11}) \bwh_1 = (\Sop \That_{12}^p-\That_{22}^p) \bwh_2 -\bfgh_-,
    \end{align*}
    where $\That^p_{ij}$, $i,j=1,2$ denote the corresponding blocks of $\Tom\Pp$.
    We have that
\begin{align*}
\Sop \That_{12}^p-\That_{22}^p&= \frac{1}{2|\bomega|}\frac{|\bomega|^2-\beta_J}{\beta_J}
    \begin{pmatrix}
    -\omega_2+\imath\omega_3 & - \omega_2 -\imath\omega_3\\
    -(\imath\omega_2+ \omega_3) 
    & \imath\omega_2 -\omega_3
    \end{pmatrix},\\
    \That^p_{21}- \Sop \That^p_{11} &=\frac{1}{2|\bomega|} \frac{|\bomega|^2+\beta_J}{\beta_J}
    \begin{pmatrix}
    -(\imath\omega_2+\omega_3) & \imath \omega_2 -\omega_3\\
    \omega_2-\imath \omega_3 
    & \omega_2 + \imath\omega_3
    \end{pmatrix}.\notag
\end{align*}
Therefore,
\begin{align*}
    \bwh_1 = \imath \frac{|\bomega|^2-\beta_J}{|\bomega|^2+\beta_J} \diag(-1,1) \bwh_2
     - \frac{1}{|\bomega|}\frac{\beta_J}{|\bomega|^2+\beta_J}
     \begin{pmatrix}-\omega_3+\imath\omega_2 & \omega_2+\imath\omega_3\\-(\omega_3+\imath\omega_2) & \omega_2-\imath\omega_3 \end{pmatrix}
     \bfgh_-    .
\end{align*}
We observe that
\begin{align*}
    \left|\begin{pmatrix}-(\omega_2+\imath\omega_3) & \imath\omega_2-\omega_3\\-(\omega_2-\imath\omega_3) & -(\imath\omega_2+\omega_3) \end{pmatrix} \right|=\sqrt{2}|\bomega|.
\end{align*}
Moreover, using \eqref{eq:case2_aux2} and \eqref{eq:case2_aux1}, we compute
\begin{align*}
    \frac{\beta_J}{|\bomega|^2+\beta_J} 
    = -\frac{\kappa+s}{2\kappa} \frac{1- \left(\frac{\kappa-s}{\kappa+s}\right)^2\ga}{1- \left(\frac{\kappa-s}{\kappa+s}\right)\ga}= \frac{\kappa+s}{2\kappa} \frac{2-d_2(s/\kappa)s/\kappa}{2-d_1(s/\kappa)s/\kappa}+O(|s/\kappa|)\leq 1,
\end{align*}
where we used Lemma~\ref{lem:s_asymptotics} in the last step.
With \eqref{eq:case2_aux0} and \eqref{eq:case2_aux1} we obtain
\begin{align*}
    \frac{|\bomega|^2-\beta_J}{|\bomega|^2+\beta_J} &= -\frac{s}{\kappa} \frac{1+\frac{\kappa-s}{\kappa+s}\ga}{1-\frac{\kappa-s}{\kappa+s}\ga}
    =- \left(\frac{s}{\kappa} \right)^2 \frac{d_1(s/\kappa)}{2-d_1(s/\kappa)s/\kappa} +O (|s/\kappa|),
\end{align*}
where we used Lemma~\ref{lem:s_asymptotics} in the last step. Since $|\bwh_2|^2=|v_{12}|^2+|v_{22}|^2$ and $|\bwh_1|^2=|v_{11}|^2+|v_{21}|^2$ the assertion is proven.
\end{proof}

\begin{lemma}\label{lem:case3_property_ii_holds}
Given the symmetrizer symbol $\hat{\cR}_3$ stated in \eqref{eq:case3_def_symmetrizer_R} and let $d>1$ and $\nu$ be sufficiently small.
Then there are constants $\delta_1=\frac{1}{8}$ and $C>0$ such that
for any $\byh\in\mathbb{C}^4$ and $\bfgh_-\in\CC^2$ that satisfy $\byh_-=\Sop\byh_+ - \bfgh_-$, it holds that
\begin{align}\label{eq:case3_property_ii_holds}
    \byh^\rH \cRh_{3} \byh \geq \delta_1(\Tom^{-1}\byh)^\rH \diag( \frac{|\bzeta|^2}{|s|^2},1,\frac{|\bzeta|^2}{|s|^2},1) \Tom^{-1}\byh - C |\frac{\kappa}{s}\bfgh_-|^2
\end{align}
for all $(\bomega,s)\in\cU_3.$
\end{lemma}

\begin{proof}
Denote by $\bvh=\Tom^{-1} \byh$. Using the bounds $|\Im(z_1z_2)|$, $|\Re(z_1z_2)|\leq |z_1||z_2|$ for any $z_1,z_2\in\mathbb{C}$, the Cauchy-Schwarz inequality with $\varepsilon>0$ and  \eqref{eq:V1H}, we obtain 
\begin{equation*}
    \begin{split}
        \byh^\rH\cRh_{3} \byh &=  \sum_{j=1}^2 \bvh_j^\rH \Vhat_j^{-\rH}\diag(-1,d)\Vhat_j^{-1} \bvh_j\\
    &\geq \frac{(d-1)}{4}\sum_{j=1}^2  \left(\left|\frac{\kappa}{s} v_{j1}\right|^2 -2 \left|\frac{\kappa}{s} v_{j1}\right|\, |v_{j2}|\frac{d+1}{d-1} +  |v_{j2}|^2\right)\\
    &\geq \frac{(d-1)}{4}\sum_{j=1}^2  \left(\left(1-\frac{1}{\varepsilon}\frac{d+1}{d-1}\right)\left|\frac{\kappa}{s} v_{j1}\right|^2  +  \left(1-{\varepsilon}\frac{d+1}{d-1}\right)|v_{j2}|^2\right)\\
    &\geq \frac{(d-1)}{4}\frac{1}{2}\sum_{j=1}^2  \left(|\frac{\kappa}{s} v_{j1}|^2  +  |v_{j2}|^2\right)\\
    &\quad +\frac{(d-1)}{4}\sum_{j=1}^2  \left(\left(\frac{1}{2}-\frac{1}{\varepsilon}\frac{d+1}{d-1}\right)\left|\frac{\kappa}{s} v_{j1}\right|^2  +  \left(\frac{1}{2}-{\varepsilon}\frac{d+1}{d-1}\right)|v_{j2}|^2\right).
    \end{split}
\end{equation*}
Taking $\varepsilon=\frac{1}{4}\frac{d-1}{d+1}$ the above inequality becomes
\begin{equation}
    \label{eq:case3_prop_iii_aux}
    \begin{split}
        \byh^\rH\cRh_{3} \byh 
    &\geq \frac{(d-1)}{8} \sum_{j=1}^2  \left(\left|\frac{\kappa}{s} v_{j1}\right|^2  +  |v_{j2}|^2\right)\\
    &\quad +\frac{(d-1)}{4}\sum_{j=1}^2  \left(\left(\frac{1}{2}-4\left(\frac{d+1}{d-1}\right)^2\right)\left|\frac{\kappa}{s} v_{j1}\right|^2  +  \frac{1}{4}|v_{j2}|^2\right).
    \end{split}
\end{equation}
Since $\frac{1}{2}-4(\frac{d+1}{d-1})^2$ is negative,
using Lemma~\ref{lem:case3_bounding_on_boundary} we can bound $|v_{j1}|$ in terms of $|v_{j2}|$ and  we obtain from \eqref{eq:vj1-vj2}
\begin{align*}
    &\sum_{j=1}^2 \left(\frac{1}{2}-4\left(\frac{d+1}{d-1}\right)^2\right) \left|\frac{\kappa}{s} v_{j1}\right|^2  +  \frac{1}{4}|v_{j2}|^2 \\
    &\ge  C_3\left(\frac{1}{2}-4 \left(\frac{d+1}{d-1}\right)^2\right)|\frac{\kappa}{s}\bfgh_-|^2 \\
    &\quad + \sum_{j=1}^2 \left(C_3\left(\frac{1}{2}-4\left(\frac{d+1}{d-1}\right)^2 \right)\left|\frac{s}{\kappa}\right|^2+\frac{1}{4} \right) \left|v_{j2}\right|^2 \\
    &\ge  C_3\left(\frac{1}{2}-4 \left(\frac{d+1}{d-1}\right)^2\right)|\frac{\kappa}{s}\bfgh_-|^2,
\end{align*}
by decreasing $\nu$ if necessary and observing that 
$|s/\kappa|=|s'|/|\kappa'|\leq 2\nu$ by \eqref{eq:claim1-B4}.
Finally, using \eqref{eq:case3_prop_iii_aux}, we have
\begin{equation*}
    \begin{split}
        \byh^\rH\cRh_{3} \byh 
    \geq \frac{(d-1)}{8}\sum_{j=1}^2  \left(\left|\frac{\kappa}{s} v_{j1}\right|^2  +  |v_{j2}|^2\right)
     - C_3\left|\frac{1}{2}-4 \left(\frac{d+1}{d-1}\right)^2\right| |\frac{\kappa}{s}\bfgh_-|^2.
    \end{split}
\end{equation*}
Using \eqref{eq:claim1-B4} again, we have that $|\kappa/s|=|\kappa'|/|s'|\geq |\bzeta|/(2|s|)$, which shows the claim.
\end{proof}

The corresponding generalization of property \ref{enum:symm-3} is next verified.
\begin{lemma}\label{lem:case3_property_iii_holds}
For all $d\geq 1$, $(\bomega,s)\in\cU_3$ and $\byh\in\CC^4$, the symmetrizer symbol $\hat{\cR}_3$ stated in \eqref{eq:case3_def_symmetrizer_R} satisfies property (iii) with $\delta_2=\frac{1}{8}$
    \begin{align}\label{eq:case3_property_iii_holds}
   \byh^{\rH}  \Re(\cRh_{3} \Mhat) \byh \geq \delta_2\eta  (\Tom^{-1}\byh)^{\rH}\diag(\left|\frac{\bzeta}{s}\right|^2,1,\left|\frac{\bzeta}{s}\right|^2,1) \Tom^{-1}\byh.
\end{align}
\end{lemma}
\begin{proof}
For $j=1,2$, we observe that
    $\Vhat_j^{-\rH}\Vhat_j^{-1} = \frac{1}{2}\diag( |\kappa|^2/|s|^2,1)$.
Using $d\geq 1$ and \eqref{eq:M1V1}, we consequently obtain
\begin{align*}
    \Re(\cRh^j \Mhat_j) &= \Re{\left( \Vhat_j^{-\rH}\diag(\kappa,d\kappa)\Vhat_j^{-1}\right)}\\
    &\geq \Re(\kappa)\min\{1,d\} \Vhat_j^{-\rH}\Vhat_j^{-1}\\
    &\geq \frac{1}{2}\eta \diag( \left|\frac{\kappa}{s}\right|^2,1)
    \geq \frac{1}{8}\eta \diag( \left|\frac{\bzeta}{s}\right|^2,1)
\end{align*}
and the claim follows from the definition of $\cRh_{3}$.
\end{proof}

For later reference, we consider \eqref{eq:rho-o}. Since the bounds \eqref{eq:case3_property_ii_holds} and \eqref{eq:case3_property_iii_holds} involve the variable $\bvh=\Tom^{-1}\byh$, we apply the same transformation to \eqref{eq:rho-o}.
We then obtain that
\begin{align}\label{eq:case3_expression_rho_o}
    s\hat{\brho}_{\s}&=
    \imath \left( (\omega_2 \EE- \omega_3 \JJ)\That_{11}  + (  \omega_2 \DD- \omega_3 \II) \That_{12} )\right)\bvh_1 \notag\\
    &\quad +
    \imath \left( (\omega_2 \EE- \omega_3 \JJ)\That_{21}  + (  \omega_2 \DD- \omega_3 \II) \That_{22} )\right)\bvh_2\\
    &=\frac{|\bomega|}{\sqrt{2}} \begin{pmatrix}\imath\\-1\end{pmatrix}v_{11} + \frac{|\bomega|}{\sqrt{2}} \begin{pmatrix}\imath\\1\end{pmatrix}v_{21}. \notag
\end{align}

\subsection{Global symmetrizer}
The symmetrizer symbols $\hat{\cR}_j$ on the subdomains $\cU_j,\,j=1,2,3$, can be combined to define a global symmetrizer $\cR$ using a partition of unity subordinate to $\{\cU_j:\,1\leq j\leq 3\}$ with $\cU_j$ given by \eqref{eq:cases}.

\begin{definition}[Partition of Unity of $\cU$]
    \label{def:partition-of-psi-j}
    Let $\psih_j \in  S_{1,0}^0(\RR^3),\, j=1,2,3$ be such that $\supp(\psih_j)\subset \cU_j,$ and $0\leq \psih_j\leq 1$ with $\sum_{j=1}^3\psih_j=1$ on $\cU$.
    The partition of unity pseudo-differential operators $\psi_j \in \OpsS_{1,0}^0(\overline{V}\times\RR^3)$, $j=1,2,3$, are defined as
    \begin{equation*}
        \begin{split}
         &   \psi_j(\bx,t,\imath D_{\bxb}, \imath D_t+\eta)
        (e^{-\eta t}\bv(\bx,t) ) \\
        &\quad =(\cL\cF)_3^{-1}[\psih_j(\bomega,s)\bvh(x_1,\bomega,s)](\bx,t)
        \quad \forall\, \bv \in (C_{(0)}^\infty(\overline{V}))^\ell,
        \end{split}
    \end{equation*}
    with $\ell\in \NN$ and $\hat{\bv} (x_1,\bomega,s)
    =\int_{\RR^3} e^{-\imath \xi t-\imath\bomega\cdot\bxb-\eta t}\bv(\bx,t)\,\mathrm{d}\bxb \,\mathrm{d}t $
    and $\bv(\bx,t)=\bm{0}$ for $(\bx,t)\in \RR^4\setminus \overline{V}$.
\end{definition}

The pseudo-differential operators $\psi_j$ can be extended to bounded linear operators  $\psi_j:\,(C_b^\infty(\overline{V}))^\ell \to (C_b^\infty(\overline{V}))^\ell$ for all $\bv \in (C_b^\infty(\overline{V}))^n$ with $(C_b^\infty(\overline{V}))^n$ the Fr\'echet space
 $C_b^\infty(\overline{V})=\bigcap_{k=0}^{\infty}C_b^k(\overline{V})$, and $\|\psi_j\|_{C_b^k(\overline{V})}=\sup_{|\balpha|\leq k}\sup_{\bx\in \overline{V}}|\partial_{\bx}^{\balpha}\psi_j(\bx)|$, see e.g. \cite{Abels2012}, Theorem 3.28, or \cite{Kumano1982}, Chapter 2, and using $C_0^{\infty}$ cut-off functions for the restriction to $\overline{V}$, see \cite{Grubb1996, Wloka1987}.

\begin{definition}
    [Global symmetrizer]
    The global symmetrizer $\cR$ is defined as
    \begin{equation}
        \label{eq:total-symmetrizer}
        \begin{split}
     &   \cR(\bx,t,\imath D_{\bxb}, \imath D_t+\eta)
    (e^{-\eta t}\bv(\bx,t)) \\
    &=\sum_{j=1}^3 \cR_j(\bx,t,\imath D_{\bxb}, \imath D_t+\eta) \psi_j(\bx,t,\imath D_{\bxb}, \imath D_t+\eta)
    (e^{-\eta t}\bv (\bx,t)),\quad \forall \bv \in (C_{(0)}^\infty(\overline{V}))^\ell,
    \end{split}
    \end{equation}
    with $\cR_j$ the symmetrizers related to the symbols $\cRh_j$, $j=1,2,3$, stated in, respectively, \eqref{eq:symmetrizer-P}, \eqref{eq:case2_def_symmetrizer_R}, and \eqref{eq:case3_def_symmetrizer_R} together with \eqref{eq:yR3y} and \eqref{eq:R3-submat}.
\end{definition}

The symmetrizer $\cR$ is zero homogeneous and  can be extended to a bounded linear operator $\cR:(L^2(V))^n\to(L^2(V))^n$, $n \in \NN$, see e.g. \cite{Grubb1996}, Page 169-171, and using $C_0^\infty$ cut-off functions for the restriction to $V$ \cite{Grubb1996,Wloka1987}.

\begin{definition}
\label{def:pseudo-MN}
    a) The pseudo-differential operator $M\in (\OpsS_{1,0}^1(\overline{V}\times \RR^3))^{4\times 4}$ related to the matrix symbol $\Mhat$ given by \eqref{eq:matrix-M} is defined as
    \begin{align*}
        &M(\bx,t,\imath D_{\bar{x}},\imath D_t+\eta)(e^{-\eta t}\bv(\bx,t))\\
        &\quad =(\cL\cF)_3^{-1}[\Mhat(\bomega,s)\bvh(x_1,\bomega,s)](\bx,t;\eta),
    \end{align*}
    for all $\bv\in (C_{(0)}^\infty(\overline{V}))^4$ with $\bv(\bx,t)=\bm{0}$ for $(\bx,t)\in \RR^4\setminus \overline{V}$.
    
    b) The pseudo-differential operators $N_1,\,N_2\in (\OpsS_{1,0}^0(V\times \RR^3))^{2\times 2}$ are  defined as
    \begin{align*}
        &N_1(\bx,t,\imath D_{\bar{x}},\imath D_t+\eta)(e^{-\eta t}\bv(\bx,t))
         =(\cL\cF)_3^{-1}[\frac{\imath}{s}(\omega_2\EE-\omega_3\JJ)\bvh(x_1,\bomega,s)](\bx,t;\eta), \\
        &N_2(\bx,t,\imath D_{\bar{x}},\imath D_t+\eta)(e^{-\eta t}\bv(\bx,t))
        =(\cL\cF)_3^{-1}[\frac{\imath}{s}(\omega_2\DD-\omega_3\II)\bvh(x_1,\bomega,s)](\bx,t;\eta),
    \end{align*}
    for all $\bv \in (C_{(0)}^\infty(\overline{V}))^2$ with $\bv(\bx,t)=\bm{0}$ for $(\bx,t)\in \RR^4\setminus \overline{V}$.
\end{definition}

\section{Stability bounds for the half-space problems}\label{sec:half-space}
We turn to proving an $L^2$ a-priori bound for the half-space domain $\Omega=\{\bx\in\RR^3:\, x_1<0\}$.
For the stability analysis, we split the domain $\Omega$ into two parts, see Section \ref{ssec:main-results}. On $\Omega\setminus\Omega_\epsilon$, $\epsilon>0$ we use a standard energy estimate, while on $\Omega_\epsilon$ we use the symmetrizer $\cR$ discussed in Section \ref{sec:symmetrizer} to obtain a stability bound. To distinguish the different domains, we introduce the function $\chi_\epsilon$.

Let $\chi_\epsilon\in C_b^\infty(\overline{\Omega})$ be a partition of unity function such that $\chi_\epsilon(\bx)=1$ for $-\frac{\epsilon}{2}<x_1<0$, $\chi_\epsilon(\bx)=0$ for $x_1<-\epsilon$ and $0\leq \chi_\epsilon\leq 1$.
Using the partition function $\chi_\epsilon$ we define
$$
\widetilde{\bm{f}} = (\bm{f}(1-\chi_\epsilon)+A_1\bm{u} \frac{\partial \chi_\epsilon}{\partial x_1})e^{-\eta t},
$$
and state the following lemma.
\begin{lemma}[Energy estimate in half-space]
    \label{lem:stability-half-space}
    Let Assumption~\ref{ass:1} hold.
    Let $\bu\in (H_\eta^1(\Omega\times(0,\infty)))^6$ be the 
    solution to the Maxwell equations \eqref{eq:hyperbolic-sys} on the half-space $\Omega\times(0,\infty)$. 
    Then there exists a constant $C>0$ such that for all $\eta> 0$ it holds that
    \begin{equation}\label{eq:interior_bound_half_space}
    \begin{split}
        \eta  \|(1-\chi_{\epsilon}))\bu\|_{0,\Omega\times[0,\infty)}^2
        \leq C_B \bigg(\frac{1}{\eta} \|\bf\|_{0,\Omega\times[0,\infty)}^2
        +\frac1\eta C_\chi^2\|\bu\|_{0,\Omega\times(0,\infty),\eta}^2
        +\|\bu_0\|_{0,\Omega}^2\bigg).
    \end{split}
    \end{equation}
\end{lemma}
\begin{proof}
Set  $\bv:\Omega\times(0,\infty)\to\RR^6$, with $\bv(\bx,t)=e^{-\eta t}(1-\chi_\epsilon(\bx))\bu(\bx,t)$.
Multiplying \eqref{eq:hyperbolic-sys} by $\e^{-\eta t}(1-\chi_\epsilon)$
yields for $\bv(\bx,t)$
\begin{align*}
 \frac{\partial}{\partial t} B\bv + \eta B\bv
     + \sum_{i=1}^3  A_i \frac{\partial \bv}{\partial x_i}\bv  =\widetilde{\bf}.
\end{align*}
Multiplication of this identity by $\bv$ and integration over $\Omega$ yields for any $t>0$
\begin{align*}
\frac{1}{2}\frac{\mathrm{d}}{\mathrm{d}t} \|B^{1/2}\bv(t)\|_{0,\Omega}^2 +
\eta \|B^{1/2}\bv(t)\|_{0,\Omega}^2 = \left( \tilde \bf(t),\bv(t)\right)_{\Omega},
\end{align*}
where we used integration-by-parts in $x_1$ and that $\bv$ vanishes on $\{x_1=0\}$.
An application of the Cauchy-Schwarz inequality and integration over $t>0$ yields \eqref{eq:interior_bound_half_space},
where we used $\bv(0)=\bu(0)(1-\chi_\epsilon)=\bu_0$, $B$ is uniformly positive and bounded, and $\sup_{\bx\in \Omega}\left|\frac{\partial \chi_\epsilon}{\partial x_1}\right|\leq C_\chi.$
\end{proof}

For the stability analysis on $\Omega_\epsilon$, we 
multiply \eqref{eq:hyperbolic-sys} by $\chi_\epsilon$, and obtain 
\begin{equation}
    \label{eq:chi-hyperbolic}
    \frac{\partial}{\partial t} (\chi_\epsilon \bu) 
     + \sum_{i=1}^3  A_i \frac{\partial}{\partial x_i}(\chi_\epsilon \bu) 
     = - A_1 \bu \frac{\partial\chi_\epsilon}{\partial x_1}=\bh_\chi,\quad\text{in }\Omega\times(0,\infty),
\end{equation}
where we used $B=I_6$ on $\Omega_\epsilon$, and Assumption \ref{ass:1}, which implies that  $\bf\chi_\epsilon=\bm{0}$ on $\Omega\times(0,\infty)$ and $\chi_\epsilon \bu (\cdot,0)=\bm{0}$ on $\Omega$.

After multiplying \eqref{eq:chi-hyperbolic} with $e^{-\eta t},\,\eta>0$ and
applying the bilateral Laplace-Fourier transform, we obtain with $\brho=S^\top \bu$, $\bh=S^\top \bh_\chi$ 
\begin{subequations}
   \begin{alignat*}{3}
    \partial_{x_1}\chi_\epsilon\widehat{\brho}_{\pm} &+ \Mhat(\bomega,s) \chi_\epsilon \widehat{\brho}_{\pm} = \bhh_{\pm}, &\quad& \text{for }x_1<0,\\
    \chi_\epsilon\widehat{\brho}_{\s}&=
    \frac\imath{s}( \omega_2 \EE- \omega_3 \JJ)\chi_\epsilon\widehat{\brho}_+
    +\frac\imath{s}( \omega_2 \DD-\omega_3 \II)\chi_\epsilon\widehat{\brho}_-, &\quad &\text{for }x_1<0,
\end{alignat*} 
\end{subequations}
where $\brho=(\brho_+^\top,\brho_\s^\top,\brho_-^\top)^\top$,
$\brhopm=(\brho_+^\top,\brho_-^\top)^\top$, and $\bh=(\bh_+^\top,\bh_\s^\top,\bh_-^\top)^\top$. The crucial observation is that $\bh_\chi$ is in the range of $A_1$, which implies that $\bh_\s=\bm{0}$.

Finally, applying the inverse bilateral Laplace-Fourier transform and using Definition \ref{def:pseudo-MN} gives the pseudo-differential formulation
\begin{subequations}\label{eq:Maxwell-system-half_space_inhom}
   \begin{alignat}{3}
    \partial_{x_1}(\chi_\epsilon e^{-\eta t}\brhopm) &+ M(\bx,t,\imath D_{\bxb},\imath D_t+\eta) (\chi_\epsilon e^{-\eta t}\brhopm)\notag \\
    &\quad= -\frac{\partial \chi_\epsilon}{\partial x_1}e^{-\eta t}\brhopm, &\quad& \text{on }\Omega\times(0,\infty),\label{eq:Maxwell-form-4_inhom}\\
    \label{eq:rho-o_inhom}
    \chi_\epsilon\brho_{\s}&=
    N_1(\bx,t,\imath D_{\bxb},\imath D_t+\eta)(\chi_\epsilon e^{-\eta t}\brho_+)\notag\\
    &\quad +N_2(\bx,t,\imath D_{\bxb},\imath D_t+\eta)(\chi_\epsilon e^{-\eta t}\brho_-), &\quad &\text{on }\Omega\times(0,\infty),
\end{alignat} 
\end{subequations}
where $M,\, N_1,\,N_2$ are defined in Definition \ref{def:pseudo-MN}.
In the stability analysis, we will use that
\begin{equation}
    \label{eq:rho-psi-j}
    \brho(\bx,t)=\sum_{j=1}^3 \brho_j(\bx,t),
\end{equation}
with $\brho_j=\psi_j(\bx,t)\brho(\bx,t)$ the solution related to each partition $\cU_j$ and $\psi_j$ the partition functions stated in Definition \ref{def:partition-of-psi-j}.

For the stability bound, we define for $(\bomega,s)\in \cU_3$ the variables
\begin{equation}
    \label{eq:vLFinverse}
    \bv = (\cL\cF)_3^{-1}[\Tom^{-1}\chi_\epsilon\widehat{\brho}_{\pm}](\bx,t;\eta),
\end{equation}
$\bv_1=(v_1,v_3)$ and $\bv_2=(v_2,v_4)$ with $v_j$, $j=1,...,4$ the components of $\bv$, and obtain the following bound.

\begin{lemma}[Boundary estimate]
    \label{lem:boundary-est}
    Let Assumption \ref{ass:1} hold.
    Assume that the solution $\bu$ 
    to the Maxwell equations \eqref{eq:hyperbolic-sys} complemented with the HW-NRBCs \eqref{eq:hw-nrbcs} at $\Gamma$ satisfies $\bu \in (H_\eta^1(\Omega\times(0,\infty)))^6$ with trace in $ (H_\eta^1(\Gamma\times(0,\infty)))^6$. Then there exist positive constants $\delta_1,\,\delta_2$, $C_g$, $\nu$ and $\eta_0$, independent of $\bu$, $\eta$, such that for all $\eta\geq \eta_0$, it holds that
    \begin{equation}
        \begin{split}
            &\delta_1 \| \bu\|_{0,\Gamma\times[0,\infty),\eta}^2
            +\eta\delta_2 \|\chi_\epsilon \bu\|_{0,\Omega\times[0,\infty),\eta}^2\\
            &+\delta_1 \|(\imath D_t+\eta)^{-1}( \psi_3\bv_1)\|_{1,\Gamma\times[0,\infty),\eta}^2\\
            &+\eta\delta_2\|(\imath D_t+\eta)^{-1}(\chi_\epsilon \psi_3\bv_1)\|_{1,\Omega\times[0,\infty),\eta}^2 \\
            &\leq \frac{32}{\nu^2}\big(C_g \|\bg_-\|_{1,\Gamma\times[0,\infty),\eta}^2
            +\frac{C_\chi^2}{\eta\delta_2} \|\bu\|_{0,\Omega\times[0,\infty),\eta}^2\big).
        \end{split}
    \end{equation}
\end{lemma}
\begin{proof}
    1) After taking the inner product of \eqref{eq:Maxwell-form-4_inhom} with $\cR_j(\chi_\epsilon e^{-\eta t}\psi_j \brhopm)$, and integrating over $\Omega\times[0,\infty)$ we obtain after integration by parts with respect to $x_1$ using the fact that $\cR_j$ does not depend on $x_1$,
    \begin{align*}
        &\frac{1}{2}\Re{( e^{-\eta t}\psi_j\brhopm, \cR_j (\chi_\epsilon e^{-\eta t}\psi_j\brhopm))_{0,\Gamma\times[0,\infty)}}\\
        &+\Re{(\chi_\epsilon e^{-\eta t}\psi_j\brhopm, \cR_j M (\chi_\epsilon e^{-\eta t}\psi_j\brhopm))_{0,\Omega\times[0,\infty)}}\\
        &=-\Re{(\cR_j(\chi_\epsilon e^{-\eta t}\psi_j\brhopm), \frac{\partial\chi_\epsilon}{\partial x_1} e^{-\eta t}\psi_j\brhopm)_{0,\Omega\times[0,\infty)}},
    \end{align*}
    for $j=1,2,3$, 
    where we used that $\chi_\epsilon=1$ at $\Gamma$.

    2) Next, we obtain bounds on each subdomain $\cU_j$, which has symmetrizer $\cR_j.$
    Using \eqref{eq:total-symmetrizer}, Lemmas \ref{lem:case1_symmetrizer_fulfills_property_ii} and \ref{lem:case2_property_ii_holds} for Property $\mathrm{ii})$, and Lemmas \ref{lem:case1_symmetrizer_fulfills_property_iii} and \ref{lem:case2_symmetrizer_prop2} for Property $\mathrm{iii})$ for the symmetrizer symbols $\cRh_j$, $j=1,2$, together with G{\aa}rding's inequality Lemma \ref{lem:garding-2} with $m=0$, which requires $\eta\geq\eta_0$, we obtain for $j=1,2$ the lower bound
    \begin{align*}
        &\frac{1}{2}\Re{(e^{-\eta t}\psi_j\brhopm, \cR_j ( e^{-\eta  t}\psi_j\brhopm))_{0,\Gamma\times[0,\infty)}}\\
        &+\Re{(\chi_\epsilon e^{-\eta t}\psi_j\brhopm, \cR_j M (\chi_\epsilon e^{-\eta t}\psi_j\brhopm))_{0,\Omega\times[0,\infty)}}\\
        &\geq \frac{1}{4}\delta_1 \| \psi_j\brhopm\|_{0,\Gamma\times[0,\infty),\eta}^2
        +\frac{1}{2}\delta_2 \eta\|\chi_\epsilon \psi_j\brhopm\|_{0,\Omega\times[0,\infty),\eta}^2
        -C_{g,j} \|\bg_-\|_{0,\Gamma\times[0,\infty],\eta}^2.
    \end{align*}
    
    3) Using \eqref{eq:vLFinverse}, Lemma \ref{lem:case3_property_ii_holds}, together with G{\aa}rding's inequality Lemma \ref{lem:garding-3} with $m=2$, which requires $\nu$ to be sufficiently small, for the $\bv_1$-contribution, and Lemma \ref{lem:garding-2} with $m=0$ and $\eta\geq \eta_0$ for the $\bv_2$-contribution gives
    \begin{align*}
        &\frac{1}{2}\Re{( e^{-\eta t}\psi_3\brhopm, \cR_3 ( e^{-\eta  t}\psi_3\brhopm))_{0,\Gamma\times[0,\infty)}}\\
        &\geq \frac{1}{4}\delta_1 \left(
        \|(\imath D_t+\eta)^{-1}( \psi_3\bv_1)\|_{1,\Gamma\times[0,\infty),\eta}^2
        +\|\psi_3\bv_1\|_{0,\Gamma\times[0,\infty),\eta}^2
        \right) \\
        &\quad + \frac{1}{4}\delta_1 
        \| \psi_3\bv_2\|_{0,\Gamma\times[0,\infty),\eta}^2
        -C_{g,3} \|\bg_-\|_{1,\Gamma\times[0,\infty],\eta}^2.
    \end{align*}
    Similarly, using \eqref{eq:rho-psi-j}, Lemma \ref{lem:case3_property_iii_holds}, together with G{\aa}rding's inequality Corollary \ref{cor:garding-4} with $m=2$, and $\eta_1=\frac{1}{2}\eta$ for the $\bv_1$-contribution, and Lemma \ref{lem:garding-2} with $m=0$ and $\eta\geq \eta_0$ for the $\bv_2$-contribution gives
    \begin{align*}
        & \Re{(\chi_\epsilon e^{-\eta t}\psi_3\brhopm, \cR_3 M (\chi_\epsilon e^{-\eta  t}\psi_3\brhopm))_{0,\Omega\times[0,\infty)}}\\
        &\geq \frac{1}{2}\delta_2 \eta \left(
        \|(\imath D_t+\eta)^{-1}(\chi_\epsilon \psi_3\bv_1)\|_{1,\Omega\times[0,\infty),\eta}^2
        +\|\chi_\epsilon \psi_3\bv_1\|_{0,\Omega\times[0,\infty),\eta}^2
        \right) \\
        &\quad + \frac{1}{2}\delta_2 \eta
        \|\chi_\epsilon \psi_3\bv_2\|_{0,\Omega\times[0,\infty),\eta}^2.
    \end{align*}
    
    4) Using $\sup_{\bx\in \Omega}\left|\frac{\partial\chi_\epsilon}{\partial x_1}\right|\leq C_\chi$, the boundedness of $\cR_j$ and the Cauchy-Schwarz inequality gives
    \begin{align*}
        & -\Re{(\cR_j(\chi_\epsilon e^{-\eta t}\psi_j\brhopm),  \frac{\partial\chi_\epsilon}{\partial x_1}e^{-\eta  t} \psi_j\brhopm)_{0,\Omega\times[0,\infty)}}\\
        &\qquad \leq  \frac{1}{4}\delta_2\eta \|\chi_\epsilon \psi_j\brhopm\|_{0,\Omega\times[0,\infty),\eta}^2
        +\frac{C_\chi^2}{\delta_2\eta} \|\brhopm\|_{0,\Omega\times[0,\infty),\eta}^2.
    \end{align*}

    Steps 2-4 then give the bound
    \begin{equation}
    \label{eq:L2-stab-02}
        \begin{split}
            & \frac{1}{4}\delta_1 \sum_{j=1}^3\| \psi_j\brhopm\|_{0,\Gamma\times[0,\infty),\eta}^2
        +\frac{1}{4}\delta_2 \eta \sum_{j=1}^3\|\chi_\epsilon \psi_j\brhopm\|_{0,\Omega\times[0,\infty),\eta}^2\\
        &+\frac{1}{4}\delta_1 
        \|(\imath D_t+\eta)^{-1}( \psi_3\bv_1)\|_{1,\Gamma\times[0,\infty),\eta}^2
        +\frac{1}{2}\delta_2 \eta\|(\imath D_t+\eta)^{-1}(\chi_{\epsilon}\psi_3\bv_1)\|_{1,\Omega\times[0,\infty),\eta}^2
        \\
        &\leq  \sum_{j=1}^3 C_{g,j} \|\bg_-\|_{1,\Gamma\times[0,\infty],\eta}^2
        +\frac{C_\chi^2}{\delta_2\eta} \|\brhopm\|_{0,\Omega\times[0,\infty),\eta}^2.
        \end{split}
    \end{equation}
    
    5a) A bound on the standing waves $\brho_{\s}$ can be obtained from \eqref{eq:rho-o_inhom} for $(\bomega,s)\in \cU_1\cup \cU_2$ using $\left|\frac{\bomega}{s}\right|=\left|\frac{\bomega'}{s'}\right|\leq \frac{2}{\nu}$ and the fact that the norms of the matrices $\EE,\DD,\II$ and $\JJ$ are bounded by $\frac{\sqrt{2}}{2}$,
    \begin{equation}
    \label{eq:L2-stab-03}
        \sum_{j=1}^2 \|\chi_\epsilon \psi_j\brho_{\s}\|_{0,\Omega\times[0,\infty),\eta}^2
        \leq \frac{8}{\nu^2} \sum_{j=1}^2 \|\chi_\epsilon \psi_j\brhopm\|_{0,\Omega\times[0,\infty),\eta}^2.
    \end{equation}

    5b) For $(\bomega,s)\in\cU_3$ we obtain using \eqref{eq:case3_expression_rho_o}
    \begin{equation}
        \label{eq:L2-stab-04}
        \|\chi_\epsilon \psi_3\brho_{\s}\|_{0,\Omega\times[0,\infty),\eta}^2
        \leq  \|(\imath D_t+\eta)^{-1}(\chi_\epsilon \psi_3\bv_1)\|_{1,\Omega\times[0,\infty],\eta}^2.
    \end{equation}

    The right-hand side of \eqref{eq:L2-stab-03} and \eqref{eq:L2-stab-04} can now be bounded by the right-hand side of \eqref{eq:L2-stab-02}, which gives
    \begin{equation*}
        \frac{1}{4}\delta_2\eta\sum_{j=1}^3 \|\chi_\epsilon \psi_j\brho_{\s}\|_{0,\Omega\times[0,\infty),\eta}^2
        \leq 
        \frac{8}{\nu^2}\big(\sum_{j=1}^3 C_{g,j} \|\bg_-\|_{1,\Gamma\times[0,\infty],\eta}^2
        +\frac{C_\chi^2}{\delta_2\eta} \|\brhopm\|_{0,\Omega\times[0,\infty),\eta}^2\big).
    \end{equation*}
    A similar bound for $\frac{1}{4}\delta_1\sum_{j=1}^3 \| \psi_j\brho_{\s}\|_{0,\Gamma\times[0,\infty),\eta}^2$ is obtained on $\Gamma\times[0,\infty)$ since the trace of $\brho_{\s}$ is in $(H_\eta^1(\Gamma\times(0,\infty)))^2$ by assumption and we have $\chi_\epsilon=1$ at $\Gamma$.

    Next, we use the triangle inequality, giving
    \begin{equation}
        \label{eq:chi-rho-tri}
        \|\chi_\epsilon  \brho\| \leq \sum_{j=1}^3 \left( \|\chi_\epsilon \psi_j \brhopm\|+\|\chi_\epsilon \psi_j \brho_{\s}\| \right),
    \end{equation}
    combine all terms related to the contributions of the right-hand side of \eqref{eq:chi-rho-tri} and use $\|\bv_1\|_0^2+\|\bv_2\|_0^2=\|\bv\|_0^2 = \|\brhopm\|_2^2$ since $\Tom$ in \eqref{eq:T-mat-2} is unitary.
    
    Finally, using $\bu=S\brho$, with $S$ being a unitary matrix, which gives $|\brho|=|\bu|$, the claim is proven.
\end{proof}

Combining the results of Lemmas \ref{lem:stability-half-space} and \ref{lem:boundary-est} we can now state the main result.

\begin{theorem}[$L^2$-stability]
\label{thm:L2-stability}
    Let Assumption \ref{ass:1} hold.
    Assume that the solution $\bu$ of the Maxwell equations with HW-NRBCs \eqref{eq:hw-nrbcs} satisfies $\bu\in (H_\eta^1(\Omega\times(0,\infty)))^6$ with trace in $(H_\eta^1(\Gamma\times(0,\infty)))^6$. Then, there are $\eta_0>0$, $\nu>0$ such that for $\eta\geq \eta_0$
    there exists a constant $C>0$, independent of $\bu$ and $\eta$, such that the following  estimate holds
    \begin{equation}
        \begin{split}
            &\| \bu\|_{0,\Gamma\times[0,\infty),\eta}^2
            +\eta \|\bu\|_{0,\Omega\times[0,\infty),\eta}^2\\
            &+\|(\imath D_t+\eta)^{-1}( \psi_3\bv_1)\|_{1,\Gamma\times[0,\infty),\eta}^2
            +\eta\|(\imath D_t+\eta)^{-1}(\chi_\epsilon \psi_3\bv_1)\|_{1,\Omega\times[0,\infty),\eta}^2 \\
            &\leq C \left( \frac{1}{\eta}\|\bf\|_{0,\Omega\times[0,\infty),\eta}^2
            +\|\bu_0\|_{0,\Omega,\eta}^2+\|\bg_-\|_{1,\Gamma\times[0,\infty),\eta}^2 \right).
        \end{split}
    \end{equation}
\end{theorem}
\begin{proof}
    The result follows directly from combining Lemmas \ref{lem:stability-half-space} and \ref{lem:boundary-est}, using  $\|\bu\|_0^2\leq 2(\|\chi_\epsilon\bu\|_0^2+\|(1-\chi_\epsilon)\bu\|_0^2)$, and choosing $\eta$ large enough.
\end{proof}
\begin{remark}
    The parameter $\eta$ in Theorem \ref{thm:L2-stability} must be sufficiently large such that $\Big(\frac{32}{\nu^2\delta_2}+C_B\Big)\frac{C_\chi^2}{\eta}\leq \delta_2\eta$, where $C_\chi$ depends on the thickness of the $\Omega_\epsilon$ boundary. 
    When $\epsilon$ in Assumption \ref{ass:1} approaches zero, 
    $C_\chi \sim {1}/{\epsilon}$, therefore, the stability result in Theorem \ref{thm:L2-stability} holds only for a source and initial data with support in fixed domains that have a positive distance from $\Gamma\times(0,\infty)$ and $\Gamma$, respectively.
\end{remark}

\appendix

\section{HW nonreflecting boundary conditions for the Maxwell equations}
\label{app:deriv-H-W-new}
In \cite{Hagstrom2004} the Hagstrom-Warburton nonreflecting boundary conditions are derived for first-order hyperbolic systems, see also \cite{Hagstrom2020}. In this appendix we provide the details of this approach for the Maxwell equations and also indicate an extra condition which is necessary to ensure stability of the HW-NRBCs. 

We consider a semi-infinite domain $\Omega_n=\{\bx\in \RR^3:\, \bn\cdot\bx<0\}$ with $\bn$ the exterior normal vector at the boundary surface $\Gamma_n.$
The surface \(\Gamma_n\) is extended into a neighborhood \(\Gamma_{n,\epsilon} \subset \mathbb{R}^3\) using a family of surfaces parallel to \(\Gamma_n\),
\[
    \Gamma_{n,\epsilon} = \left\{ \bm{y} \in \mathbb{R}^3:\ \bm{y}=\bm{x}+s \bm{n}, 
    \bm{x}\in \Gamma, -\epsilon \leq s \leq 0 \right\}.
\]
The vector field \(\bm{n}(\bm{x}):\Gamma_{n,\epsilon} \to \mathbb{R}^3\) is defined as the extension of the field of normal vectors to the surface \(\Gamma\).

We introduce the Cartesian \((n,\sigma,\tau)\) coordinate system
associated to \(\Gamma_n\) with \(n\) the coordinate in the direction
of \(\bm{n}\). The \((n,\sigma,\tau)\)-coordinate system is
related to the (\(x_1,x_2,x_3\))-coordinate system using
\[
    (n,\sigma,\tau)^\mathrm{T}=
    T(x_1,x_2,x_3)^\mathrm{T},
\]
with \(T \in \mathbb{R}^{3\times 3}\) a constant coefficient unitary rotation matrix. The gradients are then related as \\
\[
    \nabla _n = T^{-\mathrm{T}} \nabla .
\]
Define 
\begin{align*}
    A_n = \sum_{j=1}^{3} T_{1j}A_j,\
    A_\sigma = \sum_{j=1}^{3} T_{2j}A_j,\
    A_\tau = \sum_{j=1}^{3} T_{3j}A_j.
\end{align*}
The Maxwell equations \eqref{eq:hyperbolic-sys-1} with respect to the \((n,\sigma, \tau)\) coordinate system in \(\Gamma_{n,\epsilon}\) then transform into
\begin{equation}
    \label{eq:maxwell-op}
    L \bm{u}(n,\sigma,\tau,t)= \frac{\partial \bm{u}}{\partial t}+A_n \frac{\partial \bm{u}}{\partial n}
    +A_\sigma \frac{\partial \bm{u}}{\partial \sigma} 
    +A_\tau \frac{\partial \bm{u}}{\partial \tau} =\bm{0},
\end{equation}
where we have assumed that in \(\Gamma_{n,\epsilon}\)
\(\varepsilon_r, \mu_r\) are constant and the right-hand side
\(\bm{f}\) in \eqref{eq:hyperbolic-sys-1} is zero.

We introduce the infinite HW sequence, which is an auxiliary variable
reformulation of the Higdon boundary conditions
\cite{Higdon1986}, see also \cite{Givoli2003}. Let the auxiliary variables \(\bm{\phi}^e_j:\Gamma_{n,\epsilon} \times (0,\infty) \to \mathbb{R}^6\), \(j=1,2, \dots\), satisfy
\begin{subequations}
    \label{eq:vector-HW}
\begin{align}
    \left( a_0 \frac{\partial }{\partial t} + \frac{\partial }{\partial n} \right)\bm{u} &= \frac{\partial \bm{\phi}^e_1}{\partial t} ,\label{eq:vector-HW-1}\\  
    \left( a_j \frac{\partial }{\partial t} + \frac{\partial }{\partial n} \right)\bm{\phi}^e_j &= \left( a_j\frac{\partial }{\partial t} - \frac{\partial }{\partial n} \right)\bm{\phi}^e_{j+1},\quad
    j = 1,2,\ldots \label{eq:vector-HW-2} 
\end{align}
\end{subequations}
with \(\bm{\phi}^e_j(\cdot ,0)=\bm{0}\) in \(\Gamma_{n,\epsilon}\). Since the support of data has a positive distance from \(\Gamma_n\), we also have $\frac{\partial }{\partial t} \bm{\phi}^e_j(\cdot ,0)=\bm{0}$ in $\Gamma_{n,\epsilon}$, see Assumption \ref{ass:1}.

For the derivation of the NRBCs, we need the following property.

\begin{lemma}
    \label{lem:HW-1}
    Given the Maxwell operator \(L\) in \eqref{eq:maxwell-op} on
    the domain \(\Gamma_{n,\epsilon}\times (0,\infty)\), the auxiliary variables
    \(\bm{\phi}^e_j, j=1,2,\cdots,\) satisfy for all \(t> 0\),
    \[
        L \bm{\phi}^e_j(\cdot ,t) \Big|_{\Gamma_n}=\bm{0}\quad \text{at } {\Gamma_n}.
    \]
\end{lemma}
\begin{proof}
    Using \eqref{eq:vector-HW-1} and \eqref{eq:maxwell-op}, there holds
    \[
        \frac{\partial (L \bm{\phi}^e_1)}{\partial t} =L \left( \frac{\partial \bm{\phi}^e_1}{\partial t}   \right)
        = L \left( a_0 \frac{\partial }{\partial t} + \frac{\partial }{\partial n} \right)\bm{u} 
        = \left( a_0 \frac{\partial }{\partial t} + \frac{\partial }{\partial n} \right) L \bm{u}
        =\bm{0} \ \text{in } \Gamma_{n,\epsilon}\times (0,\infty).
    \]
    Since \(\bm{\phi}^e_1(\cdot ,0)=\frac{\partial
    \bm{\phi}^e_1}{\partial t} (\cdot ,0) =\bm{0}\) in \(\Gamma_{n,\epsilon}\), we have 
    \begin{equation}
        \label{eq:vector-HW-3}
        (L\bm{\phi}^e_1)(\cdot ,t)|_{\Gamma_n}=(L\bm{\phi}^e_1)(0,\sigma,\tau ,t)=\bm{0}.
    \end{equation}

    For \(j=1\), we have using \eqref{eq:vector-HW-2}
    \[
        L\left( a_1\frac{\partial }{\partial t} + \frac{\partial }{\partial n}  \right) \bm{\phi}^e_{1} 
        =L\left( a_1\frac{\partial }{\partial t} - \frac{\partial }{\partial n}  \right) \bm{\phi}^e_{2},
        \quad \text{in }\Gamma_{n,\epsilon}\times (0,\infty),
    \]
    which implies that 
    \[
        \left( a_1\frac{\partial }{\partial t} + \frac{\partial }{\partial n}  \right) L\bm{\phi}^e_{1} 
        =\left( a_1\frac{\partial }{\partial t} - \frac{\partial }{\partial n}  \right)L \bm{\phi}^e_{2},
        \quad \text{in }\Gamma_{n,\epsilon}\times (0,\infty).
    \]
    Hence, using \eqref{eq:vector-HW-3} we obtain \(\left( a_1\frac{\partial }{\partial t} - \frac{\partial }{\partial n}  \right)(L\bm{\phi}^e_2)(\cdot ,0) =\bm{0}\) in \(\Gamma_{n,\epsilon}\). 
    Since \(\frac{\partial \bm{\phi}^e_2}{\partial t}(\cdot ,0)=\bm{0}\) in \(\Gamma_{n,\epsilon}\) and \(-\frac{\partial \bm{\phi}^e_2}{\partial n}|_{\Gamma_n}\) is in the negative direction (inward \(\Omega\)) with zero initial condition we have
    \[
        (L\bm{\phi}^e_2)(\cdot ,t)|_{\Gamma_n}=\bm{0} .
    \]

    For \(j> 2\), the result follows directly by induction.
\end{proof}

Define the operators \(\mathcal{L}_j, \mathcal{K}_j: \Gamma_n \times (0,\infty)\to \mathbb{R}^6\) as
\begin{subequations}
\label{eq:op-Lj-Kj}
    \begin{align}
    \mathcal{L}_j &= (a_jA_n -I_6)\frac{\partial }{\partial t} -A_{\sigma} \frac{\partial }{\partial \sigma}
        -A_{\tau} \frac{\partial }{\partial \tau}, \quad j=0,1,2,...\label{eq:op-Lj}\\
        \mathcal{K}_j &= (a_jA_n +I_6)\frac{\partial }{\partial t} +A_{\sigma} \frac{\partial }{\partial \sigma}
        +A_{\tau} \frac{\partial }{\partial \tau},\quad j=1,2,...\label{eq:op-Kj}
\end{align}
\end{subequations}
with parameters \(0<a_0\leq 1,\, 0<a_j< 1,\, j=1,2,\dots,\) to be specified.

Multiplying \eqref{eq:vector-HW-1} with \(A_n\)
and using \(L\bm{u}=\bm{0}\) in \(\Gamma_{n,\epsilon}\times(0,\infty)\) gives
\begin{align}
    a_0 A_n\frac{\partial \bm{u}}{\partial t} + A_n\frac{\partial \bm{u}}{\partial n} &= A_n\frac{\partial \bm{\phi}^e_{1}}{\partial t}, \label{eq:hwnrbc-01} \\
    \frac{\partial \bm{u}}{\partial t} + A_n \frac{\partial \bm{u}}{\partial n}+A_\sigma \frac{\partial \bm{u}}{\partial \sigma}+A_\tau \frac{\partial \bm{u}}{\partial \tau} &= \bm{0},
    \label{eq:hwnrbc-02}
\end{align}
at \({\Gamma_n}\). Note the matrix $A_n$ has rank $4$ hence two conditions on $\frac{\partial \bm{\phi}_1^e}{\partial t}$ in \eqref{eq:hwnrbc-01} are lost, which need to be separately imposed. 
Subtracting \eqref{eq:hwnrbc-02} from \eqref{eq:hwnrbc-01} and restricting \(\bm{\phi}_1^e\) to \({\Gamma_n}\) gives the first HW-NRBC, 
\begin{equation}\label{eq:hw01-2}
    (a_0 A_n-I_6)\frac{\partial \bm{u}}{\partial t} 
    -A_\sigma \frac{\partial \bm{u}}{\partial \sigma}
    -A_\tau \frac{\partial \bm{u}}{\partial \tau} 
    = 
    A_n\frac{\partial \bm{\phi}_{1}}{\partial t} .
\end{equation}

Next, after multiplying \eqref{eq:vector-HW-2} with \(A_n\) and using \(L\bm{\phi}^e_j\Big|_{\Gamma_n}=L\bm{\phi}^e_{j+1}\Big|_{\Gamma_n}=\bm{0}\) by Lemma \ref{lem:HW-1} we obtain
\begin{alignat*}{3}
    &a_j A_n\frac{\partial \bm{\phi}^e_j}{\partial t} + A_n\frac{\partial \bm{\phi}^e_j}{\partial n} = a_jA_n\frac{\partial \bm{\phi}^e_{j+1}}{\partial t} - A_n\frac{\partial \bm{\phi}^e_{j+1}}{\partial n},\quad &&\text{in }\Gamma_{n,\epsilon}\times(0,\infty),\\
    &-\frac{\partial \bm{\phi}^e_j}{\partial t} - A_n \frac{\partial \bm{\phi}^e_j}{\partial n}-A_\sigma \frac{\partial \bm{\phi}^e_j}{\partial \sigma}-A_\tau \frac{\partial \bm{\phi}^e_j}{\partial \tau} \\
    &\qquad\qquad = 
    \frac{\partial \bm{\phi}^e_{j+1}}{\partial t} + A_n \frac{\partial \bm{\phi}^e_{j+1}}{\partial n}+A_\sigma \frac{\partial \bm{\phi}^e_{j+1}}{\partial \sigma}+A_\tau \frac{\partial \bm{\phi}^e_{j+1}}{\partial \tau},
    \quad && \text{in }\Gamma_{n,\epsilon}\times(0,\infty),
\end{alignat*}
for \(j=1,2,...\)
Adding the above two identities and denoting  $\bm{\phi}_j =\bm{\phi}_j^e |_{\Gamma_n}$,  we obtain for the boundary auxiliary variables
\begin{equation}\label{eq:hw02-2}
    (a_j A_n-I_6)\frac{\partial \bm{\phi}_j}{\partial t} 
    -A_\sigma \frac{\partial \bm{\phi}_j}{\partial \sigma}
    -A_\tau \frac{\partial \bm{\phi}_j}{\partial \tau} 
    = 
    (a_j A_n+I_6)\frac{\partial \bm{\phi}_{j+1}}{\partial t} 
    + A_\sigma \frac{\partial \bm{\phi}_{j+1}}{\partial \sigma}+A_\tau \frac{\partial \bm{\phi}_{j+1}}{\partial \tau}.
\end{equation}

We now obtain the HW-NRBCs purely in terms of the auxiliary variables $\bm{\phi}_j$, which depend only on the surface coordinates $\sigma,\tau$ and time $t$.
Using \eqref{eq:op-Lj-Kj}, we can then express \eqref{eq:hw01-2} and \eqref{eq:hw02-2} as the infinite HW-NRBCs for $\bm{\phi}_j$ at $\Gamma_n$  for $t>0$
\begin{subequations}
    \label{eq:app-hw-nrbcs}
    \begin{alignat}{3} 
        A_n\frac{\partial \bm{\phi}_{1}}{\partial t}-\mathcal{L}_0\bm{u}\mid_{\Gamma_n}
        &=\bm{0}, \label{eq:app-rho2-max-3-1}\\
        \mathcal{L}_j\bm{\phi}_j - \mathcal{K}_j\bm{\phi}_{j+1}&=\bm{0}, \quad j=1,2, \dots,\label{eq:app-rho2-max-3-2}
    \end{alignat}
\end{subequations} 
with initial conditions
\begin{equation}
    \label{eq:app--rho2-max-3-4}
    \bm{\phi}_j(\bm{x},0)=\bm{0},\quad \forall \bx\in \Gamma_n, \quad j=1,2,\dots
\end{equation}

In order to truncate the infinite sequence of HW-NRBCs, we introduce characteristic variables
using the matrix decomposition
\begin{equation}
    A_n=S(\bm{n}) \Lambda S(\bm{n})^\top, \quad \text{with } \Lambda=\begin{pmatrix}
        \Lambda_+ & &\\ & \Lambda_0 & \\ & & \Lambda_-
    \end{pmatrix}=\mr{diag}([1,1,0,0,-1,-1]),
\end{equation}
and $S(\bm{n})$ the unitary matrix specified in Appendix \ref{sec:app:mat-decomp}.

The characteristic variables $\brho:\Omega \times (0,\infty)\to \mathbb{R}^6$, $\brho_j:\Gamma_n \times (0,\infty)\to \mathbb{R}^6$, $j=1,2,\dots$, then are defined as
\begin{equation}
    \label{eq:rho-rho-j}
    \brho=S(\bm{n})^{\top}\bm{u},\quad 
    \brho_j=S(\bm{n})^{\top}\bm{\phi}_j,
\end{equation}
where $\brho=\begin{pmatrix}
        \brho_{+}^\mathrm{T},& \brho_{\s}^\mathrm{T}, &
        \brho_{-}^\mathrm{T}
    \end{pmatrix}^\mathrm{T}$, $\brho_j=\begin{pmatrix}
        \brho_{j,+}^\mathrm{T},& \brho_{j,\s}^\mathrm{T}, &
        \brho_{j,-}^\mathrm{T}
\end{pmatrix}^\mathrm{T}$, which decomposes $\brho$ and $\brho_j$ into out- and in-moving components, indicated, respectively, with $+$ and $-$. The standing component is indicated by $\s$.

Using \eqref{eq:rho-rho-j} and after multiplying \eqref{eq:app-hw-nrbcs} with $S(\bm{n})^{\top}$, we obtain at the boundary surface $\Gamma_n$ that
\begin{subequations}
    \label{eq:app-hw-nrbcs-2}
    \begin{alignat}{3} 
        \Lambda \frac{\partial \brho_{1}}{\partial t}-\widetilde{\mathcal{L}}_0\brho
        &=\bg_-,\quad t>0, \label{eq:app-hw-nrbcs-2-1}\\
        \widetilde{\mathcal{L}}_j\brho_j - \widetilde{\mathcal{K}}_j\brho_{j+1}&=\bm{0},\quad~~ t>0,\label{eq:app-hw-nrbcs-2-2}\\
        \brho_j(\bm{x},0)&=\bm{0},  \label{eq:app-hw-nrbcs-2-3}
    \end{alignat}
\end{subequations} 
for $j=1,2,\ldots$, with possibly inhomogeneous data $\bg_-=B_{\Gamma}^-(\bn)S(\bn)^\top\bg$ with $B_{\Gamma}^-(\bn)$ stated in \eqref{eq:B-Gamma}, and
\begin{align*}
    \widetilde{\mathcal{L}}_j &=S(\bm{n})^{\top} \mathcal{L}_j S(\bm{n}) = (a_j \Lambda -I_6)\frac{\partial }{\partial t} 
    -\widetilde{A}_\sigma \frac{\partial }{\partial \sigma}-\widetilde{A}_{\tau} \frac{\partial }{\partial \tau},\\
    \widetilde{\mathcal{K}}_j &=S(\bm{n})^{\top} \mathcal{K}_j S(\bm{n}) = (a_j \Lambda +I_6)\frac{\partial }{\partial t}
    +\widetilde{A}_\sigma \frac{\partial }{\partial \sigma}+\widetilde{A}_{\tau} \frac{\partial }{\partial \tau}.
\end{align*}
To truncate the infinite system of HW-NRBCs \eqref{eq:app-hw-nrbcs-2}, we assume that
$\brho_j=\bm{0}$ for all $j \geq  J+2$, and set the incoming
components of $\brho_{J+1}$ to zero,
\begin{equation}
    \label{eq:TD-0}
    \brho_{J+1,-} = \bm{0} \quad \text{at }\Gamma_n \times (0,\infty).
\end{equation}
Using $\brho_{J+1,-}=\cB^-(\bm{n})\bm{\phi}_{J+1}$, the HW-NRBCs \eqref{eq:app-hw-nrbcs} are then truncated at $J+1$ with the condition
\begin{equation}
    \label{eq:TD-0-transformed}
    \cB^-(\bm{n})\bm{\phi}_{J+1} = \bm{0}.
\end{equation}

The system of boundary conditions \eqref{eq:app-hw-nrbcs-2} is, however, not complete. We need to impose two extra conditions on $\brho_{1,\s}$ since the block $\Lambda_0$ in $\Lambda$ in \eqref{eq:app-hw-nrbcs-2-1} is a $2$-by-$2$ zero block, which follows from the fact that the matrix $A_n$ in \eqref{eq:hwnrbc-01} has rank $4$ instead of $6$.
To form a complete system of boundary conditions, we must keep the equation that lies in the kernel of $A_1$ as part of the HW-NRBCs
\begin{equation}
    \label{eq:TD-1-transformed}
    \cB^0(\bm{n})\cL_1\bm{\phi}_{1} = \bm{0},
\end{equation}
see also Remark \ref{rem:nonuniqueness}. The HW-NRBCs are now given by \eqref{eq:app-hw-nrbcs}-\eqref{eq:app--rho2-max-3-4} and \eqref{eq:TD-0-transformed}-\eqref{eq:TD-1-transformed}.

\section{Matrix decomposition}\label{sec:app:mat-decomp}
The matrix \(A_n=\sum_{i=1}^{3}n_i A_i\) with $\bm{n}=(n_1,n_2,n_3)^\top $ has  decomposition 
\begin{equation}
    A_n =S(\bm{n}) \Lambda S(\bm{n})^\top, \label{eq:An-decomp}
\end{equation}
with
\begin{enumerate}[wide]
    \item[a)] for \(n_1 \neq 0\): 
\begin{equation}
    \label{eq:mat-Sn}
    S(\bm{n}) = \frac{1}{\sqrt{2}\alpha_n}\begin{pmatrix}
        n_1 n_3     & -n_2          & \sqrt{2}\alpha_n n_1  & 0           & -n_1 n_3     & n_2      \\
        n_2 n_3     & n_1           & \sqrt{2}\alpha_nn_2  & 0           & -n_2 n_3     & -n_1       \\
        -\alpha_n^2 & 0             & \sqrt{2}\alpha_nn_3  & 0           & \alpha_n^2   & 0         \\
        -n_2        & -n_1n_3       & 0            & \sqrt{2}\alpha_nn_1 & -n_2         & -n_1n_3   \\
        n_1         & -n_2n_3       & 0            & \sqrt{2}\alpha_nn_2 & n_1          & -n_2n_3   \\
        0           & \alpha_n^2    & 0            & \sqrt{2}\alpha_nn_3 & 0            & \alpha_n^2
    \end{pmatrix} ,
\end{equation}
with \(\alpha_n = \sqrt{n_1^2 + n_2^2}\);
    \item[b)] for \(n_1=0\):
\begin{equation}
\label{eq:mat-Sn0}
    S(\bm{n}) = \frac{1}{\sqrt{2}}\begin{pmatrix}
        0     & -1            & 0           & 0           & 0     &  1        \\
        n_3  & 0             & \sqrt{2}n_2 & 0           & -n_3   & 0         \\
        -n_2   & 0             & \sqrt{2}n_3 & 0           &n_2   & 0           \\
        -1     & 0             & 0           & 0           & -1     & 0         \\
        0     & -   n_3       & 0           & \sqrt{2}n_2 & 0     & -   n_3   \\
        0     & n_2           & 0           & \sqrt{2}n_3 & 0     & n_2        
    \end{pmatrix}.
\end{equation}
Matrix \eqref{eq:mat-Sn0} is obtained by setting $n_1=0$ in \eqref{eq:mat-Sn} and multiplying with $|n_2|/n_2$.
\end{enumerate}

In order to define the boundary conditions, we also introduce the following matrices. 
Let $P^+, P^- \in \mathbb{R}^{2 \times 6}$ be zero matrices except for $P^+(1,1)=P^+(2,2)=1$ and $P^-(5,5)=P^-(6,6)=1$. Then the matrices $\cB^+(\bm{n})$, $\cB^-(\bm{n}) \in \mathbb{R}^{2 \times 6}$ are defined as
\begin{equation}
    \label{eq:B-Gamma}
    \cB^+(\bm{n})=P^+S(\bm{n})^\top\quad \text{and}\quad
    \cB^-(\bm{n})=P^-S(\bm{n})^\top
\end{equation}
and $\cB(\bm{n})=\begin{pmatrix}
    \cB^+(\bm{n}) \\ \cB^-(\bm{n})
\end{pmatrix}\in \mathbb{R}^{4 \times 6}$.

\section{Some technical Lemmas}
\label{sec:app:tech-lemmas}
The following results will be useful in the analysis.
The next lemma is a strengthened version of \cite[Lemma 2]{Kreiss2006}, which states that there is a constant $\delta>0$ such that $\Re(\kappa)\geq \delta\Re(s)$ for $(\bomega,s)\in\cU$.
We show that $\delta=1$, which is sharp, as the choice $\bomega=\bm{0}$ and $s\in\mathbb{R}$ shows.
\begin{lemma}
\label{lem:Rek}
Let $(\bomega,s)\in\cU$.
Then $\kappa=\sqrt{s^2+|\bomega|^2}$ satisfies $\Re(\kappa)\geq\Re(s)$ and $\Im(\kappa)\Im(s)\geq 0$.
\end{lemma}
\begin{proof}
    The statement $\Im(\kappa)\Im(s)\geq 0$ is obvious, and we only show the first statement.
    If $\Im(s)=0$, the result is clear. 
    Suppose $\Im(s)>0$ 
    and write $\kappa=|\kappa| \e^{\imath\psi}$ 
    and $s=|s|\e^{\imath\phi}$ with $\psi,\phi\in(0,\pi/2)$.
    Since $\kappa^2=s^2+\bomega^2$, we have $0\leq\psi\leq\phi$.\\
    \textbf{Case 1. $|\kappa|\geq |s|$.} Using the monotonicity of $\cos(t)$ for $t\in [0,\pi/2]$, we have that
    \begin{align*}
        \Re(\kappa)=|\kappa|\cos(\psi)\geq |s|\cos(\phi)=\Re(s).
    \end{align*}
     \textbf{Case 2. $|\kappa|\leq |s|$.} Using the monotonicity of $\sin(t)$ for $t\in [0,\pi/2]$, we have that
    \begin{align*}
        \Im(\kappa)=|\kappa|\sin(\psi)\leq |s|\sin(\phi)=\Im(s).
    \end{align*}
    Since $\Im(\kappa^2)=\Im(s^2)$, we have $\Re(\kappa)\Im(\kappa)=\Re(s)\Im(s)$ and thus
    \begin{align*}
        \Re(\kappa)=\frac{\Re(s)\Im(s)}{\Im(\kappa)}\geq \Re(s).
    \end{align*}
    The case $\Im(s)<0$ can be treated similarly.
\end{proof}

\begin{lemma}\label{lem:strengthened_inequality}
There exists a continuous bijection $f:[0,1]\to[0,1]$ with $f(0)=1$ and $f(1)=0$ such that
\begin{align*}
    |z-s| \leq f\left(\frac{2\Re(s\overline{z})}{|s|^2+|z|^2}\right) |z+s|,
\end{align*}
for all nonzero $s,z\in\mathbb{C}$ with $\Re(s\overline{z})\geq 0$.
\end{lemma}
\begin{proof}
    We have
    \begin{align*}
        |z+s|^2 &= |z|^2 + |s|^2 + 2 \Re(s\overline{z}),\\
        |z-s|^2 &= |z|^2 + |s|^2 - 2 \Re(s\overline{z}).
    \end{align*}
    By assumption $\Re(s\overline{z})\geq 0$, it therefore follows that $|z-s|\leq |z+s|$.
    Suppose that $\Re(s\overline{z})>0$.
    By squaring the previous identities and rearranging terms,
    we see that $|z-s|\leq f |z+s|$ for some $f\in (0,1)$ if and only if
    \begin{align*}
        |z|^2 + |s|^2 \leq 2\frac{1+f^2}{1-f^2}\Re(s\overline{z}).
    \end{align*}
    By choosing $f$ as the solution to 
    \begin{align*}
        \frac{1-f^2}{1+f^2} = 2 \frac{\Re(s\overline{z})}{|z|^2+|s|^2},
    \end{align*}
    we deduce that the claimed inequality holds.
    Note that $f$ can be continuously extended with $f(0)=1$.
    This concludes the proof.
\end{proof}
\begin{lemma}\label{lem:bound_as-k}
    Let $0<a\leq 1$ and $\bomega\in\mathbb{R}^2$.
    Let $s\in\mathbb{C}$ with $\Re(s)> 0$ and
    set $\kappa=\sqrt{s^2+|\bomega|^2}$.
    Assume that $|s|^2+|\bomega|^2=1$.
    Then it holds
\begin{align*}
    \left|\frac{as-\kappa}{as+\kappa}\right| \leq 
    f\left(\frac{2\Re(as\overline{\kappa})}{|as|^2+|\kappa|^2}\right)\leq f\left(\frac{2a}{1+a^2}(\Re(s))^2\right),
\end{align*}
where $f$ is the function from Lemma~\ref{lem:strengthened_inequality}.
\end{lemma}
\begin{proof}
Using that $|s|\leq 1$ and $|\kappa|\leq 1$ as well as $\Re(\kappa)\geq\Re(s)$ by Lemma~\ref{lem:Rek}, we obtain
\begin{align*}
 \frac{\Re(as\overline{\kappa})}{|as|^2+|\kappa|^2} \geq a \frac{\Re(s)\Re(\kappa)}{a^2+1}\geq  \frac{a}{a^2+1} \Re(s)^2,
\end{align*}
and the result follows from Lemma~\ref{lem:strengthened_inequality}.
\end{proof}
\begin{lemma}\label{lem:lower_bound_omega}
For any $0<\nu<1/(2\sqrt{2})$ the following inequalities hold,
    \begin{alignat}{3}
        |\kappa'| &\geq \frac{1}{2} &\quad&\text{for all } (\bomega,s)\in\cU_3.\label{eq:claim1-B4}\\
        |\bomega'|&\geq\frac{1}{4}&\quad&\text{for all } (\bomega,s)\in\cU_2\cup\cU_3. \label{eq:claim2-B4}\\
        |\eta'|^2 &\leq 2\nu^2 &\quad&\text{for all } (\bomega,s)\in\cU_2,\label{eq:claim3-B4}\\
        \left| |\bomega'|^2-|\xi'|^2\right| &\leq 3\nu^2 &\quad&\text{for all } (\bomega,s)\in\cU_2,\label{eq:claim4-B4}
    \end{alignat}
\end{lemma}
\begin{proof}
\textbf{Case 1.}
Suppose $(\bomega,s)\in\cU_3$, i.e., $|s'|<\nu$.
From $|\bzeta'|^2=1$, we deduce that $|\bomega'|^2\geq 1-\nu^2$, which gives \eqref{eq:claim2-B4} for $(\bomega,s)\in \cU_3$ since $\nu<\frac{1}{2\sqrt{2}}$.
Using $\kappa^2=|\bomega|^2+s^2$, the second triangle inequality implies that
\begin{align*}
    |\kappa'|^2 = ||\bomega'|^2+(s')^2|\geq |\bomega'|^2-|s'|^2\geq 1-2\nu^2,
\end{align*}
i.e., \eqref{eq:claim1-B4} holds, because $\nu<1/(2\sqrt{2})$.

\textbf{Case 2.}
Suppose $(\bomega,s)\in\cU_2$, i.e., $|\kappa'|<\nu$. 
Writing $s=\eta+\imath\xi$, we have  $\kappa^2=|\bomega|^2+\eta^2-\xi^2+2\imath\eta\xi$. Considering the real and imaginary parts of $\kappa$, we obtain the bounds
\begin{align}\label{eq:aux_omega_bound}
    \nu^2> \big| |\eta'|^2+|\bomega'|^2-|\xi'|^2\big|,\quad \nu^2> 2|\xi'\eta'|.
\end{align}
From the first inequality in \eqref{eq:aux_omega_bound}, we deduce that
\begin{align}\label{eq:aux_omega_bound2}
    1 = |\eta'|^2+|\bomega'|^2 + |\xi'|^2= |\eta'|^2+|\bomega'|^2 -|\xi'|^2 + 2|\xi'|^2 <\nu^2 + 2|\xi'|^2,
\end{align}
i.e., $(1-\nu^2)/2 < |\xi'|^2$. From the second inequality in \eqref{eq:aux_omega_bound}, it follows that $|\eta'|^2< \nu^2/(1-\nu^2)$.
Using that $\nu<1/(2\sqrt{2})$ this shows that $|\eta'|^2\leq 2 \nu^2$, i.e., \eqref{eq:claim3-B4} holds.
The triangle inequality and the first inequality in \eqref{eq:aux_omega_bound} then show that \eqref{eq:claim4-B4} holds.
Finally, we prove the lower bound \eqref{eq:claim2-B4} for $(\bomega,s)\in \cU_2$.
The second triangle inequality, \eqref{eq:claim4-B4} and the already shown bound $|\xi'|^2>(1-\nu^2)/2$ imply that
\begin{align*}
    |\bomega'|^2 = |\xi'|^2 + |\bomega'|^2-|\xi'|^2\geq |\xi'|^2- \left| |\bomega'|^2 - |\xi'|^2\right|\geq \frac{1-7\nu^2}{2}\geq \frac{1}{16},
\end{align*}
where we used $\nu<1/(2\sqrt{2})$ in the last step, which gives \eqref{eq:claim2-B4} for $(\bomega,s)\in \cU_2$
\end{proof}

\section{G{\aa}rding's inequalities}
\label{sec:garding}
In this Appendix, we summarize several versions of G{\aa}rding's inequality that are used in Section \ref{sec:half-space} to prove the stability of Maxwell equations with HW-NRBCs. 
Special attention is given to symbols that contain a $1/|s|^2$ singularity.

Recall the pseudo-differential operator $P\in (\OpsS_{1,0}^m(V\times\RR^{n}))^{\ell\times \ell} $ with $m\in \RR$, $\ell,\,n\in \NN$. Then
\[
    P(e^{-\eta t}\bv(\bx,t))=\frac{1}{(2\pi)^{n}}\int_{\RR^{n}}
     e^{\imath \xi t+\imath \bomega\cdot\bxb} \symp (\bx,t,\bomega,\xi) \bvh(x_1,\bomega,\xi)\,\mathrm{d}\bomega\,\mathrm{d}\xi,\quad \forall \bv\in (C_{(0)}^\infty(\overline{V}))^\ell,
\]
where \(\bvh(x_1,\bomega,\xi)=\int_{\RR^{n}}e^{-\imath\xi t-\imath \bomega\cdot \bxb-\eta t}\bv(\bx,t)\,\mathrm{d}\bxb \,\mathrm{d}t\) and $\bv(\bx,t)=\bm{0}$ for $(\bx,t)\in \RR^{n+1}\setminus \overline{V}$. The space-time domain $V$ is given by \eqref{eq:space-V}, see Section \ref{ssec:main-results}, with $V=\RR^{n+1},\RR^{n+1}_+,$ or a smooth bounded domain in $\RR^{n+1}$ with the segment property.

We first state an auxiliary lemma.
\begin{lemma}[\cite{Majda1975}, Lemma 3.1]
\label{lem:commu}
    The operator $(\imath D_t+\eta)^{-1}\in \OpsS_{1,0}^{-1}(\RR)$, defined by
    $$
    (\imath D_t+\eta)^{-1}v(\cdot,t)=\int_{\gamma=0}^t e^{-\eta(t-\gamma)}v(\cdot,\gamma)\,\mathrm{d}\gamma, \quad \eta>0,
    $$
    satisfies
    $$
    \|(\imath D_t+\eta)^{-1}v\|_{0,V,\eta}
    \leq \frac{1}{2\eta}\|v\|_{0,V,\eta}.
    $$
\end{lemma}

Recall $\langle \bv\rangle=(1+|\bv|^2)^{1/2}$.
We have the following G{\aa}rding inequality for the pseudo-differential operator $P$.
\begin{lemma}[G{\aa}rding's inequality]
\label{lem:garding-1}
    Let $K$ be a compact set in $\RR^{n+1}$.
    Given the pseudo-differential operator $P\in (\OpsS_{1,0}^m(K\times\RR^{n}))^{\ell\times \ell}$ with $m\in \RR,$ $\ell,n\in \NN$. Assume that $P$ has a Hermitian symbol $\symp$ and there exist constants $C_0>0$ and $M\geq 0$ such that the symbol $\symp$ satisfies
    \begin{equation*}
        \symp(\bx,t,\bomega,\xi)\geq C_0 \langle (\bomega,\xi) \rangle^m I_\ell\quad \text{for all } |(\bomega,\xi)|\geq M,\ (\bx,t)\in K,
    \end{equation*}
    with $I_\ell$ the $\ell\times\ell$ identity matrix.
    Then there exist constants $C_1>0$ and $C_2$ such that for all $\rmq\in \RR$,
    \begin{equation}
        \Re(P(e^{-\eta t}\bv),e^{-\eta t}\bv)_{K}\geq C_1 \|\bv\|_{m/2,K,\eta}^2
        - C_2 \|\bv\|_{m/2-\rmq,K,\eta}^2,\quad 
        \forall \bv\in (C_{(0)}^\infty(K))^{\ell},
        \label{eq:garding-1}
    \end{equation}
    with norm $\|\bv\|_{m,K,\eta}^2=\|v_1\|_{m,K,\eta}^2+...+\|v_\ell\|_{m,K,\eta}^2$.
\end{lemma}
\begin{proof}
    The proof is a straightforward extension of the proof in \cite{Taylor1974Pseudo}, Page 44, to Hermitian matrices of pseudo-differential operators using the $\|\cdot\|_{m,K,\eta}$ norm on the space-time domain $K$.
\end{proof}

For $\eta$ large enough, we can remove the second term on the right-hand side of \eqref{eq:garding-1}. For $C_2\leq 0$ this is trivial, therefore, we consider $C_2>0$.
Using the pseudo-differential operator $\Lambda_\eta^m \in (\OpsS_{1,0}^m(\RR^{n}))^{\ell\times \ell}$, $m\in \RR$, $\ell\in \NN$ with symbol
\begin{equation*}
    \hat{\lambda}_\eta^m=(1+\eta^2+\xi^2+|\bomega|^2)^{m/2}I_{\ell},
\end{equation*}
we have the following relation
\begin{equation}
    \label{eq:norm-equal}
    \|\bv\|_{m,\eta} = \|\Lambda_\eta^m \bv\|_{0}.
\end{equation}
Using \eqref{eq:norm-equal} and the approach outlined in \cite{BenzoniGavage2006}, Appendix C.2, we obtain

\begin{lemma}
    \label{lem:garding-2}
    Given the open domain $V=\Sigma\times(0,\infty)\subset \RR^{n+1}$ with boundary satisfying the segment condition and $\Sigma\subseteq \RR^n$. Under the conditions of Lemma \ref{lem:garding-1} we have for all $\bv\in (H_\eta^m(V))^\ell,\,m\in\RR,\,\ell\in \NN$ and $\eta\geq \eta_0=\sqrt{\max(0,4C_2^2/C_1^2-1)}$ the G{\aa}rding inequality
    \begin{equation}
        \Re(P(e^{-\eta t}\bv),e^{-\eta t}\bv)_{V}\geq \frac{1}{2}C_1 \|\bv\|_{m/2,V,\eta}^2,
        \label{eq:garding-2}
    \end{equation}
    where $C_1,\,C_2$ are the constants in \eqref{eq:garding-1}.
\end{lemma}
\begin{proof}
    Using \eqref{eq:norm-equal},
    G{\aa}rding's inequality \eqref{eq:garding-2}  follows from \eqref{eq:garding-1} if we can prove that there exists  $\alpha>0$ such that for $\eta $ sufficiently large the following inequality holds
    \begin{equation}
        \label{eq:lambda-eta}
        C_1 \hat{\lambda}_\eta^m- C_2 \hat{\lambda}_\eta^{m-2\rmq} \geq \alpha \hat{\lambda}_\eta^m.
    \end{equation}
   Condition \eqref{eq:lambda-eta} is equivalent to
    $$
    (C_1-\alpha)(1+\eta^2+\xi^2+|\bomega|^2)^{\rmq}\geq C_2\quad \text{for all}\ \xi\in \RR,\, \bomega\in\RR^2,\,\rmq\in\RR.
    $$
    Choosing $\rmq=\frac{1}{2}$  and $\alpha=\frac{1}{2}C_1$ gives the condition $\eta\ge \sqrt{\max(0,4C_2^2/C_1^2-1)}$ and we obtain on the compact set $K\subset \overline{V}$ the G{\aa}rding inequality
    \begin{equation}
        \Re(P(e^{-\eta t}\bv),e^{-\eta t}\bv)_{K}\geq \frac{1}{2}C_1 \|\bv\|_{m/2,K,\eta}^2, \quad \forall \bv\in (C_{(0)}^\infty(K))^\ell.
        \label{eq:garding-1-02}
    \end{equation}
    Next, we  cover $\overline{V}$ with a countable collection of compact sets $K_j$ such that $\cup_jK_j^\circ=V$, with $K_j^\circ$ being the open non-empty interior of $K_j$, see e.g. \cite{Grubb1996} Lemma 2.2. Using a partition of unity subordinate to $K_j,$ we then obtain that \eqref{eq:garding-1-02} applies to any $\bv\in C_{(0)}^\infty(\overline{V})$. The claim then follows from the fact that $C_{(0)}^\infty(\overline{V})$ is dense in $H_\eta^1(V)$ if $V\subseteq \RR^{n+1},\,\RR_+^{n+1},$ or a domain $V\subset \RR^{n+1}$ with boundary satisfying the segment condition, \cite{Grubb1996} Theorem 4.10, or \cite{Wloka1987} Theorem 3.6.
\end{proof}

For the stability analysis using the symmetrizer symbol $\cRh_3$, which has a $1/|s|^2$ singularity, we need the G{\aa}rding inequality for the pseudo-differential operator $Q\in (\OpsS_{1,0}^{-1,m}(V\times\RR^{n}))^{\ell\times \ell},\, m\in \RR, \ell \in \NN$ of the form
\[
    Q(e^{-\eta t}\bv(\bx,t))=\frac{1}{(2\pi)^{n}}\int_{\RR^{n}}
     e^{\imath \xi t+\imath \bomega\cdot\bxb} \symq(\bx,t,\bomega,\xi) \bvh(x_1,\bomega,\xi)\,\mathrm{d}\bomega\,\mathrm{d}\xi 
     \quad \forall\bv \in (C_{(0)}^\infty(\overline{V}))^\ell,
\]
with \(\symq\in (S_{1,0}^{-1,m}(V\times\RR^{n}))^{\ell\times \ell}\)  defined as
$$
\symq(\bx,t,\bomega,\xi)=|s|^2 \symp(\bx,t,\bomega,\xi)
$$
where $\symp\in (S_{1,0}^{m}(V\times\RR^{n}))^{\ell\times \ell}$ and $\bvh(x_1,\bomega,\xi)=\int_{\RR^n}e^{-\imath\xi t-\imath \bomega\cdot\bxb-\eta t}\bv(\bx,t)\mathrm{d}\bxb\mathrm{d}t$ and $\bv(\bx,t)=\bm{0}$ for $(\bx,t)\in \RR^{n+1}\setminus \overline{V}$.

\begin{lemma}
    \label{lem:garding-3}
    Given the space-time domain $V=\Sigma\times (0,\infty)\subset \RR^{n+1}$ with  boundary satisfying the segment condition and $\Sigma\subseteq \RR^n$. Given the pseudo-differential operator $Q\in (\OpsS_{1,0}^{-1,m}(V\times\RR^{n}))^{\ell\times \ell}$, $m\in \RR,\,m\geq 2,\,\ell\in\NN$. Assume that there exist constants $C_0> 0$ and $M\geq 0$ such that
    \begin{equation}
    \label{eq:qxt-Il}
        \symq(\bx,t,\bomega,\xi)\geq C_0 \langle (\bomega,\xi) \rangle^m I_\ell\quad \text{for all}\ |(\bomega,\xi)|\geq M,\ (\bx,t)\in V.
    \end{equation}
    Then for $\eta \geqslant 1,\left|s^{\prime}\right|=|s|/|\bzeta|\leq \nu$, with $ \nu>0$ sufficiently small, there exists a constant $C_1>0$ such that $\forall\, \bv \in\left(H_\eta^m(V)\right)^\ell$
    \begin{equation}
        \label{eq:garding-3}
    \begin{split}
    & \Re {\left(Q(\imath D_t+\eta)^{-1}\left(e^{-\eta t} \bv\right),\left(i D_t+\eta\right)^{-1}\left(e^{-\eta t} \bv\right)\right)_{V}} \\
    & \geq \frac{1}{2} C_1 \left( \|(\imath D_t+\eta)^{-1} \bv\|_{m / 2,V, \eta}^2 +\|\bv\|_{m / 2-1, V, \eta}^2 \right).
    \end{split}
    \end{equation}
\end{lemma}
\begin{proof}

(1) Using Lemma \ref{lem:commu} with $\rmq=1$, which is possible since $\symq$ satisfies \eqref{eq:qxt-Il}, we obtain for a compact set $K\subset \overline{V}$
\begin{equation}
    \label{eq:garding-001}
    \begin{split}
        & \operatorname{Re}\left(Q(\imath D_t+\eta)^{-1}\left(e^{-\eta t} \bv\right),(\imath D_t+\eta)^{-1}\left(e^{-\eta t} \bv\right)\right)_{K}\\
&\geq C_1 \left\|(\imath D_t+\eta)^{-1} \bv\right\|_{m / 2, K,\eta}^2-C_2\|(\imath D_t+\eta)^{-1}\bv\|_{m / 2-1,K, \eta}^2, \ \forall \bv \in (C_{(0)}^\infty(K))^\ell,
    \end{split}
\end{equation}
where $C_1,C_2$ are the constants in \eqref{eq:garding-1}. We assume $C_2>0$, the case $C_2\leq 0$ is trivial.

2) Next, we choose $\alpha>0$ such that for $\eta \geq 1$ the following
inequality holds for $m \geq 2$
\begin{equation}
    \label{eq:garding-lambda-eta}
    \frac{1}{2} C_1 \frac{\hat{\lambda}_\eta^m}{| s|^2}-C_2 \frac{\hat{\lambda}_\eta^{m-2}}{|s|^2} \geqslant \alpha \hat{\lambda}_\eta^{m-2}.
\end{equation}
Since $\eta \geq 1$ we have $|s| \geq 1$, $ \hat{\lambda}_\eta \geq 1$ and \eqref{eq:garding-lambda-eta} implies
$$
\frac{1}{2} C_1 \frac{\hat{\lambda}_\eta^2}{|s|^2}-C_2 \geq \alpha.
$$
Using $\frac{|s|}{\hat{\lambda}_\eta}\leq\frac{|s|}{|\bzeta|}=|s^{\prime}|<\nu$ with $0<\nu<\frac{1}{2\sqrt{2}}$ and choosing $\alpha=\frac{1}{2} C_1$ we obtain the lower bound $C_1 \geq \frac{2\nu^2 }{1-\nu^2}C_2$. Hence, \eqref{eq:garding-lambda-eta} gives for $\nu$ small enough
\begin{equation}
    \label{eq:garding-002}
    \begin{split}
        \frac{1}{2} C_1 \| (\imath D_t+\eta)^{-1} \bv \|_{m / 2, K,\eta}^2
        &-C_2 \| (\imath D_t+\eta)^{-1} \bv \|_{m / 2-1,K, \eta}^2 \\
        &\qquad\geq \frac{1}{2} C_1\|\bv\|_{m / 2-1, K,\eta}^2, \quad \forall \bv \in (C_{(0)}^\infty(K))^\ell.
    \end{split}
\end{equation}
Combining \eqref{eq:garding-001} and \eqref{eq:garding-002} proves the result for $\bv\in \left(C_0^\infty(K)\right)^\ell$, $K$ compact. The remaining steps to obtain the bound on $\bv\in \left(H_\eta^m(V)\right)^\ell$ are identical to those in Lemma \ref{lem:garding-2}.
\end{proof}

If $C_0$ in \eqref{eq:qxt-Il} is written as $C_0=C_0^{\prime}\eta$, with $\eta,C_0'>0$, we obtain the following Corollary.
\begin{corollary}[\cite{Majda1975}, Lemma 3.3]
\label{cor:garding-4}
    Given the space-time domain $V=\Sigma\times(0,\infty)\subset \RR^{n+1}$ with boundary satisfying the segment condition and $\Sigma\subseteq\RR^n$.
    Given the pseudo-differential operator $Q \in (\OpsS_{1,0}^{-1,m}\left(V \times \RR^{d}\right) )^{\ell\times \ell},\, m \in \RR,\, m \geq 2$. Assume that there exist constants $C_0^{\prime}>0$ and $M \geq 0$, such that
    $$
    \begin{array}{r}
    \symq(\bx, t, \bomega, \xi)
    \geq C_0^{\prime} \eta\langle(\bomega, \xi)\rangle^m I_\ell,
    \quad \text{ for all}\ |(\bomega, \xi)| \geq M,\ (\bx,t)\in V.
    \end{array}
    $$    
    Then for $\eta\ge 1$, $\nu$ sufficiently small, and any $\eta_1\in \left(0,\eta\right)$, we have $\forall\bv\in (H_\eta^m(V))^\ell$ the G{\aa}rding inequality 
    $$
    \begin{aligned}
    & \operatorname{Re}\left(Q (i D_t+\eta)^{-1}\left(e^{-\eta t} \bv\right),(\imath D_t+\eta)^{-1}\left(e^{-\eta t} \bv\right)\right)_{V} \\
    & \quad \geq C_1^{\prime}\left(\eta-\eta_1\right)\left(\left\|(\imath D_t+\eta)^{-1} \bv\right\|_{m / 2, V, \eta}^2+\|\bv\|_{m / 2-1, V, \eta}^2\right).
    \end{aligned}
    $$
\end{corollary}
\begin{proof}
    The proof is nearly identical to the proof of Lemma \ref{lem:garding-3}. Only the right-hand side of \eqref{eq:garding-lambda-eta} needs to be replaced with 
    $C_1'(\eta-\eta_1)\hat{\lambda}_\eta^{m-2}$.
\end{proof}

\vspace{10pt}
 \section*{Declarations}

   \subsection*{Data Availability}
   The manuscript has no associated data.

   \subsection*{Competing interests}
    The authors declare that they have no known competing financial interests or personal relationships that could have appeared to influence the work reported in this paper.

\subsection*{Acknowledgements}

The research of K. Liu was partially supported by the China Scholarship Council (No. 201806180078) and the National Natural Science Foundation of China (No. 12401549).
The research of M. Schlottbom and J.J.W. van der Vegt was partially supported by the University of Science and Technology of China, Hefei, Anhui, China.
The research of Y. Xu was partially supported by NSFC grant No. 12471347 and Laoshan Laboratory  (No.LSKJ202300305).

\providecommand{\bysame}{\leavevmode\hbox to3em{\hrulefill}\thinspace}
\providecommand{\MR}{\relax\ifhmode\unskip\space\fi MR }
\providecommand{\MRhref}[2]{%
  \href{http://www.ams.org/mathscinet-getitem?mr=#1}{#2}
}
\providecommand{\href}[2]{#2}

\end{document}